# The Natural Display Topos of Coalgebras

Colin Zwanziger

May 2023

Department of Philosophy

Carnegie Mellon University

Pittsburgh, PA 15213

**Thesis Committee:**

Steve Awodey, Chair

Jeremy Avigad

Jonas Frey

Clive Newstead

*Submitted in partial fulfillment of the requirements for the degree of*

*Doctor of Philosophy in*

*Logic, Computation, and Methodology.*

## Abstract

A classical result of topos theory holds that the category of coalgebras for a Cartesian comonad on a topos is again a topos (Kock and Wraith, 1971).

It is natural to refine this result to a topos-theoretic setting that includes universes. To this end, we introduce the notions of natural display topos and natural Cartesian display comonad, and show that the natural model of coalgebras for a natural Cartesian display comonad on a natural display topos is again a natural display topos. As an application, this result extends the approach to universes of Hofmann and Streicher (1997) from presheaf toposes to sheaf toposes with enough points.

Whereas natural display toposes provide a categorical semantics for a form of extensional Martin-Löf type theory, we also prove our main result in the more general setting of natural typoses, which encompasses models of intensional Martin-Löf type theory.

A natural Cartesian display comonad on a natural typos may also be used as a model for dependent type theory with an S4 box operator, or comonadic modality, as introduced by Nanevski et al. (2008). Modal contexts, which have been regarded as tricky to handle semantically, are interpreted as contexts of the natural typos of coalgebras. We sketch an interpretation within this approach.

As part of the framework in which the above takes place, we introduce a refinement of the notion of natural model (see Awodey, 2018), which is (strictly 2-)equivalent to the notion of full, split comprehension category (see Jacobs, 1993), rather than the notion of category with attributes (Cartmell 1978).

## Acknowledgments

This document would not have taken shape without the input, help, and understanding of many.

Foremost, I would like to thank my advisor, Steve Awodey, for his years of patient tutelage. It was a gift to learn so much from one teacher. I hope that the present work will be some testament to that influence.

Jonas Frey was instrumental in the latter part of my studies, and it was only through his generous and detailed commentary, as well as encouragement, that this project came together.

Jeremy Avigad and Clive Newstead provided valued advice throughout my time at Carnegie Mellon. This work also benefited from many discussions with Mathieu Anel, Andrew Swan, Ulrik Buchholtz, Jon Sterling, Mike Shulman, and other members of and visitors at the homotopy type theory group at CMU.

Peter Dybjer and Thierry Coquand provided valued input during a visit to Gothenburg, Sweden in the fall of 2019.

Special thanks are due to the encouragement, input, and support of Ansten Klev, Mandy Simons, Chris Barker, Justyna Grudzińska, Benedikt Ahrens, Kevin Kelly, Adam Bjorndahl, and Maria Polinsky.

This material is based upon work supported by the US Air Force Office of Scientific Research under award number FA9550-21-1-0009.

Any errors in the document are my own, of course.

# Contents

# Chapter 1

# Introduction

The titular contribution of this work, the construction of the **natural display topos of coalgebras**, along with the related coalgebra constructions that we will develop in parallel, can be motivated by at least two applications. My primary motivation was to give a categorical semantics for S4 dependent type theory (DTT), one of the most fundamental modal DTTs. A second application – which, in the end, is more completely treated within this thesis – is the construction of universes in settings that arise from our coalgebra constructions. This generalizes the construction of the universes in presheaf toposes in Hofmann and Streicher (1997).

#### Outline

We introduce the present work by way of our two main applications in Sections 1.1 and 1.2, respectively. Section 1.3 lists notable contributions that appear in the text.

## 1.1 Semantics of S4 DTT

S4 DTT is a fundamental modal DTT of importance to philosophy and linguistics, computer science, and mathematics (especially via homotopy type theory). In this section, we recall some of the relevance of S4 DTT, and broach our proposal for the categorical semantics of this system. Philosophically, I find S4 DTT and systems of S4 predicate logic inspired by it (like Zwanziger, 2017) to be of particular interest because they belatedly resolve the Quinean criticism of quan-tified modal logic (see Quine, 1980). Our basic semantic proposal will be to interpret S4 DTT into a comonad of **natural models** (see Awodey, 2018); we will make use of a novel variant of the notion of natural model which gives rise to the correct notion of comonad of natural models. Technical exposition of S4 DTT is deferred to Chapter 7, which also contains the sketch of our semantic interpretation.

### 1.1.1 S4 Modality

S4 modal logic arose originally in the propositional case (Lewis and Langford, 1932). S4 propositional logic is perhaps the most extensively studied of all modal logics, and provides a simple and appealing axiomatization of modal necessity. The so-called S4 modal box operator

$$\phi \mapsto \Box\phi \quad ,$$

which we may think of as representing necessity, is axiomatized in the propositional case by adding the principles

(T). $\Box\phi \to \phi$

(4). $\Box\phi \to \Box\Box\phi$

to the so-called normal axioms for a box operator.

In the cases of predicate logic and DTT, the S4 box operator acts on predicate formulas and dependent types, respectively, rather than propositional formulas, but comparable principles obtain.

Though the idea of (S4) predicate modal logic is natural, it has proven thorny to formulate properly. Famously, concern that a modal operator could not be modularly integrated into predicate logic led Quine (see 1980) to reject predicate modal logic as incoherent. While this radical conclusion was unwarranted, the substance of Quine's criticisms cannot be so easily waived away.

Indeed, Barcan (1946)'s original formulation of predicate modal logic omitted standard resources of predicate logic like function symbols. Montague (1973)'s influential typed system was more expressive, allowing a treatment of *de dicto* and *de re* readings, but at the cost of aban-doning the usual rules for $\forall$, $\exists$, $\lambda$, and $=$. The same issue continues to affect modern textbook treatments such as Fitting and Mendelsohn (2012).

It was only with the introduction of S4 DTT by Nanevski et al. (2008) that a modal system appeared that includes a non-modal component with the expressivity of predicate logic and that is governed by the usual rules of inference. The author and collaborators (Awodey et al., 2015; Zwanziger, 2017) subsequently used this breakthrough to extract systems of S4 predicate logic that directly meet Quine's demands.[1] Since these systems are formulated as intensional logics in the sense of Montague (1973), they can also be used to put the widespread use of intensional logic in linguistics on a firmer logical foundation.

As for S4 DTT, its box operator can be viewed alternatively as the analog of the box operator of S4 predicate logic under the propositions-as-types principle or as the dependently-typed generalization of Montague (1973)'s simply-typed operator

$$a \mapsto \langle s, a \rangle \quad ,$$

which maps a type $a$ to the 'type of intensions of type $a$.'

From the perspective of a linguist or philosopher of language, DTT provides a sophisticated tool for modelling contextual and anaphoric dependencies in natural language. In particular, following Sundholm (1989)'s compositional translation of the donkey sentence into DTT, DTT provides the most elegant approach to dynamic semantics. S4 DTT is thus a tool that allows a combined treatment of Montagovian intensional semantics and dynamic semantics. S4 DTT has also been used by the author (Zwanziger, 2019a) to give a computational implementation of Montagovian intensional semantics, using the implementation of S4 DTT as part of the proof assistant Agda by Andrea Vezzosi (The Agda Team, 2005-2022). An approach to hyperintensional semantics using S4 homotopy type theory also appeared at Zwanziger (2018a).

S4 DTT has been applied in homotopy type theory (HoTT), notably as part of cohesive HoTT, (Shulman, 2018), which underpins an ambitious program to give HoTT-based foundations to broader areas of modern geometry, including string theory. Licata et al. (2018) also used S4 DTT to construct universes in models of cubical type theory. The introduction of S4 DTT (Nanevski et al., 2008) was motivated by – and remains influential in – computer science; however, I will leave discussion of computer science applications to the more knowledgeable.

[1]For a philosophical exposition, see Zwanziger (2018b, §3.1).

### 1.1.2 Semantics of S4 Modality

In view of these widespread applications of S4 DTT, the categorical semantics of DTT should be extended to the S4 case.

We take the view that the categorical semantics of an S4 modal box operator should be given by a suitable comonad; the counit of the comonad validates the principle (T), while the comultiplication of the comonad validates the principle (4). For example, the notion of finite meet-preserving comonad on a Heyting algebra provides a serviceable model of S4 propositional logic. In my masters thesis (Zwanziger, 2017), I developed a semantics for S4 predicate logic (and Montagovian intensional logic) in a Cartesian comonad on a topos.

In the categorical semantics of DTT, by a version of the propositions-as-types principle, consideration of subobjects in a category is replaced by consideration of so-called types in an appropriate categorical structure such as a **natural model** (see Awodey, 2018). Indeed, our proposal for the categorical semantics of S4 DTT will be a comonad of natural models, or, more properly, a comonad in the 2-category of natural models.

However, care must be taken in formulating the 2-category of natural models. Indeed, we will use a novel refinement of the notion of natural model in the sense of Awodey (2018), introduced in Chapter 3. This new notion of natural model is most naturally conceived as having equivalent 1-category theory to the original notion, but distinct and more felicitous 2-category theory. Its 2-category theory is equivalent to that of the notion of full, split comprehension category (see Jacobs, 1993), rather than Cartmell (1978)'s notion of category with attributes.[2]

The interpretation of a basic S4 DTT in a comonad

$$\Box : \mathbf{C} \to \mathbf{C}$$

of natural models is sketched in Chapter 7. The **natural model of coalgebras $\mathbf{C}^{\Box}$** constructed in

[2]It has been emphasized by Lumsdaine (2018, 2021) that categories with attributes cannot be 2-categorically identified with full, split comprehension categories.

Chapter 4 is used to interpret the so-called modal contexts of S4 DTT.

Towards interpretations for versions of S4 DTT with more features such as $\Sigma$- and $\Pi$-types, we are led to consider structured natural models that we call **natural typoses** and **natural display toposes**, as well as appropriate comonads thereof (Chapter 3). Such interpretations would accordingly make use of the **natural typos of coalgebras** and **natural display topos of coalgebras**, as constructed in Chapter 5.

## 1.2 Universes

### 1.2.1 Display Toposes

The notion of (elementary) topos (Lawvere, 1970) abstracts to the level of categorical algebra several aspects of the category of sets. However, it is natural to assume the existence in one's category of sets of a Grothendieck universe (or, equivalently, strongly inaccessible cardinal), which is not reflected in the topos axioms.

In remedy of this, various notions of universe in a topos have been introduced. Bénabou (1973)'s axioms are already reasonable, though we will make use of the stronger axioms of Streicher (2005). We will call a topos equipped with an impredicative universe in the sense of Streicher—really a suitable class of display maps (Taylor, 1999) that fixes a universe up to equivalence of internal categories—a **display topos**.

If the notion of display topos is worth its salt, much of the category theory of toposes must be replicable at the level of display toposes.[3] In particular, we would hope that

- the display category of sheaves $\mathcal{E}_J$ for a Cartesian display reflector

  $$J : \mathcal{E} \to \mathcal{E}$$

  on a display topos $\mathcal{E}$ is again a display topos;

[3] In seeking to replicate topos-theoretic results at the level of display toposes, we are inspired by the methodology of Moerdijk and Palmgren (2002).

- the display category of coalgebras $\mathcal{E}^{\Box}$ for a Cartesian display comonad

$$\Box : \mathcal{E} \to \mathcal{E}$$

on a display topos $\mathcal{E}$ is again a display topos.

Indeed, following Streicher (2005, §3), when $J : \mathcal{E} \to \mathcal{E}$ is such a Cartesian display reflector, a universe

$$p : \mho_{\bullet} \to \mho$$

in $\mathcal{E}$ (viewed as a display map) immediately yields the universe

$$a(p) : a(\mho_{\bullet}) \to a(\mho)$$

in $\mathcal{E}_J$, where $a : \mathcal{E} \to \mathcal{E}_J$ is the sheafification display functor, ensuring that $\mathcal{E}_J$ is a display topos.

By contrast, when

$$\Box : \mathcal{E} \to \mathcal{E}$$

is such a Cartesian display comonad and

$$p : \mho_{\bullet} \to \mho$$

is a universe in $\mathcal{E}$, the display map

$$Fp : F\mho_{\bullet} \to F\mho$$

in $\mathcal{E}^{\Box}$, where $F : \mathcal{E} \to \mathcal{E}^{\Box}$ is the cofree display functor, is not a universe for $\mathcal{E}^{\Box}$ in general. Rather, given a universe $\boldsymbol{\mho}$ in $\mathcal{E}$ (viewed as an internal category), a universe in $\mathcal{E}^{\Box}$ is constructed (Corollary 5.2.1) as the object of coalgebras for the canonical induced Cartesian comonad

$$\beta : F\boldsymbol{\mho} \to F\boldsymbol{\mho}$$

in $\mathcal{E}^{\Box}$.

This construction of a universe in $\mathcal{E}^{\Box}$, confirming it as the **display topos of coalgebras**, is an instance of the main construction of the present work (carried out at Theorem 5.1.3).

In this topos-theoretic case, the construction amounts to an untruncated version of the classical construction of subobject classifier in the topos of coalgebras (Kock and Wraith, 1971). In that construction, when

$$\Box : \mathcal{E} \to \mathcal{E}$$

is a Cartesian comonad on a topos $\mathcal{E}$, the subobject classifier in the topos of coalgebras $\mathcal{E}^{\Box}$ is constructed as the internal poset of coalgebras (or, equivalently, the internal poset of fixed points) for the canonical induced Cartesian comonad

$$\beta : F\mathbf{\Omega} \to F\mathbf{\Omega}$$

in $\mathcal{E}^{\Box}$, where $\mathbf{\Omega}$ is a subobject classifier in $\mathcal{E}$ (regarded as an internal poset) and $F : \mathcal{E} \to \mathcal{E}^{\Box}$ is the cofree functor (for details, see, *e.g.*, Mac Lane and Moerdijk 1992, §V.8.).

In addition to advancing the theory of display toposes emerging from Streicher (2005), the construction of the display topos of coalgebras is the first coalgebra construction of its kind. Moerdijk and Palmgren (2002) also raised, but left open, the issue of a coalgebra construction in their setting of **stratified pseudotoposes**; our construction resolves the analogous issue in our comparable setting.

### 1.2.2 Typoses

Whereas display toposes may be used in the categorical semantics of extensional Martin-Löf type theory, one might require a more general categorical notion which encompasses models of intensional Martin-Löf type theory. To this end, we also introduce a notion of **typos**, generalizing the notion of display topos.[4] Naturally, we show that the display category of coalgebras for a Cartesian display comonad on a typos is again a typos, as a variant of our main result.

[4]Taylor (1987) introduced the term ‘typos’ for what we would call a typos such that all morphisms of the underlying category are display maps.

### 1.2.3 Strict Properties

A display topos (or typos) fixes a universe up to equivalence of internal categories. But, mainly when interpreting DTT, one may wish for a universe to satisfy properties which are not equivalence-invariant (see Gratzer et al., 2022).

Most notably, in a model of DTT, when we have two universes $\mathfrak{U}$ and $\mathfrak{U}^+$ and a full and faithful internal functor

$$\mathfrak{U} \xrightarrow{i} \mathfrak{U}^+ \quad ,$$

we would like the type-forming operations at each universe to commute on the nose with the action of $i$. In particular, we would like $i$ to preserve the Cartesian closed structure of $\mathfrak{U}$ on the nose.

Another strict property, which does not require multiple universes to formulate, is **realignment**, which may be formulated thus:

**Definition 1.2.1.** *We say that a universe $\mathfrak{U}$ for a typos $\mathcal{E}$ (See §3.1.1) satisfies the* ***realignment*** *property if, for any monomorphism $m : \Delta \rightarrowtail \Gamma$ of $\mathcal{E}$ and pasting diagram*

$$\begin{array}{ccc} \Delta & \xrightarrow{\quad A \quad} & \mathfrak{U} \\ {\scriptstyle m}\big\downarrow & \phi \cong & \nearrow_{B} \\ \Gamma & & \end{array}$$

*in $\mathbf{Cat}(\mathcal{E})$, there exists some $B' : \Gamma \to \mathfrak{U}$ and $\phi' : B' \cong B$ such that $A = B' \circ m$ and $\phi = \phi' \circ m$.*

A similar abstract formulation of realignment was first extracted by Shulman (2013) from the construction of the universal Kan fibration of Kapulkin, Lumsdaine, and Voevodsky (Kapulkin and Lumsdaine, 2012) and emphasized by Gratzer et al. (2022). Realignment played an important role in the topos-theoretic axioms for models of cubical type theory of Orton and Pitts (2018) and the dissertation Sterling (2021), *inter al.*

Any such strict properties are preserved by our main construction; for example, when $\Box : \mathcal{E} \to \mathcal{E}$ is a Cartesian display comonad on a typos $\mathcal{E}$ and $\mathfrak{U}$ is a universe for $\mathcal{E}$ satisfying realignment, then our universe $(F\mathfrak{U})^\beta$ for $\mathcal{E}^\Box$ also satisfies realignment.

### 1.2.4 Universes of Sheaves

Our main construction, in particular, extends the well-known approach to universes of Hofmann and Streicher (1997) from the setting of presheaf toposes to that of sheaf toposes with enough points (see Chapter 6). A sheaf topos $\mathrm{Sh}(\mathbf{C}, J)$ with enough points is appropriately comonadic over $\mathbf{Set}^X$ for some set $X$, so our construction produces a universe in $\mathrm{Sh}(\mathbf{C}, J)$ from the one in $\mathbf{Set}^X$.

Although the universe in $\mathbf{Set}^{\mathbf{C}^{\mathrm{op}}}$ *à la* Hofmann and Streicher (1997) has the strict properties one might ask for, its aforementioned sheafification *à la* Streicher (2005) lacks them in general. The payoff of our construction here is thus to have a universe in $\mathrm{Sh}(\mathbf{C}, J)$ that satisfies the same strict properties as the usual one in $\mathbf{Set}^{\mathbf{C}^{\mathrm{op}}}$, provided that $\mathrm{Sh}(\mathbf{C}, J)$ has enough points.

Gratzer et al. (2022) do construct such strict universes in arbitrary sheaf toposes. However, our construction is more elementary (avoiding recourse to a small object argument) and constructive.

### 1.2.5 Natural Display Toposes and Natural Typoses

Given our focus on strict properties, it will be natural to employ variations on the notion of display topos and typos that fix a universe only up to isomorphism of internal categories.

We thus define the notions of **natural display topos** and **natural typos**, in which the universe is fixed up to isomorphism by a natural model structure (*cf.* Moggi, 1991). We term this notion of universe a **type classifier**. The strict 2-category of display toposes (resp. typoses) and Cartesian display functors turns out to be 2-equivalent, but not strictly 2-equivalent, to the strict 2-category of natural display toposes (resp. natural typoses) and natural Cartesian display functors (Corollary 3.2.14).

Naturally, in a variant of our main construction, we show that the natural model of coalgebras $\mathcal{E}^{\square}$ for a natural Cartesian display comonad

$$\square : \mathcal{E} \to \mathcal{E}$$

on a natural display topos $\mathcal{E}$ is again a natural display topos (Theorem 5.1.3), and similarly for natural typoses.

Interestingly, it is not the case in general that the natural model of sheaves $\mathcal{E}_J$ for a natural Cartesian display reflector

$$J : \mathcal{E} \to \mathcal{E}$$

on a natural display topos $\mathcal{E}$ is again a natural display topos. Under certain conditions on $J$, one can show that, if it exists, the type classifier in $\mathcal{E}_J$ is the object of sheaves $\mho_j$ for the canonical induced Cartesian reflector

$$j : \mho \to \mho$$

on the type classifier $\mho$ of $\mathcal{E}$. However, $\mho_j$ is, in general, a *stack*, rather than a sheaf. This last fact was already noted for the case of sheaves on a topological space by Grothendieck (1960, §3.3).

Of course, the type classifier in $\mathcal{E}^\Box$ is the object of coalgebras for the canonical induced Cartesian comonad

$$\beta : F\mho \to F\mho$$

in $\mathcal{E}^\Box$, where $F : \mathcal{E} \to \mathcal{E}^\Box$ denotes the cofree natural Cartesian display functor and $\mho$ the type classifier of $\mathcal{E}$. By virtue of its construction internal to $\mathcal{E}^\Box$, this object of coalgebras cannot fail to be a coalgebra.

## 1.3 Summary of Contributions

We note some key points of interest in the text:

- the definitions of **universe** for a display category (Definition 3.1.21), **typos** (Definition 3.1.22) and **display topos** (Definition 3.1.24);
- the definition of **natural model** (in the present sense) (Definition 3.2.2);
- the recovery (Lemma 3.2.8) from first principles of the coherence condition on 2-morphisms of categories with families stipulated by Castellan et al. (2017, Appendix B);
- the 2-equivalence of display toposes and natural display toposes, as well as related state-

ments (Theorem 3.2.13, Corollary 3.2.14), which may be compared to the 2-equivalence of locally Cartesian closed categories and certain categories with families at Clairambault and Dybjer (2011) and Castellan et al. (2017, Appendix B);

- the strict 2-equivalence of full, split comprehension categories and natural models (Theorem 3.4.7);
- the construction of the **natural model of coalgebras**, and related statements (Theorem 4.3.1, Corollary 4.4.1);
- the construction of the **natural display topos of coalgebras**, and related statements (Theorem 5.1.3);
- the extension of Hofmann and Streicher (1997)'s approach to universes in presheaf toposes to sheaf toposes with enough points (Proposition 6.3.4);
- the sketch of the interpretation of S4 DTT into a comonad of natural models (Conjecture 7.2.2).

# Chapter 2

# Preliminaries

## 2.1 Strict 2-Categories

In this section, we introduce some general conventions and propositions, mainly concerning strict 2-categories.

Let $\mathbf{0}$ denote the initial category, $\mathbf{1}$ a terminal category, and $\mathbf{2}$ a category consisting of 2 objects, 0 and 1, and one nontrivial morphism $0 \leqslant 1$.

We will use the term '2-category' for what is sometimes called a 'bicategory,' and 'strict 2-category' for such a 2-category in which all the coherence isomorphisms are identities. Accordingly, we use the term '2-functor' for what is sometimes called a 'pseudofunctor,' *etc.*

Let $\mathbf{C}$ be a strict 2-category. We will refer to 1-morphisms of $\mathbf{C}$ as morphisms and 0-morphisms as objects.

We denote by $\mathbf{Cat}$ the strict 2-category of small categories and by $\mathbf{Cat}(\mathbf{C})$ the (potentially very large) strict 2-category of categories in a category $\mathbf{C}$.

By a **small** ($\mathbf{C}$**-**)**indexed category**, we mean a 2-functor from $\mathbf{C}^{\mathrm{op}}$ to $\mathbf{Cat}$, where $\mathbf{C}$ is a small category, referred to as the **indexing category**. In this case, we denote by $\mathbf{IndCat}(\mathbf{C})$ the strict 2-category of small $\mathbf{C}$-indexed categories. We denote by $\mathbf{IndCat}$ the strict 2-category of small indexed categories over varying indexing categories, as defined, *e.g.*, by Hermida (1993). We say that a small indexed category is **strict** if it is strict as a 2-functor. We denote by $\mathbf{Cat}^{\mathbf{C}^{\mathrm{op}}}$ the usual strict 2-category of small, strict $\mathbf{C}$-indexed categories, and by $\mathbf{IndCat}_{\mathrm{s}}$ the usual strict 2-category of small, strict indexed categories over varying indexing categories.

We will elide the difference between a category $\mathbf{D}$ in a small category $\mathbf{C}$ and its **externalization**, normally written $[\mathbf{D}] \in \mathbf{Cat}^{\mathbf{C}^{\mathrm{op}}}$. We will say that a category $\mathbf{D}$ in a small category $\mathbf{C}$ **internalizes** a small, strict $\mathbf{C}$-indexed category $\mathbf{P}$ if we have $\mathbf{D} \cong \mathbf{P} \in \mathbf{Cat}^{\mathbf{C}^{\mathrm{op}}}$.

Following Johnstone (2002), we use the term **Cartesian** for categories with finite limits and functors that preserve finite limits. We denote by $\mathbf{Cart}$ the strict 2-category of small, Cartesian categories and Cartesian functors.

Let $\mathbf{C}$ be a small, Cartesian category. We will say that a small $\mathbf{C}$-indexed category is **Cartesian** if it factors through $\mathbf{Cart}$, and similarly for strict $\mathbf{C}$-indexed functors between small, Cartesian $\mathbf{C}$-indexed categories. We will denote by $\mathbf{Cart}^{\mathbf{C}^{\mathrm{op}}}$ the strict 2-category of Cartesian, small, strict $\mathbf{C}$-indexed categories and Cartesian, strict $\mathbf{C}$-indexed functors. Finally, we will say that a category in $\mathbf{C}$ is **Cartesian** if it is Cartesian as a strict indexed category, and similarly for functors in $\mathbf{C}$.

We use the term **strict** 2**-Cartesian** for strict 2-categories that admit strict 2-limits that are finite in the sense of Street (1976) and strict 2-functors that preserve the same (up to coherent isomorphism).

**Proposition 2.1.1.** *Let* $\mathbf{C}$ *be a Cartesian category. Then, the strict* 2*-category* $\mathbf{Cat}(\mathbf{C})$ *is strict* 2*-Cartesian. Moreover, let* $\mathbf{D}$ *be a Cartesian category and* $F : \mathbf{C} \to \mathbf{D}$ *a Cartesian functor. Then, the strict* 2*-functor* $\mathbf{Cat}(F) : \mathbf{Cat}(\mathbf{C}) \to \mathbf{Cat}(\mathbf{D})$ *is strict* 2*-Cartesian.*

*Proof.* For the first claim, see Street (1976). As for the second, since the construction of finite strict 2-limits in $\mathbf{Cat}(\mathbf{C})$ proceeds as a Cartesian construction in $\mathbf{C}$, it is preserved by the Cartesian functor $F$. □

## 2.2 Gray-Categories

We will have occasion to work categorically with strict 2-categories, in particular in Chapter 4, when we show that the strict 2-category of categories with attributes admits the construction of coalgebras. For this purpose, we make use of **Gray-category** theory (Gray, 1974; Gordon et al., 1995), a form of semi-strict 3-category theory.

Just as the algebra of sets and functions provides a motivation for the axioms of a category, the algebra of strict 2-categories, strict 2-functors, (potentially nonstrict) 2-natural transformations, and modifications provides a motivation for the axioms of a Gray-category.

Informally, a Gray-category is a 3-category in which composition of 1-morphisms is strictly associative and unital, but horizontal composition of 2-morphisms is defined only up to coherent isomorphism. We refer the reader to Gordon et al. (1995) for a formal introduction.

We also make minor use of the notion of **lax Gray-category**. These are more general than Gray-categories, with motivation provided by the algebra of strict 2-categories, strict 2-functors, lax 2-natural transformations, and modifications. In a lax Gray-category, composition of 1-morphisms remains strictly associative and unital, but horizontal composition of 2-morphisms is defined only up to a coherence 3-morphism in the appropriate 'lax' direction. Lax Gray-categories, rather than Gray-categories, were the earlier notion, introduced by Gray (1974, §1,4 Appendix C). For a formal introduction to lax Gray-categories, we refer the reader to this original, though noting that the term 'lax Gray-category' came into use much later (*e.g.* MacDonald and Scull, 2021; Morehouse, 2022).

Unfortunately, to make our theory applicable to large strict 2-categories, we are forced to entertain very large collections:

We denote by $\mathbf{GRAY}_\ell$ (resp. $\mathbf{GRAY}$) the (very large) lax Gray-category (resp. Gray-category) of strict 2-categories, strict 2-functors, lax 2-natural transformations (resp. 2-natural transformations), and modifications. By contrast, we denote by $\mathbf{2CAT}$ the (very large) strict 3-category of strict 2-categories, strict 2-functors, strict 2-natural transformations, and modifications.

**Definition 2.2.1.** *Let* $\mathbf{C}$ *be a lax Gray-category. We will call a (*$\{0\} \hookrightarrow \mathbf{2} \hookleftarrow \{1\}$*)-weighted lax Gray-limit of a cospan* $f : A \to C \leftarrow B : g$ *in* $\mathbf{C}$ *a* ***Gray comma object*** *and denote it by* $f \downarrow g$.

It was observed by Shulman (2012) that a Gray comma object in $\mathbf{GRAY}_\ell$ is what Gray (1969) called a '2-comma category.'

We will denote by $\mathbf{I}$ the 'walking equivalence' strict 2-category, with objects denoted by $0$ and $1$.

**Definition 2.2.2.** *Let* $\mathbf{C}$ *be a Gray-category. We will call a (*$\{0\} \hookrightarrow \mathbf{I} \hookleftarrow \{1\}$*)-weighted Gray-* *limit of a cospan $f : A \to C \leftarrow B : g$ in $\mathbf{C}$ an **equi-comma object** and denote it by $f \wr g$.*

An equi-comma object in $\mathbf{GRAY}$ is what Buckley (2014) called an 'equiv-comma.'

When $\lozenge : \mathbf{C} \to \mathbf{C}$ is a Gray-monad on a Gray-category $\mathbf{C}$, we will denote by $\mathbf{C}^{\lozenge}$ the Gray-category of algebras for $\lozenge$ in the sense of Power (2007).

**Definition 2.2.3.** *We will say that a Gray-functor $F : \mathbf{D} \to \mathbf{C}$ is Gray-**monadic** if it is equivalent (in the strict 2-category of Gray-functors with codomain $\mathbf{C}$) to the forgetful Gray-functor $U : \mathbf{C}^{\lozenge} \to \mathbf{C}$ for some Gray-monad $\lozenge : \mathbf{C} \to \mathbf{C}$.*

**Proposition 2.2.4.** *Gray-monadic Gray-functors create any equi-comma objects that exist in their codomain.*

*Proof.* Similar to Example 4.1 of Power (2007). □

**Lemma 2.2.5.** *If a morphism of Gray-functors of form*

$$\begin{array}{ccc} \mathbf{D} & \xrightarrow{H} & \mathbf{C}^{\to} \\ {\scriptstyle F}\downarrow & & \downarrow{\scriptstyle \mathrm{cod}_{\mathbf{C}}} \\ \mathbf{E} & \xrightarrow[I]{} & \mathbf{C} \end{array}$$

*is equivalent (in the strict 2-category of morphisms of Gray-functors with codomain $\mathrm{cod}_{\mathbf{C}}$) to the forgetful morphism of Gray-functors*

$$\begin{array}{ccc} (\mathbf{C}^{\to})^{(\lozenge^{\to})} & \longrightarrow & \mathbf{C}^{\to} \\ \downarrow & & \downarrow{\scriptstyle \mathrm{cod}_{\mathbf{C}}} \\ \mathbf{C}^{\lozenge} & \longrightarrow & \mathbf{C} \end{array}$$

*induced by some Gray-monad $\lozenge : \mathbf{C} \to \mathbf{C}$, then it creates any fiberwise equi-comma objects that exist in its codomain. That is, for any $\Delta \in \mathbf{E}$, the Gray-functor $F^{-1}(\{\Delta\}) \to \mathbf{C}/I\Delta$ obtained by restricting $H$ to the fiber of $F$ over $\Delta$ creates any equi-comma objects that exist in its codomain.*

*Proof.* Comparable to Proposition 2.2.4. □

## 2.3 Comonads

The notion of **comonad** in a strict 2-category is central to the present work. We will use the theory of such comonads that may be extracted (by dualization) from "The Formal Theory of Monads" (Street, 1972) and "The Formal Theory of Monads II" (Lack and Street, 2002), as well as some minimal extensions thereof.

### 2.3.1 The Formal Theory of Comonads

In this section, we collect relevant information pertaining to the application of the formal theory of comonads as described by Street (1972).

We will frequently write $\Box : C \to C$ for a comonad on an object $C$ in a strict 2-category $\mathbf{C}$ with underlying morphism $\Box : C \to C$, counit $\varepsilon : \Box \Rightarrow \mathrm{id}_C$ and comultiplication $\delta : \Box \Rightarrow \Box\Box$. We will frequently write $L \dashv R : C \to D$ for an adjunction in a strict 2-category $\mathbf{C}$ with right and left adjoints $R : C \leftrightarrows D : L$, unit $\eta : \mathrm{id}_D \Rightarrow RL$, and counit $\varepsilon : LR \Rightarrow \mathrm{id}_C$.

We recall the following basic fact.

**Proposition 2.3.1.** *Let* $\mathbf{C}$ *be a strict* 2*-category. Any adjunction* $L \dashv R : C \to D$ *in* $\mathbf{C}$ *induces a comonad* $\Box : C \to C$ *by setting*

- $\Box :\equiv LR$;
- $\varepsilon :\equiv \varepsilon$;
- $\delta :\equiv L\eta R$.

Any adjunction $L \dashv R : C \to D$ that induces (by Proposition 2.3.1) the comonad $\Box : C \to C$ is said to be a **decomposition** for $\Box$.

When $\mathbf{C}$ is a strict 2-category and $\Box_1 : C_1 \to C_1$ and $\Box_2 : C_2 \to C_2$ are comonads in $\mathbf{C}$, a **morphism of comonads** in $\mathbf{C}$, $f : \Box_1 \to \Box_2$, consists of a morphism $f : C_1 \to C_2$ between the underlying objects of $\mathbf{C}$, together with a 2-morphism $\tau : f\Box_1 \Rightarrow \Box_2 f$ satisfying the commutative diagrams

$$\begin{array}{ccc} f\Box_1 & \xrightarrow{\tau} & \Box_2 f \\ & {\scriptstyle f\varepsilon_1}\searrow \quad \swarrow{\scriptstyle \varepsilon_2 f} & \\ & f & \end{array}$$

and

$$\begin{array}{ccccc} f\Box_1 & & \xrightarrow{\tau} & & \Box_2 f \\ \downarrow{\scriptstyle f\delta_1} & & & & \downarrow{\scriptstyle \delta_2 f} \\ f\Box_1\Box_1 & \xrightarrow[\tau\Box_1]{} & \Box_2 f\Box_1 & \xrightarrow[\Box_2\tau]{} & \Box_2\Box_2 f \end{array} \quad .$$

When $\mathbf{C}$ is a strict 2-category, $\Box_1 : C_1 \to C_1$ and $\Box_2 : C_2 \to C_2$ are comonads in $\mathbf{C}$, and $f_1 : \Box_1 \to \Box_2$ and $f_2 : \Box_1 \to \Box_2$ are morphisms of comonads in $\mathbf{C}$, a **2-morphism of comonads** in $\mathbf{C}$, $\alpha : f_1 \Rightarrow f_2$, consists of a 2-morphism, also denoted $\alpha : f_1 \Rightarrow f_2$, between the underlying morphisms of $\mathbf{C}$, satisfying the commutative diagram

$$\begin{array}{ccc} f_1\Box_1 & \xrightarrow{\alpha\Box_1} & f_2\Box_2 \\ \downarrow{\scriptstyle \tau_1} & & \downarrow{\scriptstyle \tau_2} \\ \Box_2 f_1 & \xrightarrow[\Box_2\alpha]{} & \Box_2 f_2 \end{array} \quad .$$

We then write $\mathbf{Comon}(\mathbf{C})$ for the strict 2-category of comonads and morphisms and 2-morphisms of comonads in $\mathbf{C}$.

When $\mathbf{C}$ is a strict 2-category and $\Box : C \to C$ is a comonad in $\mathbf{C}$, an **object of coalgebras** for $\Box$ consists of an object $C^\Box$, such that

$$\mathrm{Hom}_{\mathbf{C}}(D, C^\Box) \cong \mathrm{Hom}_{\mathbf{Comon}(\mathbf{C})}(\mathrm{id}_D, \Box) \quad , \tag{2.1}$$

naturally in $D$.[1] We will say that a strict 2-category $\mathbf{C}$ **admits the construction of coalgebras** if an object of coalgebras exists for each comonad in $\mathbf{C}$.

When $C^\Box$ is an object of coalgebras for a comonad $\Box : C \to C$ in a strict 2-category $\mathbf{C}$, the

[1]We follow Street (1972) in using the term 'object of coalgebras.' However, the less suggestive term '(co)Eilenberg-Moore object' is more common.

**cofree morphism**

$$F : C \to C^{\Box}$$

is defined as the transpose across the isomorphism 2.1 of the morphism of comonads

$$(\Box, \delta) : (C, \mathrm{id}) \to (C, \Box) \quad ,$$

while the **forgetful morphism**

$$U : C^{\Box} \to C$$

is defined as the underlying morphism of the transpose across the isomorphism 2.1 of the morphism

$$\mathrm{id}_{C^{\Box}} : C^{\Box} \to C^{\Box} \quad .$$

As developed in Street (1972), we then have

$$U \dashv F : C \to C^{\Box} \quad ,$$

and this adjunction is a decomposition for

$$\Box : C \to C \quad ,$$

which we term the **forgetful-cofree decomposition**.

When

$$L \dashv R : C \to D$$

is an adjunction decomposing the comonad

$$\Box : C \to C \quad ,$$

the **comparison morphism**

$$K : D \to C^{\Box}$$

is defined as the transpose across the isomorphism 2.1 of the morphism of comonads

$$(L, L\eta) : (C, \mathrm{id}_C) \to (C, \Box) \quad .$$

This adjunction is **comonadic** (resp. **strictly comonadic**, **comonadic on the nose**) if the comparison morphism is an equivalence (resp. isomorphism, identity). Accordingly, any morphism $L : D \to C$ in a strict 2-category $\mathbf{C}$ is **comonadic** (resp. **strictly comonadic**) if it is the left adjoint of some comonadic (resp. strictly comonadic) adjunction.

When $\mathbf{C}$ and $\mathbf{D}$ are strict 2-categories that admit the construction of coalgebras, we say that a strict 2-functor $\Phi : \mathbf{C} \to \mathbf{D}$ **preserves** (resp. **strictly preserves**, **preserves on the nose**) **the construction of coalgebras** if, for any comonad $\Box : C \to C$ in $\mathbf{C}$, the adjunction $\Phi U \dashv \Phi F : \Phi C \to \Phi(C^{\Box})$ in $\mathbf{D}$ is comonadic (resp. strictly comonadic, comonadic on the nose).

When $\mathbf{C}$ and $\mathbf{D}$ are strict 2-categories that admit the construction of coalgebras and $F : \mathbf{C} \to \mathbf{D}$ is a (contextually salient) strict 2-functor, we may say that $\mathbf{C}$ **admits the construction of coalgebras over** $\mathbf{D}$ if $F$ preserves the construction of coalgebras on the nose.

Moreover, when $F : \mathbf{C} \to \mathbf{E}$ and $G : \mathbf{D} \to \mathbf{E}$ are (contextually salient) strict 2-functors respectively exhibiting $\mathbf{C}$ and $\mathbf{D}$ as admitting the construction of coalgebras over $\mathbf{E}$, we may say that a strict 2-functor $H : \mathbf{C} \to \mathbf{D}$ **preserves** (resp. **strictly preserves**, **preserves on the nose**) **the construction of coalgebras over** $\mathbf{E}$ if $F = GH$ and $H$ preserves (resp. strictly preserves, preserves on the nose) the construction of coalgebras.

### Internal Categories

**Proposition 2.3.2.** *Let* $\mathbf{C}$ *be a Cartesian category. Then, the strict* 2*-category* $\mathbf{Cat}(\mathbf{C})$ *admits the construction of coalgebras. Moreover, let* $\mathbf{D}$ *be a Cartesian category and* $F : \mathbf{C} \to \mathbf{D}$ *a Cartesian functor. Then, the strict* 2*-functor* $\mathbf{Cat}(F) : \mathbf{Cat}(\mathbf{C}) \to \mathbf{Cat}(\mathbf{D})$ *strictly preserves the construction of coalgebras.*

*Proof.* Objects of coalgebras can be recovered as finite strict 2-limits (*cf.* Street, 1976). We can thus apply Proposition 2.1.1. □

**Corollary 2.3.3.** *Let* $\mathbf{C}$ *be a small, Cartesian category. Then, the externalization strict* 2*-functors*

$\mathbf{Cat}(\mathbf{C}) \hookrightarrow \mathbf{Cat}^{\mathbf{C}^{\mathrm{op}}}$ *and* $\mathbf{Cart}(\mathbf{C}) \hookrightarrow \mathbf{Cart}^{\mathbf{C}^{\mathrm{op}}}$ *strictly preserve the construction of coalgebras.*

When $\mathbf{C}$ is a Cartesian category, we say that a comonad $\square : \mathbf{D} \to \mathbf{D}$ in $\mathbf{Cat}(\mathbf{C})$ is **Cartesian** if its underlying internal functor is Cartesian.

### 2.3.2 The Formal Theory of Comonads II

In this section, we collect information pertaining to the application of the extension of the formal theory of comonads described by Lack and Street (2002).

We will denote by $\mathbf{GRAY}_\square$ the (very large) Gray-category of strict 2-categories that admit the construction of coalgebras and strict 2-functors that preserve the construction of coalgebras.

**Proposition 2.3.4.** *The Gray-category* $\mathbf{GRAY}_\square$ *is Gray-monadic over* $\mathbf{GRAY}$.

*Proof.* The strict 3-monad $\mathrm{EM} : \mathbf{2CAT} \to \mathbf{2CAT}$ whose algebras are essentially strict 2-categories equipped with a choice of an object of coalgebras for any comonad, which was described (via the dual case of monads) by Lack and Street (2002), may naturally be weakened to a Gray-monad, which we denote by $\mathrm{EM} : \mathbf{GRAY} \to \mathbf{GRAY}$. □

**Corollary 2.3.5.** *The Gray-category* $\mathbf{GRAY}_\square$ *is closed under the construction of equi-comma objects in* $\mathbf{GRAY}$.

*Proof.* From Proposition 2.2.4. □

We will denote by $\mathbf{GRAY}_\square^{\to}$ the (very large) Gray-category of strict 2-functors that preserve on the nose the construction of coalgebras and morphisms of strict 2-functors such that both components preserve the construction of coalgebras, and by $\mathrm{cod}_\square : \mathbf{GRAY}_\square^{\to} \to \mathbf{GRAY}_\square$ the evident codomain Gray-functor.

**Proposition 2.3.6.** *The forgetful morphism of Gray-functors*

$$\begin{array}{ccc} (\mathbf{GRAY}^{\to})^{(\mathrm{EM}^{\to})} & \longrightarrow & \mathbf{GRAY}^{\to} \\ \big\downarrow & & \big\downarrow \mathrm{cod} \\ \mathbf{GRAY}^{\mathrm{EM}} & \longrightarrow & \mathbf{GRAY} \end{array}$$

*induced by the Gray-monad* $\mathrm{EM} : \mathbf{GRAY} \to \mathbf{GRAY}$ *is equivalent (in the strict* $2$*-category of morphisms of Gray-functors with codomain* $\mathrm{cod} : \mathbf{GRAY}^{\to} \to \mathbf{GRAY}$*) to the forgetful morphism of Gray-functors*

$$\begin{array}{ccc} \mathbf{GRAY}_{\square}^{\to} & \longrightarrow & \mathbf{GRAY}^{\to} \\ {\scriptstyle \mathrm{cod}_{\square}}\downarrow & & \downarrow{\scriptstyle \mathrm{cod}} \\ \mathbf{GRAY}_{\square} & \longrightarrow & \mathbf{GRAY} \end{array}.$$

*Proof.* The equivalence of the forgetful Gray-functors $\mathbf{GRAY}^{\mathrm{EM}} \to \mathbf{GRAY}$ and $\mathbf{GRAY}_{\square} \to \mathbf{GRAY}$ is Proposition 2.3.4. The equivalence of the forgetful Gray-functors $(\mathbf{GRAY}^{\to})^{(\mathrm{EM}^{\to})} \to \mathbf{GRAY}^{\to}$ and $\mathbf{GRAY}_{\square}^{\to} \to \mathbf{GRAY}^{\to}$ is based on the equivalence of algebras for $\mathrm{EM}^{\to} : \mathbf{GRAY}^{\to} \to \mathbf{GRAY}^{\to}$ and (strict) morphisms of algebras for $\mathrm{EM} : \mathbf{GRAY} \to \mathbf{GRAY}$. □

**Corollary 2.3.7.** *Let* $\mathbf{C}$ *be a strict* $2$*-category that admits the construction of coalgebras. The (very large) Gray-category of strict* $2$*-categories that admit the construction of coalgebras over* $\mathbf{C}$ *and strict* $2$*-functors that preserve the construction of coalgebras over* $\mathbf{C}$ *is closed under the construction of equi-comma objects in* $\mathbf{GRAY}/\mathbf{C}$.

*Proof.* From Lemma 2.2.5. □

# Chapter 3

# Natural Display Toposes

We introduce the notion of **natural display topos**, a notion of ‘topos with a universe’ derived from Streicher (2005), but in which the universe is fixed up to isomorphism, as suggested by Moggi (1991). We also introduce the more general notion of **natural typos**, which refines a notion of Taylor (1987), to encompass natural models of intensional Martin-Löf type theory. Our approach is based on a refinement of the theory of natural models (for which, see Awodey, 2018).

#### Outline

Section 3.1 gives an introduction to **display categories**, **typoses**, and **display toposes**. Section 3.2 gives our introduction to natural models, natural typoses, and natural display toposes. Section 3.3 contextualizes natural models within the literature. Section 3.4 gives an introduction to full, split comprehension categories and proves the equivalence of these with natural models in the present sense.

## 3.1 Display Categories, Typoses, and Display Toposes

> We take the position that size must be introduced by the endowment of extra data which is not universally determined.... (Street, 1980)

In this section, we discuss some theory of **display categories** (Taylor, 1999), including the theory of **typoses** and **display toposes**.

We also lay the groundwork for the treatment of **natural models** in Section 3.2. That part of the development is indebted to Shulman (2019).

### 3.1.1 Objects

We define some basic objects of display category theory, including our notions of typos and display topos.

#### Basics

In brief, a display category $\mathbf{C}$ consists of a category, which we also denote by $\mathbf{C}$, endowed with the extra data of a class of morphisms of $\mathbf{C}$, called **display maps**, closed under pullback and composition with isomorphisms.

**Definition 3.1.1** (Taylor 1999)**.** *A **display category** $\mathbf{C}$ consists of*

- *a category, also denoted $\mathbf{C}$;*
- *a class of morphisms of $\mathbf{C}$, denoted $\mathcal{D}_{\mathbf{C}}$, such that*
  - *for every $p : E \to B$ with $p \in \mathcal{D}_{\mathbf{C}}$ and $f : A \to B$, there exists a pullback $f^*p : f^*E \to A$ with $f^*p \in \mathcal{D}_{\mathbf{C}}$;*
  - *$\mathcal{D}_{\mathbf{C}}$ is closed under composition with isomorphisms.*

*We say that a morphism $p \in \mathcal{D}_{\mathbf{C}}$ is a **display map** (or, simply, **display**).*

We will henceforth assume that a display category $\mathbf{C}$ is equipped with a choice of pullback for every display $p : E \to B$ and $f : A \to B$, and reserve the notation $f^*p : f^*E \to A$ for this pullback.

We will say that a display category is **small** if its underlying category is small.

**Definition 3.1.2.** *Let $\mathbf{C}$ be a display category. We will write $\boldsymbol{\mathcal{D}}_{\mathbf{C}}$ for the $\mathbf{C}$-indexed category with $\boldsymbol{\mathcal{D}}_{\mathbf{C}}(A)$ given by the full subcategory of $\mathbf{C}/A$ on the displays, and with 2-functorial action given by our choice of pullbacks. Similarly, we will write $(\boldsymbol{\mathcal{D}}_{\bullet})_{\mathbf{C}}$ for the $\mathbf{C}$-indexed category with $(\boldsymbol{\mathcal{D}}_{\bullet})_{\mathbf{C}}(A)$ given by the category of sections of display maps with codomain $A$, and with* 2*-functorial action given by our choice of pullbacks. We will write* $\varpi_{\mathbf{C}} : (\boldsymbol{\mathcal{D}}_{\bullet})_{\mathbf{C}} \to \boldsymbol{\mathcal{D}}_{\mathbf{C}}$ *for the forgetful* $\mathbf{C}$*-indexed functor.*

**Definition 3.1.3.** *We will call a display category* $\mathbf{C}$ ***Cartesian*** *(resp.* ***Cartesian closed****) if*

- *its underlying category is Cartesian (resp. Cartesian closed);*
- *for each* $A \in \mathbf{C}$ *the category* $\boldsymbol{\mathcal{D}}_{\mathbf{C}}(A)$ *is Cartesian (resp. Cartesian closed) and the inclusion* $\boldsymbol{\mathcal{D}}_{\mathbf{C}}(A) \hookrightarrow \mathbf{C}/A$ *preserves Cartesian (resp. Cartesian closed) structure;*
- *for each* $f : A \to B \in \mathbf{C}$*, the pullback functor* $\boldsymbol{\mathcal{D}}_{\mathbf{C}}(f) : \boldsymbol{\mathcal{D}}_{\mathbf{C}}(B) \to \boldsymbol{\mathcal{D}}_{\mathbf{C}}(A)$ *preserves Cartesian (resp. Cartesian closed) structure.*

Particularly from a type-theoretic point of view, it is important to mention the appropriate senses in which a display category may admit $\Sigma$s and $\Pi$s.

The standard way of imposing that a display category admits $\Sigma$s is simply to require that its display maps are closed under composition. This can perhaps be justified to the reader via the following triviality:

**Proposition 3.1.4.** *Let* $\mathbf{C}$ *be a display category. Then, the display maps of* $\mathbf{C}$ *are closed under composition if and only if, for all display maps* $q : B \to A$ *and* $p : E \to B$*, there exists some display map* $\Sigma_q(p) : \Sigma_q(E) \to A$ *and isomorphism*

$$\mathrm{Hom}_{\mathbf{C}/B}(p, q^*(r)) \cong \mathrm{Hom}_{\mathbf{C}/A}(\Sigma_q(p), r) \quad ,$$

*natural in* $r : D \to A$*.*

*Proof.* $\Rightarrow$. We set $\Sigma_q(p) :\equiv q \circ p$.

$\Leftarrow$. The display map $\Sigma_q(p)$ is characterized up to isomorphism as $q \circ p$. This $q \circ p$ is thus a display map, as display maps are closed under composition with isomorphisms.

□

The case of $\Pi$ is dual to that of $\Sigma$.

**Definition 3.1.5** (Streicher 1991, Definition 1.26)**.** *We will say that a display category* $\mathbf{C}$ ***admits*** $\Pi$***s*** *if, for all displays* $q : B \to A$ *and* $p : E \to B$*, there exists some display* $\Pi_q(p) : \Pi_q(E) \to A$ *and isomorphism*

$$\mathrm{Hom}_{\mathbf{C}/B}(q^*(r), p) \cong \mathrm{Hom}_{\mathbf{C}/A}(r, \Pi_q(p)) \quad ,$$

*natural in $r : D \to A$.*

Finally, we recover $\Pi$s from Cartesian closure, in the presence of enough additional structure.

**Proposition 3.1.6.** *When a display category is Cartesian, Cartesian closed, and its display maps are closed under composition, it admits $\Pi$s.*

*Proof.* Let $\mathbf{C}$ be such a display category and $q : B \to A$ and $p : E \to B$ be display maps of $\mathbf{C}$. Then, we define $\Pi_q(p) : \Pi_q(E) \to A$ as the pullback

$$\begin{array}{ccc} \Pi_q(p) & \longrightarrow & (q \circ p)^q \\ \downarrow & \lrcorner & \downarrow p^q \\ \mathrm{id}_A & \longrightarrow & q^q \end{array}$$

in $\mathcal{D}(A)$, in which the bottom morphism is the transpose of the composite of

$$(\mathrm{id}_A \times q) \cong q \xrightarrow{\mathrm{id}} q \quad .$$

□

**Universes**

Given a display category $\mathbf{C}$, one can ask the following: does there exist a discrete opfibration $p : \mathcal{U}_\bullet \to \mathcal{U}$ in $\mathbf{Cat}(\mathbf{C})$ that mediates by pullback a (2-)natural equivalence $\mathrm{Hom}_{\mathbf{Cat}(\mathbf{C})}(\Gamma, \mathcal{U}) \simeq \mathcal{D}_{\mathbf{C}}(\Gamma)$ between internal functors from $\Gamma \in \mathbf{C}$ to $\mathcal{U}$ and display maps with codomain $\Gamma$? We will call such a $p : \mathcal{U}_\bullet \to \mathcal{U}$ a universe, though it is a take on the notion of internal full subcategory stemming from Bénabou (1970s). Our treatment is particularly indebted to Street (2018)'s discussion of internal full subcategories.

In a display category with enough structure (notably in a typos or display topos), admission of a universe in our sense is equivalent to admission of a generic display map.

**Definition 3.1.7.** *Let $\mathbf{C}$ be a category. A $\mathbf{C}$-indexed functor $p : \mathbf{E} \to \mathbf{B}$ is called a **$\mathbf{C}$-indexed discrete opfibration** if the component functor $p_\Gamma : \mathbf{E}(\Gamma) \to \mathbf{B}(\Gamma)$ is a discrete opfibration, for all* $\Gamma \in \mathbf{C}$.

**Proposition 3.1.8.** *Let* $\mathbf{C}$ *be a small category. A discrete opfibration in* $\mathbf{IndCat}(\mathbf{C})$ *is the same thing as a* $\mathbf{C}$*-indexed discrete opfibration with small domain and codomain.*

*Proof.* See Cigoli et al. (2020). □

**Definition 3.1.9.** *Let* $\mathbf{C}$ *be a small category. We will say that a morphism* $p : \mathbf{E} \to \mathbf{B} \in \mathbf{IndCat}(\mathbf{C})$ *is* ***strict*** *if, for any* $\alpha : \Delta \to \Gamma \in \mathbf{C}$*, the mediating isomorphism of the naturality square*

$$\begin{array}{ccc} \mathbf{E}(\Delta) & \xleftarrow{\mathbf{E}(\alpha)} & \mathbf{E}(\Gamma) \\ {\scriptstyle p_\Delta}\downarrow & \cong & \downarrow{\scriptstyle p_\Gamma} \\ \mathbf{B}(\Delta) & \xleftarrow[\mathbf{B}(\alpha)]{} & \mathbf{B}(\Gamma) \end{array}$$

*is the identity.*

The following statements are of technical interest, and we omit the proofs.

**Proposition 3.1.10.** *Let* $\mathbf{C}$ *be a small category. A morphism of* $\mathbf{Cat}^{\mathbf{C}^{\mathrm{op}}}$ *is a discrete opfibration in* $\mathbf{Cat}^{\mathbf{C}^{\mathrm{op}}}$ *if and only if it is a (necessarily strict) discrete opfibration in* $\mathbf{IndCat}(\mathbf{C})$*.*

**Proposition 3.1.11.** *Let* $\mathbf{C}$ *be a small category. A strict discrete opfibration* $p : \mathbf{E} \to \mathbf{B} \in \mathbf{IndCat}(\mathbf{C})$ *is a discrete opfibration in* $\mathbf{Cat}^{\mathbf{C}^{\mathrm{op}}}$ *if and only if* $\mathbf{B} \in \mathbf{Cat}^{\mathbf{C}^{\mathrm{op}}}$*.*

**Proposition 3.1.12.** *Let* $\mathbf{C}$ *be a small category. A morphism of* $\mathbf{Cat}(\mathbf{C})$ *is a discrete opfibration in* $\mathbf{Cat}(\mathbf{C})$ *if and only if it is a discrete opfibration in* $\mathbf{Cat}^{\mathbf{C}^{\mathrm{op}}}$*.*

**Lemma 3.1.13.** *Let* $\mathbf{C}$ *be a category and* $p : \mathbf{E} \to \mathbf{B} \in \mathbf{Cat}(\mathbf{C})$ *a discrete opfibration. Then, a commutative square in* $\mathbf{Cat}(\mathbf{C})$ *of form*

$$\begin{array}{ccc} \Delta & \longrightarrow & \mathbf{E} \\ {\scriptstyle \alpha}\downarrow & & \downarrow{\scriptstyle p} \\ \Gamma & \longrightarrow & \mathbf{B} \end{array},$$

*in which* $\alpha : \Delta \to \Gamma \in \mathbf{C}$*, is a strict* $2$*-pullback if and only if its transpose*

$$\begin{array}{ccc} \Delta & \longrightarrow & \mathbf{E}_0 \\ {\scriptstyle \alpha}\downarrow & & \downarrow{\scriptstyle p_0} \\ \Gamma & \longrightarrow & \mathbf{B}_0 \end{array}$$

*in the discrete-underlying adjunction is a pullback in* $\mathbf{C}$*.*

**Definition 3.1.14.** *Let* $\mathbf{C}$ *be a small display category. We will say that a strict discrete opfibration* $p : \mathbf{E} \to \mathbf{B} \in \mathbf{IndCat}(\mathbf{C})$ *(resp. discrete opfibration* $p : \mathbf{E} \to \mathbf{B} \in \mathbf{Cat}^{\mathbf{C}^{\mathrm{op}}}$*) is* ***representable*** *if, for any* $\Gamma \in \mathbf{C}$ *and* $b \in \mathbf{B}(\Gamma)$*, the strict* $2$*-pullback* $b^*p : b^*\mathbf{E} \to \Gamma$ *is isomorphic to a display map in* $\mathbf{IndCat}(\mathbf{C})/\Gamma$ *(resp.* $\mathbf{Cat}^{\mathbf{C}^{\mathrm{op}}}/\Gamma$*).*

Note that, when $\mathbf{C}$ is a small display category, a discrete opfibration $p : \mathbf{E} \to \mathbf{B} \in \mathbf{Cat}^{\mathbf{C}^{\mathrm{op}}}$ is representable if and only if it is representable as a strict discrete opfibration in $\mathbf{IndCat}(\mathbf{C})$.

When $\mathbf{C}$ is a small display category, we will assume that any representable strict discrete opfibration $p : \mathbf{E} \to \mathbf{B} \in \mathbf{IndCat}(\mathbf{C})$ (resp. representable discrete opfibration $p : \mathbf{E} \to \mathbf{B} \in \mathbf{Cat}^{\mathbf{C}^{\mathrm{op}}}$) comes equipped with a choice, for each $\Gamma \in \mathbf{C}$ and $b \in \mathbf{B}(\Gamma)$, of a display map isomorphic to the strict 2-pullback $b^*p : b^*\mathbf{E} \to \Gamma$ in $\mathbf{IndCat}(\mathbf{C})/\Gamma$ (resp. $\mathbf{Cat}^{\mathbf{C}^{\mathrm{op}}}/\Gamma$).

**Definition 3.1.15.** *Let* $\mathbf{C}$ *be a display category. We say that a discrete opfibration* $p : \mathbf{E} \to \mathbf{B} \in \mathbf{Cat}(\mathbf{C})$ *is* ***display*** *if its underlying morphism* $p_0 : \mathbf{E}_0 \to \mathbf{B}_0 \in \mathbf{C}$ *is a display map.*

The following statements are of technical interest, and we omit the proofs.

**Proposition 3.1.16.** *Let* $\mathbf{C}$ *be a display category and* $p : \mathbf{E} \to \mathbf{B} \in \mathbf{Cat}(\mathbf{C})$ *a discrete opfibration. Then, the following are equivalent:*

- $p$ *is display;*
- *for any* $\Gamma \in \mathbf{C}$ *and* $F : \Gamma \to \mathbf{B}$*, there exists a strict* $2$*-pullback* $F^*p : F^*\mathbf{E} \to \Gamma$ *which is a display map;*
- *for any* $\mathbf{A} \in \mathbf{Cat}(\mathbf{C})$ *and* $F : \mathbf{A} \to \mathbf{B}$*, there exists a strict* $2$*-pullback* $F^*p : F^*\mathbf{E} \to \mathbf{A}$ *which is a display discrete opfibration.*

**Corollary 3.1.17.** *Let* $\mathbf{C}$ *be a small display category. Then, a discrete opfibration* $p : \mathbf{E} \to \mathbf{B} \in \mathbf{Cat}(\mathbf{C})$ *is display if and only if it is a representable discrete opfibration in* $\mathbf{Cat}^{\mathbf{C}^{\mathrm{op}}}$*.*

Given a display category $\mathbf{C}$, a display discrete opfibration $p : \mathbf{E} \to \mathbf{B} \in \mathbf{Cat}(\mathbf{C})$, and $F : \mathbf{A} \to \mathbf{B} \in \mathbf{Cat}(\mathbf{C})$, we will henceforth reserve the notation $F^*p : F^*\mathbf{E} \to \mathbf{A}$ for the evident strict 2-pullback with action on objects given by $(F_0)^*p_0 : (F_0)^*\mathbf{E}_0 \to \mathbf{A}_0$ and on morphisms by $(F_1)^*p_1 : (F_1)^*\mathbf{E}_1 \to \mathbf{A}_1$.

The following statements are of technical interest, and we omit the proofs.

**Proposition 3.1.18.** *Let* $\mathbf{C}$ *be a small display category. Then, a strict discrete opfibration* $p$ :

$\mathbf{E} \to \mathbf{B} \in \mathbf{IndCat}(\mathbf{C})$ *(resp. discrete opfibration* $p : \mathbf{E} \to \mathbf{B} \in \mathbf{Cat}^{\mathbf{C}^{\mathrm{op}}}$*) is representable if and only if, for any* $\mathbf{A} \in \mathbf{Cat}(\mathbf{C})$ *and* $F : \mathbf{A} \to \mathbf{B}$*, the strict* 2*-pullback* $F^*p : F^*\mathbf{E} \to \mathbf{A}$ *is isomorphic to a display discrete opfibration in* $\mathbf{IndCat}(\mathbf{C})/\mathbf{A}$ *(resp.* $\mathbf{Cat}^{\mathbf{C}^{\mathrm{op}}}/\mathbf{A}$*).*

**Corollary 3.1.19.** *Let* $\mathbf{C}$ *be a small display category. Then, given a representable discrete opfibration* $p : \mathbf{E} \to \mathbf{B} \in \mathbf{Cat}^{\mathbf{C}^{\mathrm{op}}}$*, if* $\mathbf{B}$ *is internalizable, then* $p$ *is internalizable by a display discrete opfibration.*

**Proposition 3.1.20.** *Let* $\mathbf{C}$ *be a small display category and* $\mathbf{B} \in \mathbf{IndCat}(\mathbf{C})$ *(resp.* $\mathbf{Cat}^{\mathbf{C}^{\mathrm{op}}}$*,* $\mathbf{Cat}(\mathbf{C})$*). Consider the category with objects given by strict* 2*-pullback squares of form*

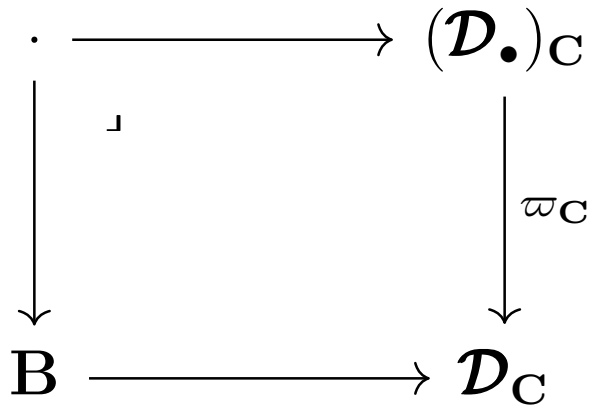

,

*in which the left vertical arrow is a representable strict discrete opfibration in* $\mathbf{IndCat}(\mathbf{C})$ *(resp. representable discrete opfibration in* $\mathbf{Cat}^{\mathbf{C}^{\mathrm{op}}}$*, display discrete opfibration in* $\mathbf{Cat}(\mathbf{C})$*), and morphisms given by pasting diagrams of form*

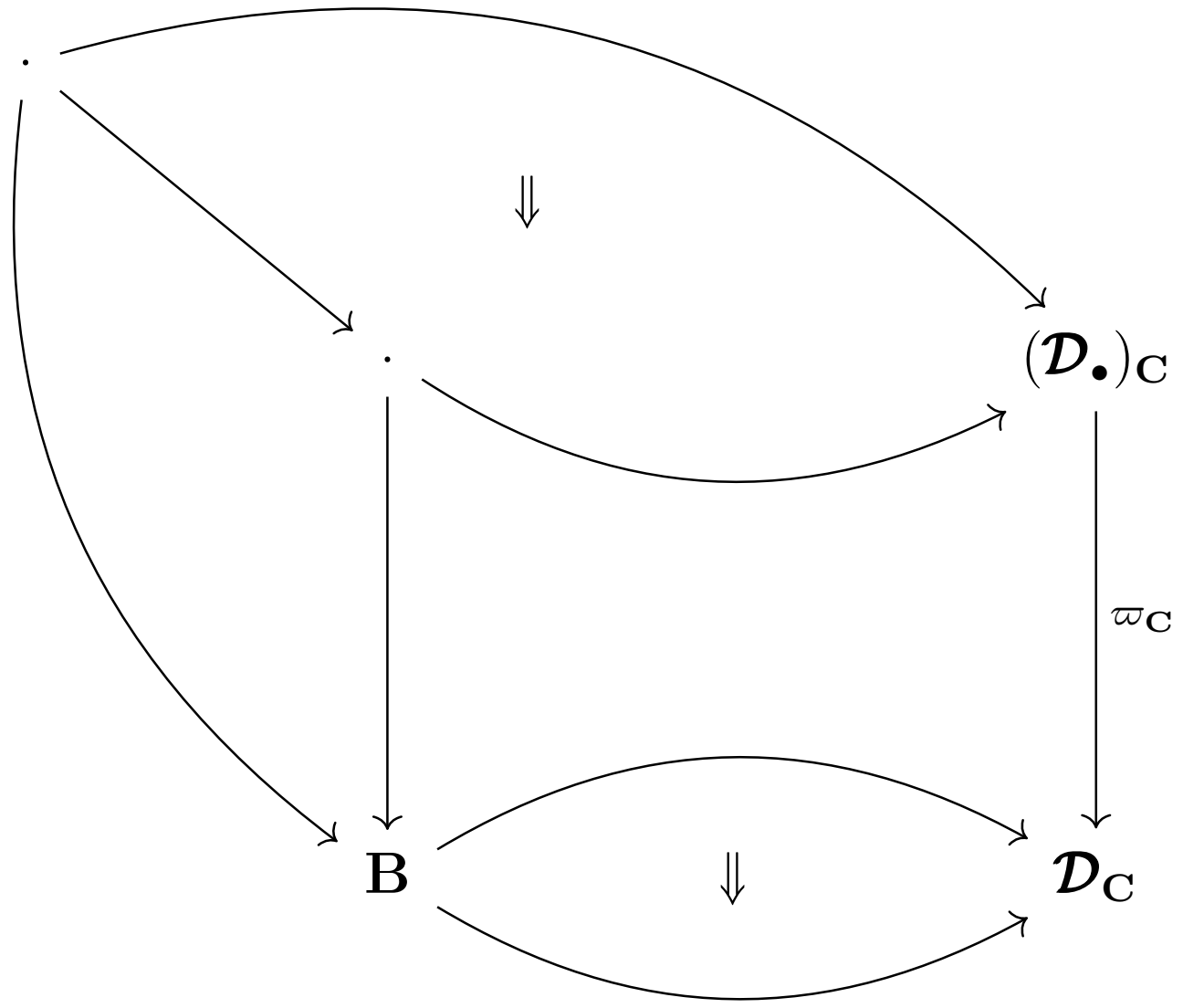

,

*in which the back and front faces are the source and target strict* 2*-pulbacks, respectively. The evident projection functors to* $\mathbf{IndCat}(\mathbf{C})(\mathbf{B}, \mathcal{D}_\mathbf{C})$ *and to the category of representable strict discrete opfibrations in* $\mathbf{IndCat}(\mathbf{C})$ *(resp. representable discrete opfibration in* $\mathbf{Cat}^{\mathbf{C}^{\mathrm{op}}}$*, display* *discrete opfibration in* $\mathbf{Cat}(\mathbf{C})$*) with codomain* $\mathbf{B}$ *are equivalences. Moreover, these equivalences are* $2$*-natural in* $\mathbf{B}$.

*Proof.* We treat the case of $\mathbf{B} \in \mathbf{IndCat}(\mathbf{C})$.

The inverse equivalence to the projection to $\mathbf{IndCat}(\mathbf{C})(\mathbf{B}, \boldsymbol{\mathcal{D}}_{\mathbf{C}})$ is given by strict 2-pullback.

As for the inverse equivalence to the other projection, we give the action on objects. Let $p : \mathbf{E} \to \mathbf{B} \in \mathbf{IndCat}(\mathbf{C})$ be a representable strict discrete opfibration. We define a $\mathbf{C}$-indexed functor $P_p : \mathbf{B} \to \boldsymbol{\mathcal{D}}_{\mathbf{C}}$ with component $P_{p,\Gamma} : \mathbf{B}(\Gamma) \to \boldsymbol{\mathcal{D}}_{\mathbf{C}}(\Gamma)$ at $\Gamma \in \mathbf{C}$ given on $b \in \mathbf{B}(\Gamma)$ by our choice of display map isomorphic to $b^*p : b^*\mathbf{E} \to \Gamma$ in $\mathbf{IndCat}(\mathbf{C})/\Gamma$. Accordingly, we define a $\mathbf{C}$-indexed functor $(P_\bullet)_p : \mathbf{E} \to (\boldsymbol{\mathcal{D}}_\bullet)_{\mathbf{C}}$ with component $(P_\bullet)_{p,\Gamma} : \mathbf{E}(\Gamma) \to (\boldsymbol{\mathcal{D}}_\bullet)_{\mathbf{C}}(\Gamma)$ at $\Gamma \in \mathbf{C}$ given on $e \in \mathbf{E}(\Gamma)$ by the section of $P_{p,\Gamma}(p_\Gamma(e))$ appearing as the dashed morphism of the diagram

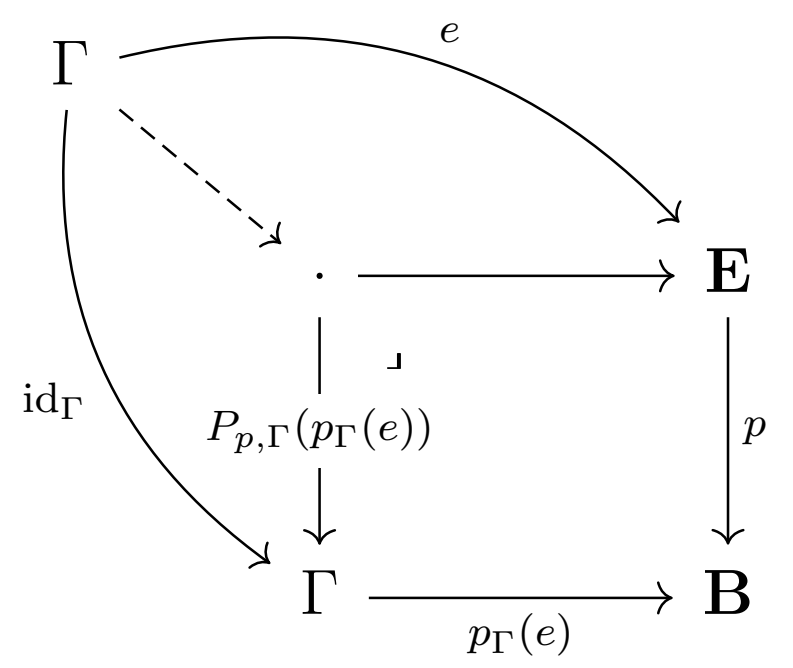

.

We thus obtain a strict 2-pullback square

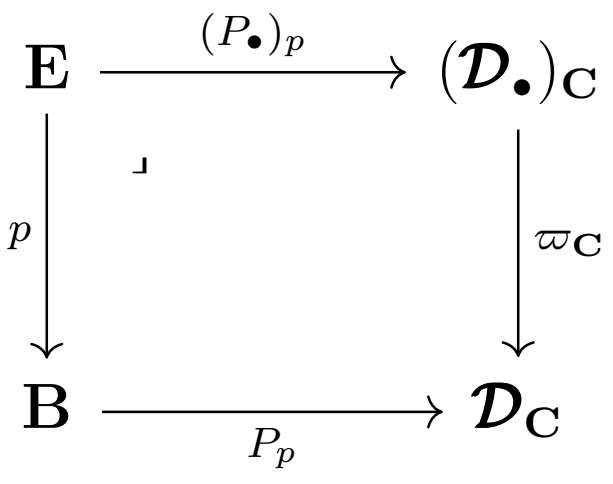

.

□

We can now define the central notions of this section:

**Definition 3.1.21.** *Let* $\mathbf{C}$ *be a small display category. We say that a representable discrete opfibration* $p : \mathbf{E} \to \mathbf{B} \in \mathbf{Cat}^{\mathbf{C}^{\mathrm{op}}}$ *(resp. display discrete opfibration* $p : \mathbf{E} \to \mathbf{B} \in \mathbf{Cat}(\mathbf{C})$*) is a* ***typing*** *(resp.* ***universe****) if the* $\mathbf{C}$*-indexed functor* $P_p : \mathbf{B} \to \boldsymbol{\mathcal{D}}_{\mathbf{C}}$ *is a* $\mathbf{C}$*-indexed equivalence.*

When $\mathbf{C}$ is a small display category, we will typically denote a typing by $p : \mathbf{Tm} \to \mathbf{Tp}$ and write $\mathrm{comp} : \mathbf{Tp} \simeq \boldsymbol{\mathcal{D}}_{\mathbf{C}}$ for $P_p : \mathbf{Tp} \simeq \boldsymbol{\mathcal{D}}_{\mathbf{C}}$ and $\mathrm{comp}_{\bullet} : \mathbf{Tm} \simeq (\boldsymbol{\mathcal{D}}_{\bullet})_{\mathbf{C}}$ for $(P_{\bullet})_p : \mathbf{Tm} \simeq (\boldsymbol{\mathcal{D}}_{\bullet})_{\mathbf{C}}$. We will typically denote a universe by $p : \mathfrak{U}_{\bullet} \to \mathfrak{U}$ and write $\mathrm{comp} : \mathfrak{U} \simeq \boldsymbol{\mathcal{D}}_{\mathbf{C}}$ for $P_p : \mathfrak{U} \simeq \boldsymbol{\mathcal{D}}_{\mathbf{C}}$ and $\mathrm{comp}_{\bullet} : \mathfrak{U}_{\bullet} \simeq (\boldsymbol{\mathcal{D}}_{\bullet})_{\mathbf{C}}$ for $(P_{\bullet})_p : \mathfrak{U}_{\bullet} \simeq (\boldsymbol{\mathcal{D}}_{\bullet})_{\mathbf{C}}$.

When $\mathbf{C}$ is a small display category, a display discrete opfibration $p : \mathfrak{U}_{\bullet} \to \mathfrak{U} \in \mathbf{Cat}(\mathbf{C})$ is a universe if and only if it is a typing. Also, in light of Corollary 3.1.19, when $\mathbf{C}$ is a small display category and $p : \mathbf{Tm} \to \mathbf{Tp} \in \mathbf{Cat}^{\mathbf{C}^{\mathrm{op}}}$ is a typing, if $\mathbf{Tp}$ is internalizable, then $p$ is internalizable by a universe in $\mathbf{Cat}(\mathbf{C})$ .

### Typoses

A typos is to be a display category-theoretic model of dependent type theory with $\Sigma$- and $\Pi$-types, as well as a universe in the current sense.

**Definition 3.1.22.** *We will say that a small, Cartesian display category is a **small typos** if*

- *it is Cartesian closed (in the sense of Definition 3.1.3);*
- *its displays are closed under composition;*
- *it admits a universe (in the sense of Definition 3.1.21).*

**Remark 3.1.23.** *Our notion of typos generalizes the eponymous one of Taylor (1987). A typos in Taylor's sense is equivalent to a typos in the present sense, in which all morphisms are display maps.*

*It was shown by Pitts and Taylor (1989) that the only typos in Taylor's sense is the terminal display category, by a version of Russell's paradox. Indeed, Taylor's notion of typos was introduced to give a categorical formulation of Russell's paradox.*

By comparing Definition 3.1.22 with Propositions 3.1.4 and 3.1.6, we observe that a typos admits $\Sigma$s and $\Pi$s in the appropriate senses.

### Display Toposes

In brief, a display topos will consist of a display category such that its underlying category is a topos and its class of displays constitutes an impredicative universe in the sense of Streicher (2005). However, we give a definition in our terms:

**Definition 3.1.24.** *We will say that a display category* ***admits a subobject classifier*** *if its underlying category admits a subobject classifier* $\top : 1 \to \Omega$ *such that both* $\top : 1 \to \Omega$ *and* $! : \Omega \to 1$ *are display maps. We will say that a small typos is a* ***small display topos*** *if it admits a subobject classifier.*

### 3.1.2 Morphisms and 2-Morphisms

We turn to morphisms and 2-morphisms of display categories, typoses, and display toposes.

**Definition 3.1.25.** *Let* $\mathbf{C}$ *and* $\mathbf{D}$ *be display categories. A* ***display functor*** $F : \mathbf{C} \to \mathbf{D}$ *consists of a functor which*

- *takes displays to displays;*
- *takes pullbacks of displays to pullbacks of displays.*

When $\mathbf{C}$ and $\mathbf{D}$ are Cartesian display categories, we say that a display functor $F : \mathbf{C} \to \mathbf{D}$ is **Cartesian** if its underlying functor is Cartesian.

We denote by $\mathbf{DC}$ the strict 2-category of small display categories, display functors, and natural transformations. We denote by $\mathbf{DC}_{\mathfrak{U}}$ the strict 2-subcategory of $\mathbf{DC}$ spanned by the small display categories that admit a universe. We denote by $\mathbf{CDC}_{\mathfrak{U}}$ the strict 2-subcategory of $\mathbf{DC}_{\mathfrak{U}}$ spanned by the small, Cartesian display categories that admit a universe and Cartesian display functors. We denote by $\mathbf{Typ}$ the strict 2-subcategory of $\mathbf{CDC}_{\mathfrak{U}}$ spanned by the small typoses. We denote by $\mathbf{DTop}$ the strict 2-subcategory of $\mathbf{Typ}$ spanned by the small display toposes.

When $\mathbf{C}$ and $\mathbf{D}$ are small display categories and $F : \mathbf{C} \to \mathbf{D}$ is a display functor, we will write $\mathcal{D}_F : \mathcal{D}_{\mathbf{C}} \to F^*\mathcal{D}_{\mathbf{D}}$ and $(\mathcal{D}_\bullet)_F : (\mathcal{D}_\bullet)_{\mathbf{C}} \to F^*(\mathcal{D}_\bullet)_{\mathbf{D}}$ for the evident induced $\mathbf{C}$-indexed functors.

When $\mathbf{C}$ and $\mathbf{D}$ are small display categories, $F : \mathbf{C} \to \mathbf{D}$ and $G : \mathbf{C} \to \mathbf{D}$ are display functors, and $\alpha : F \Rightarrow G$ is a natural transformation, we will write $\mathcal{D}_\alpha : \mathcal{D}_F \Rightarrow (\alpha^*\mathcal{D}_{\mathbf{D}} \circ \mathcal{D}_G) : \mathcal{D}_{\mathbf{C}} \to F^*\mathcal{D}_{\mathbf{D}}$ and $(\mathcal{D}_\bullet)_\alpha : (\mathcal{D}_\bullet)_F \Rightarrow (\alpha^*(\mathcal{D}_\bullet)_{\mathbf{D}} \circ (\mathcal{D}_\bullet)_G) : (\mathcal{D}_\bullet)_{\mathbf{C}} \to F^*(\mathcal{D}_\bullet)_{\mathbf{D}}$ for the evident induced $\mathbf{C}$-indexed natural transformations.

When $\mathbf{C}$ and $\mathbf{D}$ are small display categories, $F : \mathbf{C} \to \mathbf{D}$ and $G : \mathbf{C} \to \mathbf{D}$ are display functors, $\alpha : F \Rightarrow G$ is a natural transformation, and $p : E \to B \in \mathbf{C}$ is a display, we may abbreviate $((\boldsymbol{\mathcal{D}}_\alpha)_B)_p$ by $(\boldsymbol{\mathcal{D}}_\alpha)_p$.

We denote by $\boldsymbol{\mathcal{D}} : \mathbf{DC} \to \mathbf{IndCat}$ the strict 2-functor given by the assignments

$$\mathbf{C} \mapsto (\mathbf{C}, \boldsymbol{\mathcal{D}}_{\mathbf{C}})$$

$$F \mapsto (F, \boldsymbol{\mathcal{D}}_F)$$

$$\alpha \mapsto (\alpha, \boldsymbol{\mathcal{D}}_\alpha) \quad .$$

## 3.2 Natural Models, Natural Typoses, and Natural Display Toposes

> We will ... *take as fundamental concepts the notions of* $I$*-indexed families and substitution functors*.... Paré and Schumacher (1978)

To capture the setting of a (display) category equipped with abstract indexed families, we introduce a notion that we call **natural model**. This refines the eponymous notion detailed by Awodey (2018).

Our present notion of natural model has equivalent 1-category theory to the original notion, but we prefer the present notion on the basis of its 2-category theory, which is equivalent to that of full, split comprehension categories (see Section 3.4).

In brief, a (small) natural model in the present sense consists of a small display category $\mathbf{C}$, together with a typing $\mathbf{Tm}_{\mathbf{C}} \to \mathbf{Tp}_{\mathbf{C}} \in \mathbf{Cat}^{\mathbf{C}^{\mathrm{op}}}$ (see Definition 3.1.21). An object $\Gamma \in \mathbf{C}$ is called a **context**. The objects of $\mathbf{Tp}_{\mathbf{C}}(\Gamma)$ are called **types** in context $\Gamma$, though we may think of them as '$\Gamma$-indexed families of display objects.' The objects of $\mathbf{Tm}_{\mathbf{C}}(\Gamma)$ are called **terms** in context $\Gamma$, though we may think of them as '$\Gamma$-indexed families of pointed display objects.'

The 2-categorical distinction between natural models in the present sense and natural models in the original sense is simple. Given a parallel pair $F, G : \mathbf{C} \to \mathbf{D}$ of morphisms of natural models in the original sense, a 2-morphism $\alpha : F \Rightarrow G$ is equivalently given by a natural transformation between underlying functors, which we also denote by $\alpha : F \Rightarrow G$, subject to coherence conditions which imply, for each display $p : E \to B \in \mathbf{C}$, that the naturality square

$$\begin{array}{ccc} FE & \xrightarrow{\alpha_E} & GE \\ {\scriptstyle Fp}\downarrow & & \downarrow{\scriptstyle Gp} \\ FB & \xrightarrow[\alpha_B]{} & GB \end{array}$$

is a pullback.[1] By contrast, as we will see in Proposition 3.2.11, 2-morphisms of natural models in the present sense are equivalently given by arbitrary natural transformations between underlying functors.

The present, more general notion of 2-morphism of natural models is the more natural one. This is perhaps most vividly illustrated by its use in Corollary 3.2.14, which shows the 2-equivalence of typoses (resp. display toposes) and **natural typoses** (resp. **natural display toposes**). Natural typoses (resp. natural display toposes) adapt typoses (resp. display toposes) to the setting of natural models, and have the advantage of fixing a universe up to isomorphism of internal categories, rather than equivalence of internal categories.

### 3.2.1 Objects

We define some basic objects of natural model theory, up to our notion of natural display topos.

**Definition 3.2.1.** *A* ***small pre-natural model*** $\mathbf{C}$ *consists of*

- *a small display category, also denoted* $\mathbf{C}$*;*
- *a representable discrete opfibration* $p_{\mathbf{C}} : \mathbf{Tm}_{\mathbf{C}} \to \mathbf{Tp}_{\mathbf{C}} \in \mathbf{Cat}^{\mathbf{C}^{\mathrm{op}}}$.

When $\mathbf{C}$ is a small pre-natural model, we will write $\mathrm{comp}_{\mathbf{C}} : \mathbf{Tp}_{\mathbf{C}} \to \boldsymbol{\mathcal{D}}_{\mathbf{C}}$ and $(\mathrm{comp}_{\bullet})_{\mathbf{C}} : \mathbf{Tm}_{\mathbf{C}} \to (\boldsymbol{\mathcal{D}}_{\bullet})_{\mathbf{C}}$ for, respectively, the $\mathbf{C}$-indexed functors $P_{(p_{\mathbf{C}})} : \mathbf{Tp}_{\mathbf{C}} \to \boldsymbol{\mathcal{D}}_{\mathbf{C}}$ and $(P_{\bullet})_{(p_{\mathbf{C}})} : \mathbf{Tm}_{\mathbf{C}} \to (\boldsymbol{\mathcal{D}}_{\bullet})_{\mathbf{C}}$ defined in the proof of Proposition 3.1.20.

When $\mathbf{C}$ is again a small pre-natural model, we may write $\Gamma \vdash_{\mathbf{C}}$ for $\Gamma \in \mathbf{C}$, and say that $\Gamma$ is a **context** of $\mathbf{C}$. We may refer to the morphisms of $\mathbf{C}$ as **substitutions** of $\mathbf{C}$. Given $\Gamma \vdash_{\mathbf{C}}$, we may write $\Gamma \vdash_{\mathbf{C}} A$ (resp. $\Gamma \vdash_{\mathbf{C}} A \xrightarrow{f} B$) for $A : \Gamma \to \mathbf{Tp}_{\mathbf{C}}$ (resp. $f : A \Rightarrow B : \Gamma \to \mathbf{Tp}_{\mathbf{C}}$), and say that $A$ is a **type** of $\mathbf{C}$ in context $\Gamma$ (resp. that $f$ is a **function** of $\mathbf{C}$ from $A$ to $B$ in context $\Gamma$).

[1] The notion of 2-morphism we are using here is the one that may be formalistically inferred from the definition of natural model in the original sense. The corresponding notion of 2-morphism is used by Lumsdaine (2021) in the equivalent setting of categories with attributes. Of course, alternative notions of 2-morphism may be defined *ad hoc*.

Given $\Gamma \vdash_{\mathbf{C}}$, we may write $\Gamma \vdash_{\mathbf{C}} a$ for $a : \Gamma \to \mathbf{Tm}_{\mathbf{C}}$, and say that $a$ is a **term** of $\mathbf{C}$ in context $\Gamma$. We may also say that $a \Rightarrow b : \Gamma \to \mathbf{Tm}_{\mathbf{C}}$ is an **application** of $\mathbf{C}$ mapping $a$ to $b$ in context $\Gamma$. Given $\Gamma \vdash_{\mathbf{C}} A$ and $\Gamma \vdash_{\mathbf{C}} a$, we may write $\Gamma \vdash_{\mathbf{C}} a : A$ when $p_{\mathbf{C}} \circ a = A$, and say that $a$ is a term of $\mathbf{C}$ of type $A$ in context $\Gamma$. Given $\Gamma \vdash_{\mathbf{C}} a : A$, $\Gamma \vdash_{\mathbf{C}} b : B$, and $\Gamma \vdash_{\mathbf{C}} A \xrightarrow{f} B$, we may write $\Gamma \vdash_{\mathbf{C}} a \mapsto_f b$ when there exists $x : a \Rightarrow b : \Gamma \to \mathbf{Tm}_{\mathbf{C}}$ such that $p_{\mathbf{C}} \circ x = f$, and say that $f$ maps $a$ to $b$ in context $\Gamma$.

When $\mathbf{C}$ is a small pre-natural model and we have $\alpha : \Delta \to \Gamma \in \mathbf{C}$, we may write $A[\alpha]$ for $\mathbf{Tp}_{\mathbf{C}}(\alpha)(A)$, when $\Gamma \vdash_{\mathbf{C}} A$ is given. We may do similarly when given functions, terms, and applications.

When $\mathbf{C}$ is a small pre-natural model and we have $\Gamma \vdash_{\mathbf{C}} a : A$ and $\Gamma \vdash_{\mathbf{C}} A \xrightarrow{f} B$, we may write $\Gamma \vdash_{\mathbf{C}} f(a)$ for the unique term yielding $\Gamma \vdash_{\mathbf{C}} a \mapsto_f f(a)$.

When $\mathbf{C}$ is a small pre-natural model and $\Gamma \vdash_{\mathbf{C}} A$, we may write $p_A : \Gamma.A \to \Gamma$ for $\mathrm{comp}_{\mathbf{C}}(A) : \cdot \to \Gamma$. When $\Gamma \vdash_{\mathbf{C}} A \xrightarrow{f} B$, we may write $\Gamma.f : \Gamma.A \to \Gamma.B$ for the substitution underlying $\mathrm{comp}_{\mathbf{C}}(f) : \mathrm{comp}_{\mathbf{C}}(A) \to \mathrm{comp}_{\mathbf{C}}(B)$. When $\Gamma \vdash_{\mathbf{C}} a : A$, we may write $\overline{a} : \Gamma \to \Gamma.A$ for $(\mathrm{comp}_{\bullet})_{\mathbf{C}}(a) : \Gamma \to \cdot$.

When $\mathbf{C}$ is a small pre-natural model and we have $\Gamma \vdash_{\mathbf{C}} A$, we write $\Gamma.A \vdash_{\mathbf{C}} v_A$ for the unique term of $\mathbf{C}$ completing the strict 2-pullback square

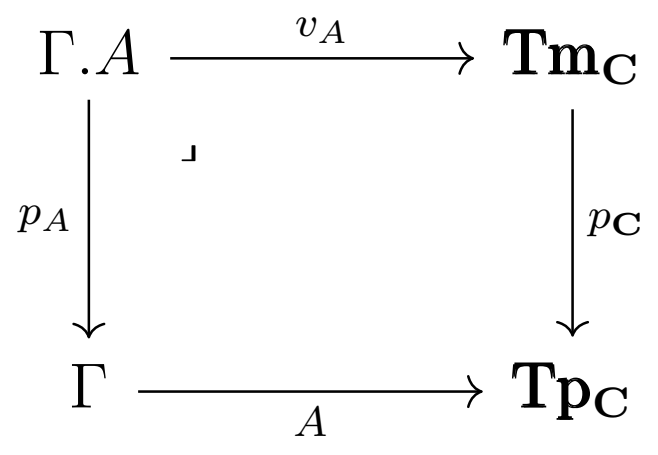

in $\mathbf{Cat}^{\mathbf{C}^{\mathrm{op}}}$.

When $\mathbf{C}$ is a small pre-natural model and we have $\Gamma \vdash_{\mathbf{C}} A$ and $\alpha : \Delta \to \Gamma \in \mathbf{C}$, we write $q(A, \alpha) : \Delta.A[\alpha] \to \Gamma.A$ for the unique substitution of $\mathbf{C}$ completing the diagram

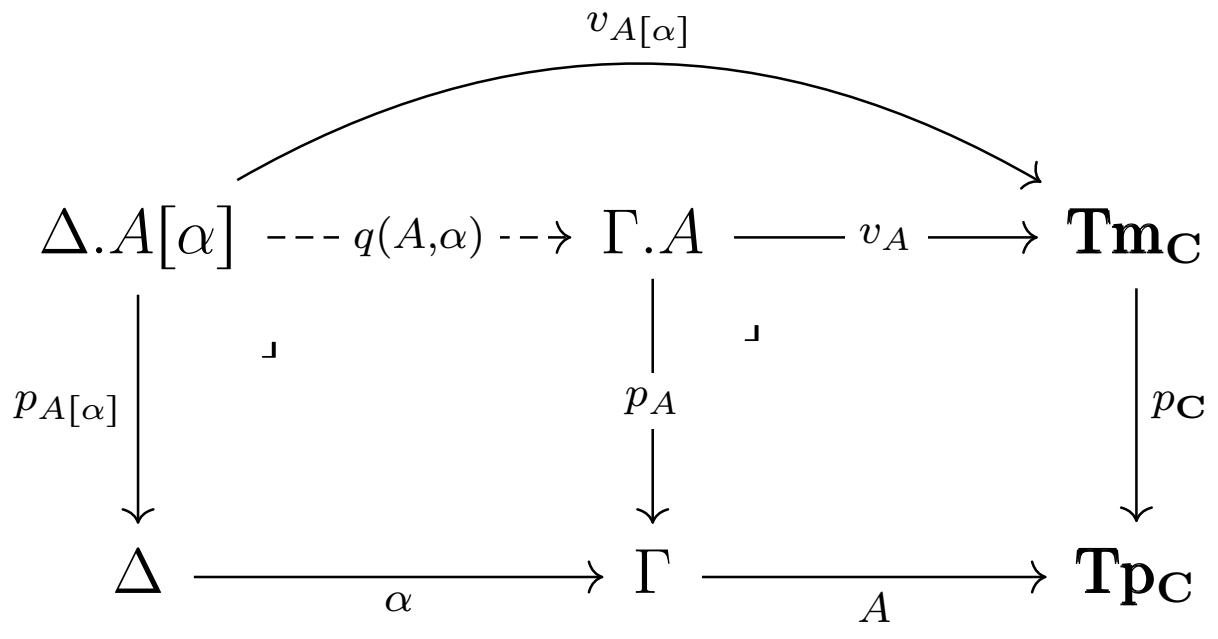

in $\mathbf{Cat}^{\mathbf{C}^{\mathrm{op}}}$.

**Definition 3.2.2.** *We will say that a small pre-natural model* $\mathbf{C}$ *is a* ***small natural model*** *if the representable discrete opfibration* $p_{\mathbf{C}} : \mathbf{Tm}_{\mathbf{C}} \to \mathbf{Tp}_{\mathbf{C}} \in \mathbf{Cat}^{\mathbf{C}^{\mathrm{op}}}$ *is a typing (for which see Definition 3.1.21).*

We will say that a small natural model is **Cartesian** (resp. **Cartesian closed**) if its underlying display category is Cartesian (resp. Cartesian closed).

**Definition 3.2.3** (*cf.* Moggi 1991, Definition 6.6)**.** *Let* $\mathbf{C}$ *be a small natural model. We will say that a category* $\mho$ *in* $\mathbf{C}$ *is a* ***type classifier*** *if it internalizes* $\mathbf{Tp}_{\mathbf{C}}$ *(in the sense of §2.1).*

It follows that, when $\mathbf{C}$ is a small display category, $\mathbf{C}$ admits a universe if and only if it admits the structure of a small natural model that admits a type classifier.

**Definition 3.2.4.** *We say that a Cartesian small natural model is a* ***small natural typos*** *if*

- *it admits a type classifier;*
- *its underlying Cartesian display category is a small typos.*

Since both conditions imply that the underlying display category admits a universe, the brevity of Definition 3.2.4 comes at the cost of some redundancy.

We say that a small natural typos is a **small natural display topos** if its underlying small typos is a small display topos.

### 3.2.2 Morphisms and 2-Morphisms

We turn to morphisms and 2-morphisms of natural models, natural typoses, and natural display toposes.

In order to define morphisms of natural models, we first define morphisms of pre-natural models (*cf.* the premorphisms of natural models of Newstead, 2018).

**Definition 3.2.5.** *Let* $\mathbf{C}$ *and* $\mathbf{D}$ *be small pre-natural models. A* ***morphism*** $F : \mathbf{C} \to \mathbf{D}$ *of small pre-natural models consists of:*

- *a display functor, also denoted* $F : \mathbf{C} \to \mathbf{D}$*;*
- *strict* $\mathbf{C}$*-indexed functors* $\mathbf{Tp}_F : \mathbf{Tp}_{\mathbf{C}} \to F^*\mathbf{Tp}_{\mathbf{D}}$ *and* $\mathbf{Tm}_F : \mathbf{Tm}_{\mathbf{C}} \to F^*\mathbf{Tm}_{\mathbf{D}}$*, such that the square*

$$\begin{array}{ccc} \mathbf{Tm_C} & \xrightarrow{\mathbf{Tm}_F} & F^*\mathbf{Tm_D} \\ {\scriptstyle p_\mathbf{C}}\downarrow & & \downarrow{\scriptstyle F^*p_\mathbf{D}} \\ \mathbf{Tp_C} & \xrightarrow[\mathbf{Tp}_F]{} & F^*\mathbf{Tp_D} \end{array}$$

*commutes in* $\mathbf{Cat}^{\mathbf{C}^{\mathrm{op}}}$.

When $F : \mathbf{C} \to \mathbf{D}$ is a morphism of small pre-natural models, we will typically write $FA$ for $(\mathbf{Tp}_F)_\Gamma A$, where $\Gamma \vdash_\mathbf{C} A$, and similarly for functions, terms, and applications.

When $F : \mathbf{C} \to \mathbf{D}$ is a morphism of small pre-natural models and we have $\Gamma \vdash_\mathbf{C} A$, we denote by $(\tau_F)_A : F(\Gamma.A) \to F\Gamma.FA$ the substitution of $\mathbf{D}$ appearing as the dashed arrow of the commutative diagram

$$\begin{array}{ccc} F\Gamma.FA & \xrightarrow{v_{FA}} & \mathbf{Tm_D} \\ {\scriptstyle p_{FA}}\downarrow & \lrcorner & \downarrow{\scriptstyle p_\mathbf{D}} \\ F\Gamma & \xrightarrow[FA]{} & \mathbf{Tp_D} \end{array} \qquad F(\Gamma.A) \xrightarrow{Fv_A} \mathbf{Tm_D}, \quad F(\Gamma.A) \overset{(\tau_F)_A}{\dashrightarrow} F\Gamma.FA, \quad F(\Gamma.A) \xrightarrow{Fp_A} F\Gamma$$

in $\mathbf{Cat}^{\mathbf{D}^{\mathrm{op}}}$. We denote by $(\tau_F)_A : Fp_A \to p_{FA}$ the morphism of $\boldsymbol{\mathcal{D}}_\mathbf{D}(F\Gamma)$ with $(\tau_F)_A : F(\Gamma.A) \to F\Gamma.FA$ underlying. Moreover, when we have $\Gamma \vdash_\mathbf{C} a : A$, we denote by $((\tau_\bullet)_F)_a : F\overline{a} \to \overline{Fa}$ the morphism of $(\boldsymbol{\mathcal{D}}_\bullet)_\mathbf{D}(F\Gamma)$ with $(\tau_F)_A : F(\Gamma.A) \to F\Gamma.FA$ underlying.

When $F : \mathbf{C} \to \mathbf{D}$ is a morphism of small pre-natural models, let

$$\begin{array}{ccc} \mathbf{Tp_C} & \xrightarrow{\mathbf{Tp}_F} & F^*\mathbf{Tp_D} \\ {\scriptstyle \mathrm{comp}_\mathbf{C}}\downarrow & \tau_F \Uparrow & \downarrow{\scriptstyle F^*\,\mathrm{comp}_\mathbf{D}} \\ \boldsymbol{\mathcal{D}}_\mathbf{C} & \xrightarrow[\boldsymbol{\mathcal{D}}_F]{} & F^*\boldsymbol{\mathcal{D}}_\mathbf{D} \end{array}$$

denote the $\mathbf{C}$-indexed natural transformation with component at context $\Gamma$ given by the assignment

$$A \mapsto (Fp_A \xrightarrow{(\tau_F)_A} p_{FA}) \quad .$$

Moreover, let

$$
\begin{array}{ccc}
\mathbf{Tm}_{\mathbf{C}} & \xrightarrow{\mathbf{Tm}_F} & F^*\mathbf{Tm}_{\mathbf{D}} \\
{\scriptstyle (\mathrm{comp}_\bullet)_{\mathbf{C}}}\big\downarrow & (\tau_\bullet)_F \Uparrow & \big\downarrow{\scriptstyle F^*(\mathrm{comp}_\bullet)_{\mathbf{D}}} \\
(\mathcal{D}_\bullet)_{\mathbf{C}} & \xrightarrow[(\mathcal{D}_\bullet)_F]{} & F^*(\mathcal{D}_\bullet)_{\mathbf{D}}
\end{array}
$$

denote the $\mathbf{C}$-indexed natural transformation with component at context $\Gamma$ given by the assignment

$$a \mapsto (F\overline{a} \xrightarrow{((\tau_\bullet)_F)_a} \overline{Fa}) \quad .$$

The following lemma is then immediate.

**Lemma 3.2.6.** *Let $F : \mathbf{C} \to \mathbf{D}$ be a morphism of small pre-natural models. Then, the pasting diagram*

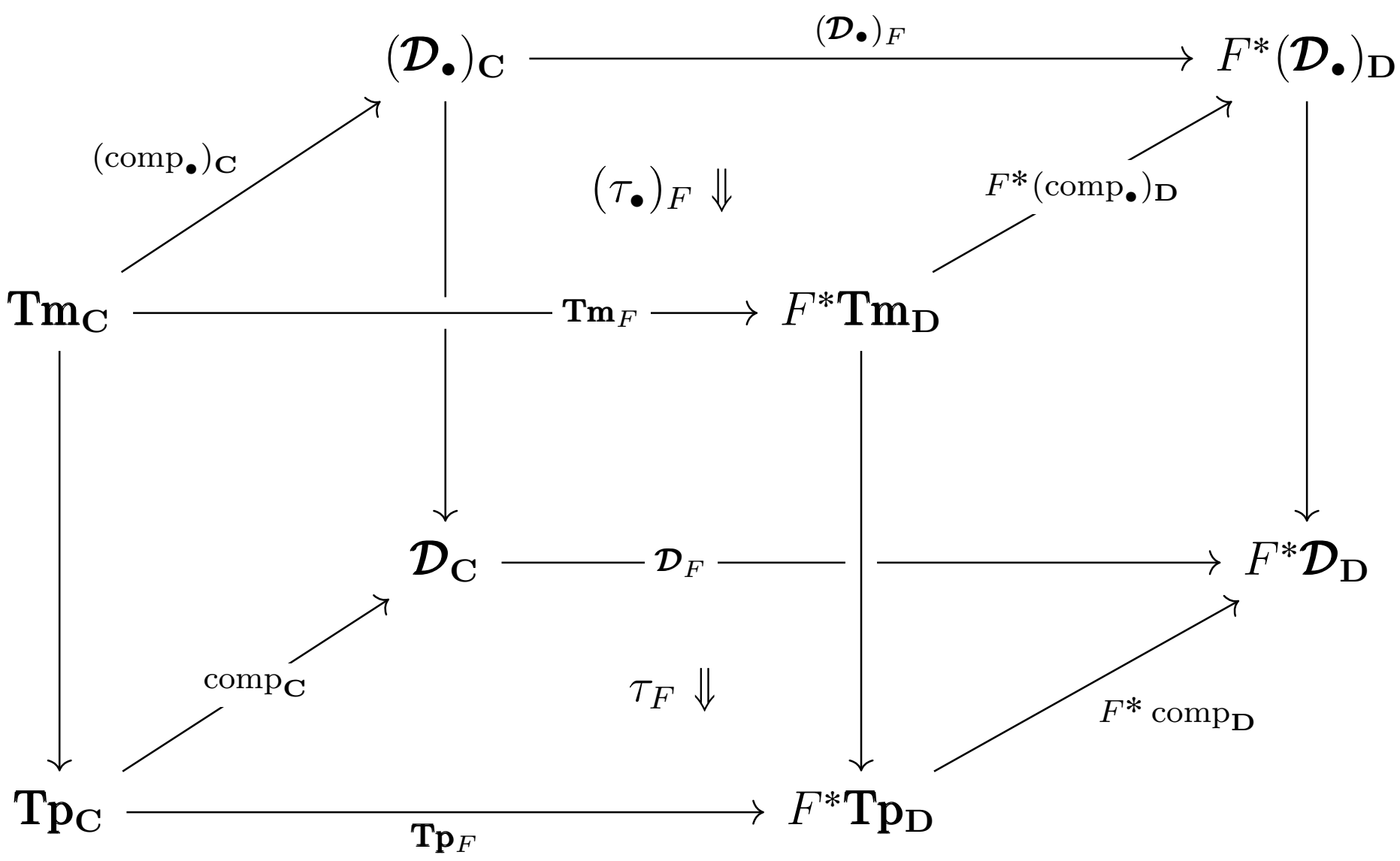

*commutes in* $\mathbf{IndCat}(\mathbf{C})$.

**Definition 3.2.7.** *Let $\mathbf{C}$ and $\mathbf{D}$ be small pre-natural models, and $F : \mathbf{C} \to \mathbf{D}$ and $G : \mathbf{C} \to \mathbf{D}$ morphisms of small pre-natural models. A* $2$***-morphism*** *of small pre-natural models $\alpha : F \Rightarrow G$ consists of:*

- *a natural transformation, also denoted $\alpha : F \Rightarrow G$;*
- $\mathbf{C}$*-indexed natural transformations $\mathbf{Tp}_\alpha : \mathbf{Tp}_F \Rightarrow (\alpha^*\mathbf{Tp}_{\mathbf{D}} \circ \mathbf{Tp}_G)$ and $\mathbf{Tm}_\alpha : \mathbf{Tm}_F \Rightarrow (\alpha^*\mathbf{Tm}_{\mathbf{D}} \circ \mathbf{Tm}_G)$, such that the pasting diagram*

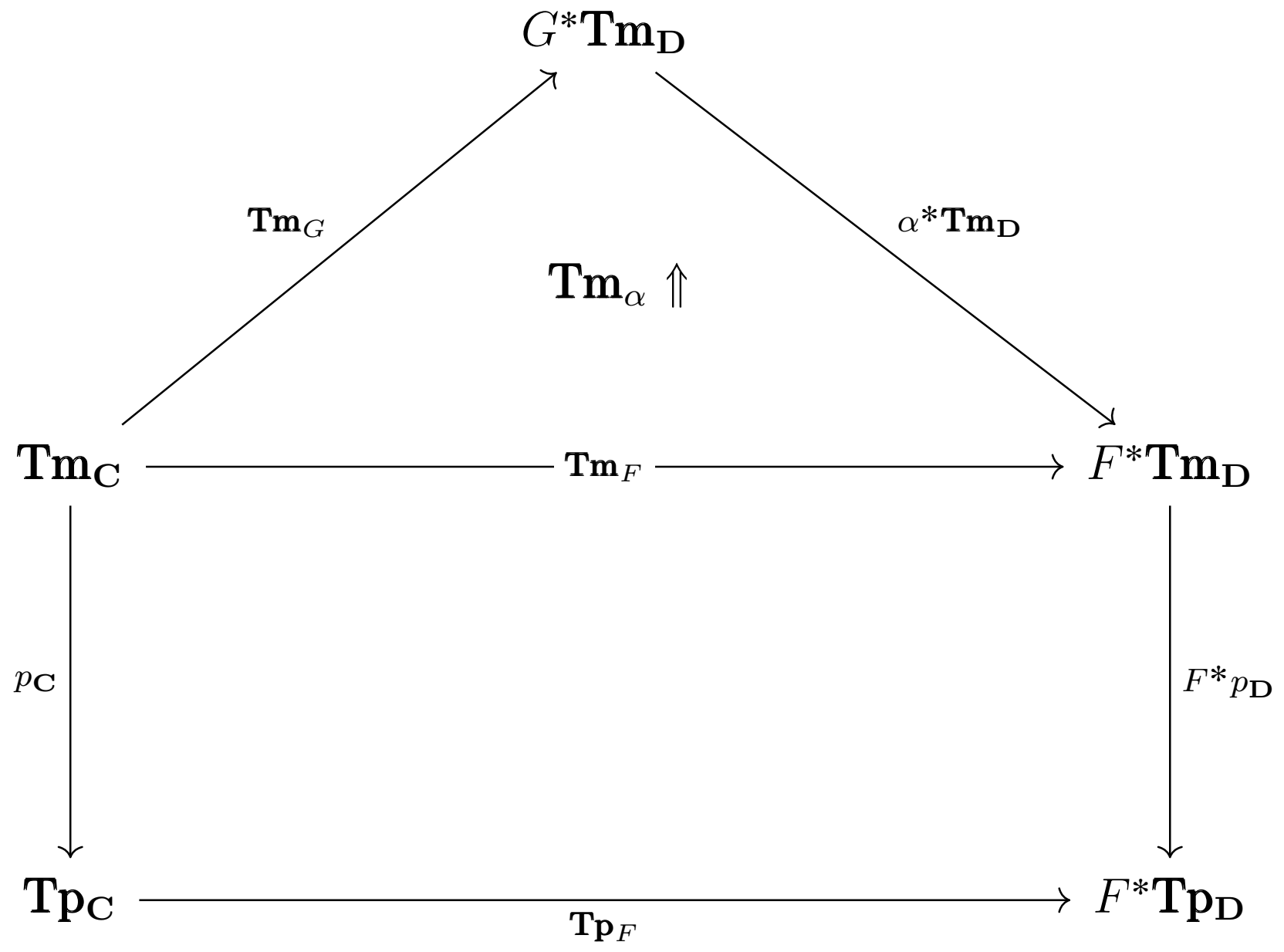

*is equal to the pasting diagram*

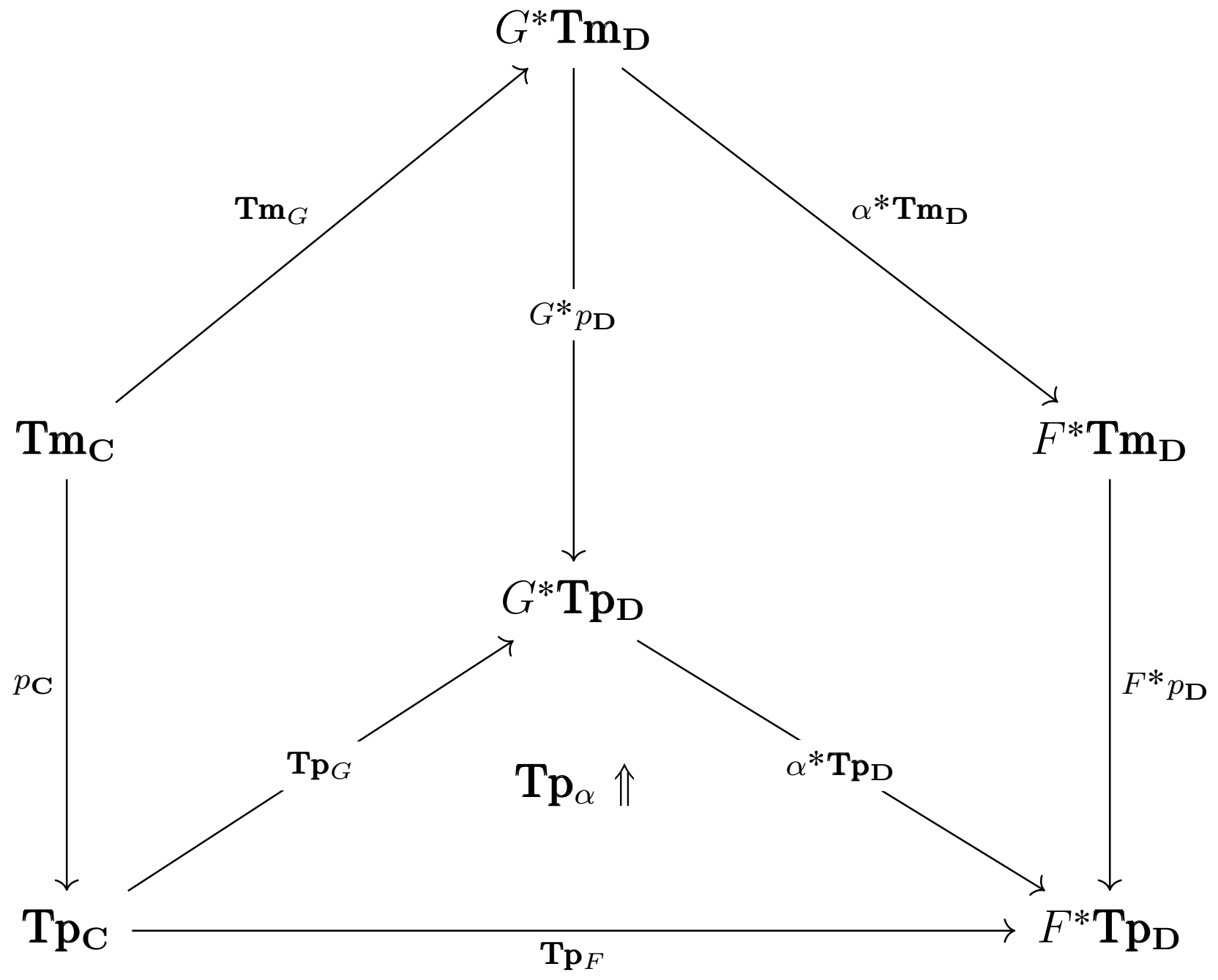

*in* $\mathbf{Cat}^{\mathbf{C}^{\text{op}}}$.

When $\alpha : F \Rightarrow G : \mathbf{C} \to \mathbf{D}$ is a 2-morphism of small pre-natural models, we will typically write $\alpha_A$ for $((\mathbf{Tp}_\alpha)_\Gamma)_A$, where $\Gamma \vdash_\mathbf{C} A$, and similarly for terms.

**Lemma 3.2.8.** *Let* $\alpha : F \Rightarrow G : \mathbf{C} \to \mathbf{D}$ *be a* 2*-morphism of small pre-natural models. Then, the pasting diagram*

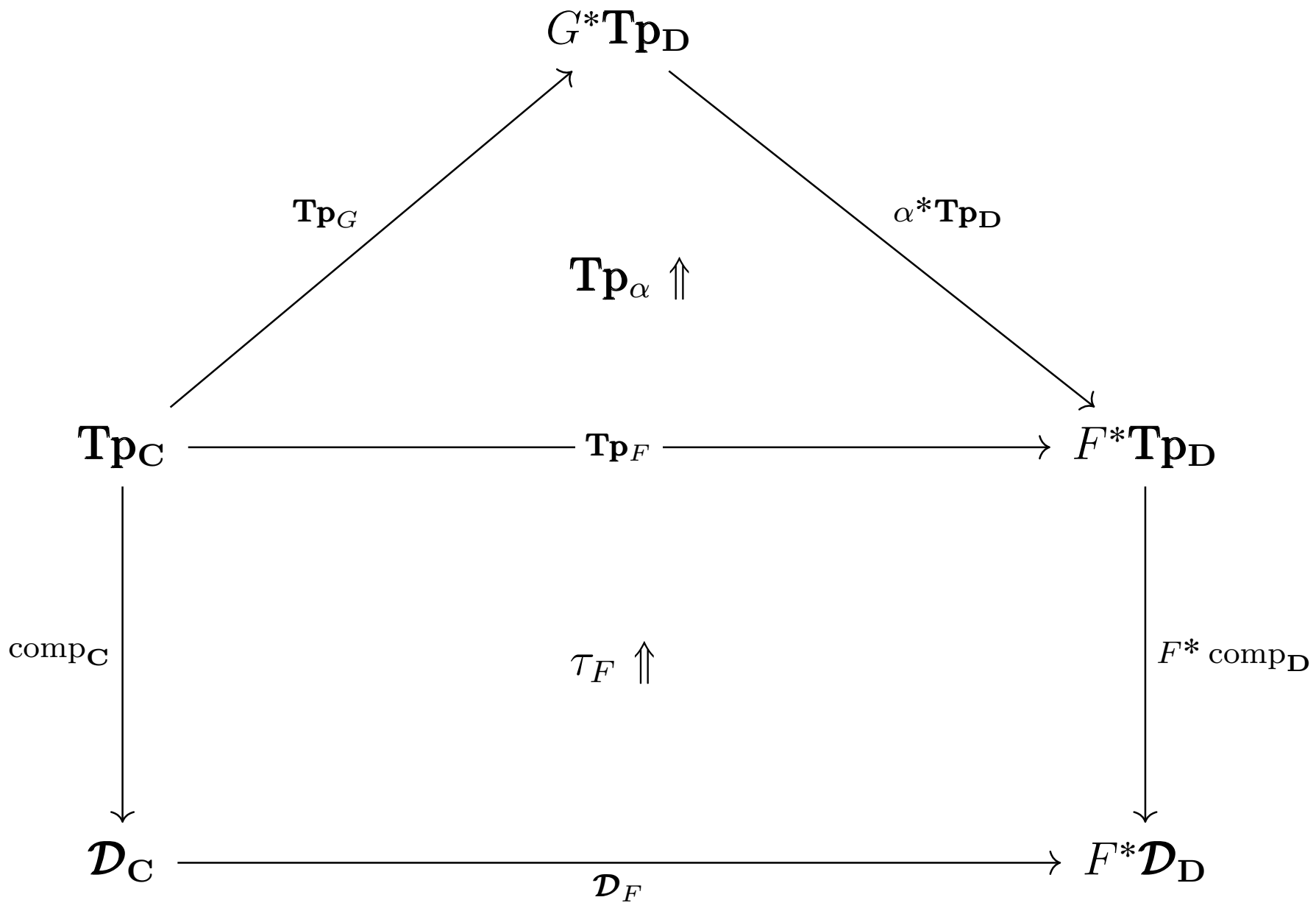

*is equal to the pasting diagram*

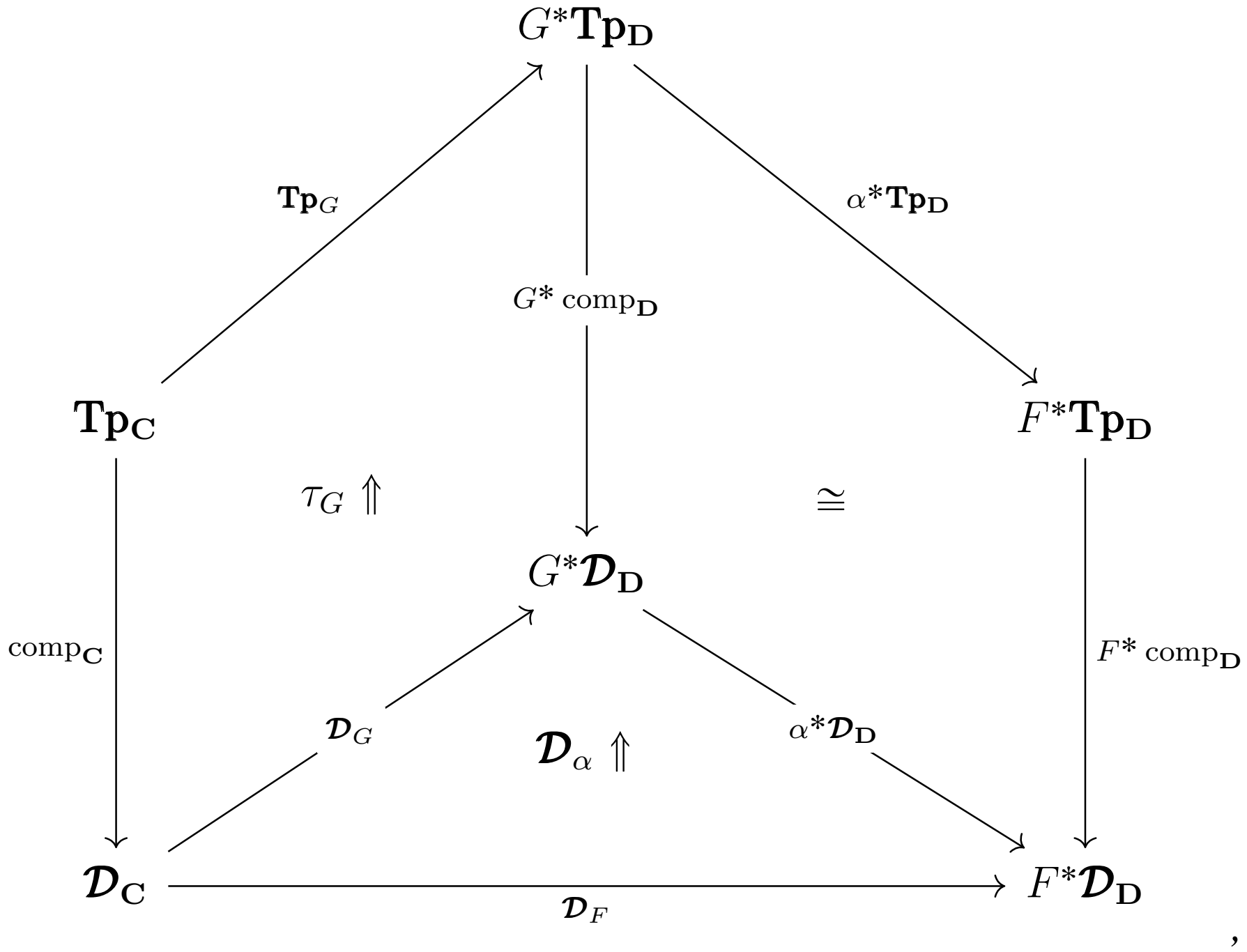

*and the pasting diagram*

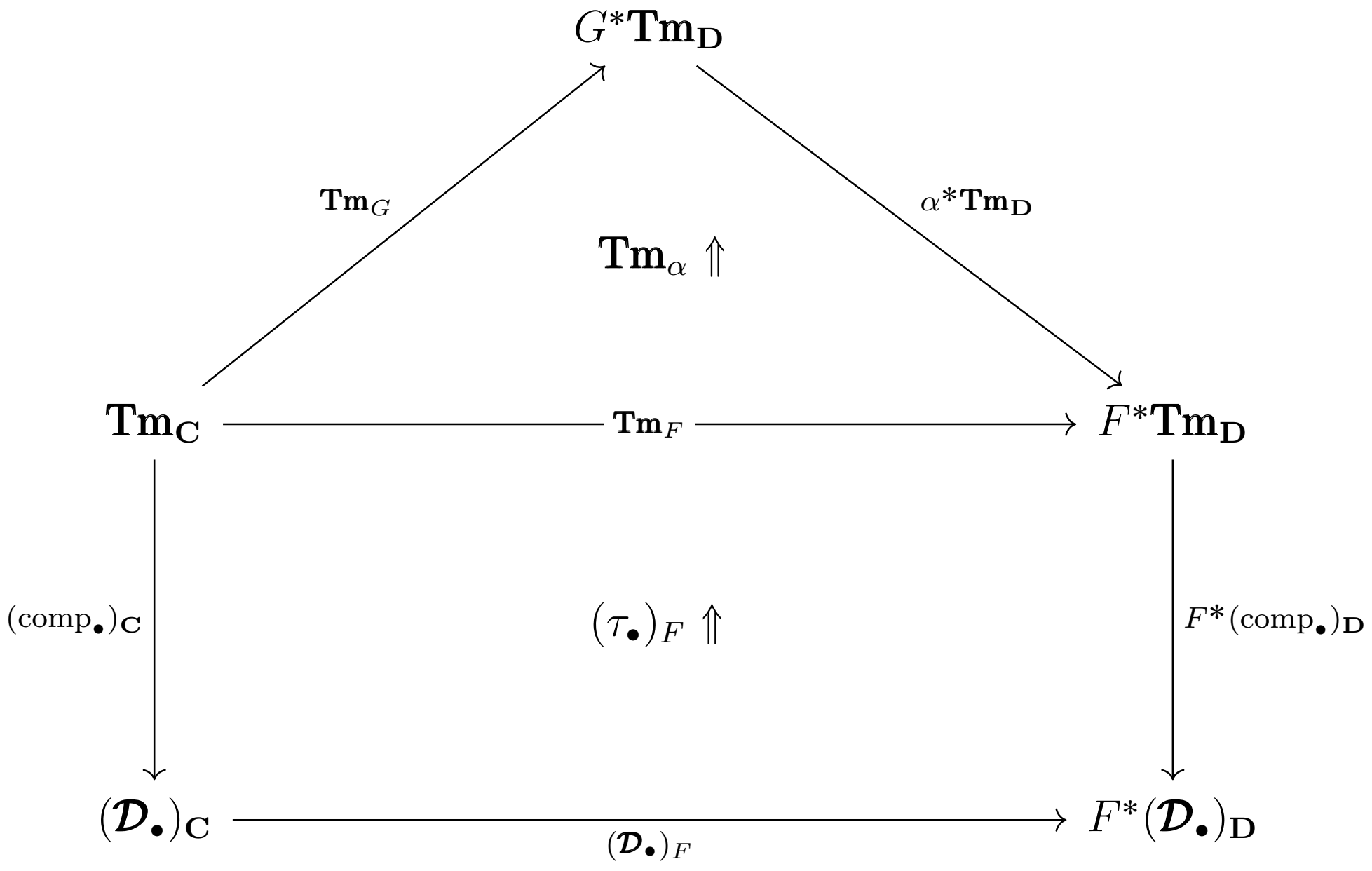

*is equal to the pasting diagram*

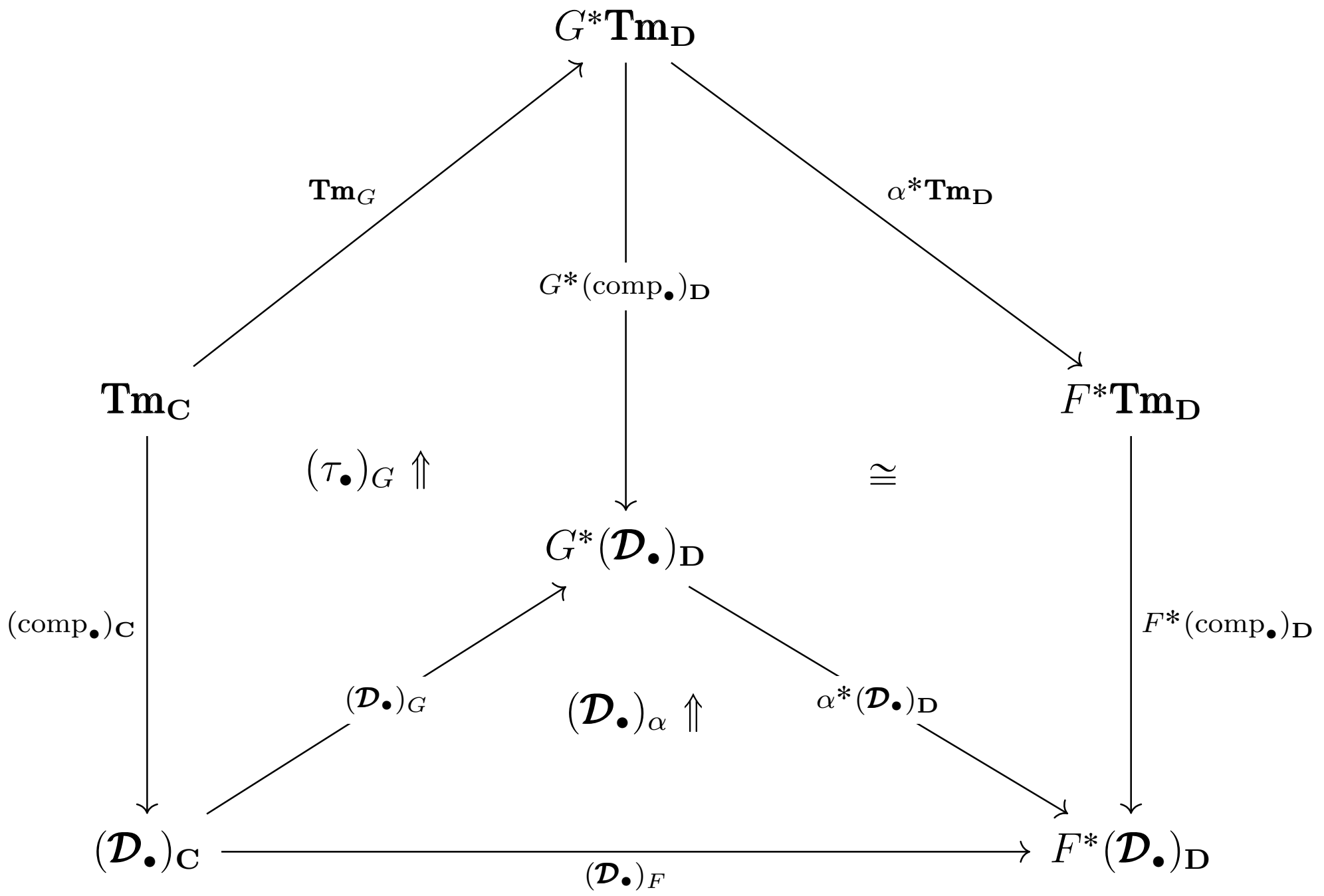

,

*both in* $\mathbf{IndCat}(\mathbf{C})$.

*Proof.* For the first equality of pasting diagrams, we must have, whenever given $\Gamma \vdash_{\mathbf{C}} A$, a commutative diagram

$$\begin{array}{ccccc}
F(\Gamma.A) & & \xrightarrow{(\tau_F)_A} & & F\Gamma.FA \\
\Big\downarrow {\scriptstyle (\boldsymbol{\mathcal{D}}_\alpha)_{(p_A)}} & & & & \Big\downarrow {\scriptstyle F\Gamma.\alpha_A} \\
(\alpha_\Gamma)^* G(\Gamma.A) & \xrightarrow[(\alpha_\Gamma)^*(\tau_G)_A]{} & (\alpha_\Gamma)^*(G\Gamma.GA) & \xrightarrow[(\rho_{(\alpha_\Gamma)})_{GA}]{\simeq} & F\Gamma.(GA[\alpha_\Gamma])
\end{array}$$

in $\mathbf{D}$. As $F\Gamma.(GA[\alpha_\Gamma])$ is a (strict 2-)pullback, it suffices to show that

$$p_{GA[\alpha_\Gamma]} \circ (\rho_{(\alpha_\Gamma)})_{GA} \circ (\alpha_\Gamma)^*(\tau_G)_A \circ (\boldsymbol{\mathcal{D}}_\alpha)_{(p_A)} = p_{GA[\alpha_\Gamma]} \circ F\Gamma.\alpha_A \circ (\tau_F)_A$$

and

$$v_{GA[\alpha_\Gamma]} \circ (\rho_{(\alpha_\Gamma)})_{GA} \circ (\alpha_\Gamma)^*(\tau_G)_A \circ (\boldsymbol{\mathcal{D}}_\alpha)_{(p_A)} = v_{GA[\alpha_\Gamma]} \circ F\Gamma.\alpha_A \circ (\tau_F)_A \quad .$$

The former equation is easy. As for the latter, we have

$$\begin{aligned}
v_{GA[\alpha_\Gamma]} \circ (\rho_{(\alpha_\Gamma)})_{GA} \circ (\alpha_\Gamma)^*(\tau_G)_A \circ (\boldsymbol{\mathcal{D}}_\alpha)_{(p_A)} &= v_{GA} \circ q(\alpha_\Gamma, GA) \circ (\rho_{(\alpha_\Gamma)})_{GA} \circ (\alpha_\Gamma)^*(\tau_G)_A \circ (\boldsymbol{\mathcal{D}}_\alpha)_{(p_A)} \\
&= v_{GA} \circ (p_{GA})^*\alpha_\Gamma \circ (\alpha_\Gamma)^*(\tau_G)_A \circ (\boldsymbol{\mathcal{D}}_\alpha)_{(p_A)} \\
&= v_{GA} \circ (\tau_G)_A \circ (Gp_A)^*\alpha_\Gamma \circ (\boldsymbol{\mathcal{D}}_\alpha)_{(p_A)} \\
&= Gv_A \circ (Gp_A)^*\alpha_\Gamma \circ (\boldsymbol{\mathcal{D}}_\alpha)_{(p_A)} \\
&= Gv_A \circ \alpha_{\Gamma.A} \\
&= (p_{\mathbf{D}} \circ \alpha_{(v_A)})(Fv_A) \\
&= (\alpha_{(p_{\mathbf{C}} \circ v_A)})(Fv_A) \\
&= (\alpha_{(A \circ p_A)})(Fv_A) \\
&= (\alpha_A \circ Fp_A)(Fv_A) \\
&= (\alpha_A \circ p_{FA} \circ (\tau_F)_A)(v_{FA} \circ (\tau_F)_A) \\
&= (\alpha_A \circ p_{FA})(v_{FA}) \circ (\tau_F)_A \\
&= v_{GA[\alpha_\Gamma]} \circ F\Gamma.\alpha_A \circ (\tau_F)_A \quad .
\end{aligned}$$

The first equality of pasting diagrams may then be lifted to the second, using the fact that $F^*\varpi_{\mathbf{D}}$ is a discrete opfibration in $\mathbf{IndCat}(\mathbf{C})$. $\square$

**Remark 3.2.9.** *In the form of Lemma 3.2.8, we have recovered the coherence condition stipulated by Castellan et al. (2017, Appendix B).*

**Definition 3.2.10.** *Let* $\mathbf{C}$ *and* $\mathbf{D}$ *be small natural models. We say that a morphism* $F : \mathbf{C} \to \mathbf{D}$ *of small pre-natural models is a* ***morphism of small natural models*** *if the* $\mathbf{C}$*-indexed natural transformation*

$$\begin{array}{ccc} \mathbf{Tp}_{\mathbf{C}} & \xrightarrow{\mathbf{Tp}_F} & F^*\mathbf{Tp}_{\mathbf{D}} \\ {\scriptstyle \mathrm{comp}_{\mathbf{C}}}\downarrow & \tau_F \Uparrow & \downarrow{\scriptstyle F^*\mathrm{comp}_{\mathbf{D}}} \\ \mathcal{D}_{\mathbf{C}} & \xrightarrow[\mathcal{D}_F]{} & F^*\mathcal{D}_{\mathbf{D}} \end{array}$$

*is a* $\mathbf{C}$*-indexed natural isomorphism.*

When $\mathbf{C}$ and $\mathbf{D}$ are Cartesian small natural models, we will say that a morphism of natural models $F : \mathbf{C} \to \mathbf{D}$ is **Cartesian** if its underlying display functor is Cartesian. When $\mathbf{C}$ and $\mathbf{D}$ are, moreover, small natural typoses, we may say that $F$ is a **natural Cartesian display functor** in this case.

We denote by $\mathbf{PreNM}$ the strict 2-category of small pre-natural models and morphisms and 2-morphisms thereof. We denote by $\mathbf{NM}$ the strict 2-subcategory of $\mathbf{PreNM}$ spanned by the small natural models and morphisms thereof. We denote by $\mathbf{NM}_{\mho}$ the strict 2-subcategory of $\mathbf{NM}$ spanned by the small natural models that admit a type classifier. We denote by $\mathbf{CNM}_{\mho}$ the strict 2-subcategory of $\mathbf{NM}_{\mho}$ spanned by the Cartesian small natural models that admit a type classifier and Cartesian morphisms of small natural models. We denote by $\mathbf{NTyp}$ the strict 2-subcategory of $\mathbf{CNM}_{\mho}$ spanned by the small natural typoses. We denote by $\mathbf{NDTop}$ the strict 2-subcategory of $\mathbf{NTyp}$ spanned by the small natural display toposes.

**Proposition 3.2.11.** *The forgetful strict* 2*-functor* $\mathbf{NM} \to \mathbf{DC}$ *is bijective on* 2*-morphisms.*

*Proof.* Given a 2-morphism $\alpha : F \Rightarrow G : \mathbf{C} \to \mathbf{D}$ of natural models, the component $\mathbf{Tp}_\alpha : \mathbf{Tp}_F \Rightarrow (\alpha^*\mathbf{Tp}_{\mathbf{D}} \circ \mathbf{Tp}_G) : \mathbf{Tp}_{\mathbf{C}} \to F^*\mathbf{Tp}_{\mathbf{D}}$ (resp. $\mathbf{Tm}_\alpha : \mathbf{Tm}_F \Rightarrow (\alpha^*\mathbf{Tm}_{\mathbf{D}} \circ \mathbf{Tm}_G) : \mathbf{Tm}_{\mathbf{C}} \to F^*\mathbf{Tm}_{\mathbf{D}}$) may be reconstructed from the eponymous underlying natural transformation $\alpha : F \Rightarrow G : \mathbf{C} \to \mathbf{D}$ as the conjugate of $\mathcal{D}_\alpha : \mathcal{D}_F \Rightarrow (\alpha^*\mathcal{D}_{\mathbf{D}} \circ \mathcal{D}_G) : \mathcal{D}_{\mathbf{C}} \to F^*\mathcal{D}_{\mathbf{D}}$ (resp. $(\mathcal{D}_\bullet)_\alpha : (\mathcal{D}_\bullet)_F \Rightarrow (\alpha^*(\mathcal{D}_\bullet)_{\mathbf{D}} \circ (\mathcal{D}_\bullet)_G) : (\mathcal{D}_\bullet)_{\mathbf{C}} \to F^*(\mathcal{D}_\bullet)_{\mathbf{D}}$). □

**Corollary 3.2.12.** *The forgetful strict* $2$*-functors* $\mathbf{NM}_{\mho} \to \mathbf{DC}_{\mho}$*,* $\mathbf{CNM}_{\mho} \to \mathbf{CDC}_{\mho}$*,* $\mathbf{NTyp} \to \mathbf{Typ}$*, and* $\mathbf{NDTop} \to \mathbf{DTop}$ *are bijective on* $2$*-morphisms.*

The notion of display category that admits a universe is *a priori* distinct from the notion of natural model that admits a type classifier. Indeed, the latter is more clearly akin to the notion of display category *equipped* with a universe. Nevertheless, these notions are equivalent in a weak sense: they give rise to 2-equivalent strict 2-categories.

**Theorem 3.2.13.** *The forgetful strict* $2$*-functor* $\mathbf{NM}_{\mho} \to \mathbf{DC}_{\mho}$ *is a* $2$*-equivalence.*

*Proof.* It suffices to show that the forgetful strict 2-functor $\mathbf{NM}_{\mho} \to \mathbf{DC}_{\mho}$ is essentially surjective on objects and morphisms and bijective on 2-morphisms. Bijectivity on 2-morphisms was already established by Corollary 3.2.12.

As for essential surjectivity on objects, let $\mathbf{C}$ be a small display category that admits a universe $p : \mho_{\bullet} \to \mho \in \mathbf{Cat}(\mathbf{C})$. Then, $p$ is a small natural model structure on $\mathbf{C}$ for which $\mho$ is a type classifier.

Towards essential surjectivity on morphisms, let $\mathbf{C}$ and $\mathbf{D}$ be small natural models that each admit a type classifier and $F : \mathbf{C} \to \mathbf{D}$ a display functor between the underlying display categories. We will extend $F$ with the structure of a morphism of small natural models.

Via a choice of type classifiers, we attain universes $p_{\mathbf{C}} : (\mho_{\bullet})_{\mathbf{C}} \to \mho_{\mathbf{C}}$ for $\mathbf{C}$ and $p_{\mathbf{D}} : (\mho_{\bullet})_{\mathbf{D}} \to \mho_{\mathbf{D}}$ for $\mathbf{D}$. Henceforth, we conflate the natural model structure on $\mathbf{C}$ with $p_{\mathbf{C}}$ and that on $\mathbf{D}$ with $p_{\mathbf{D}}$.

Since $F$ is a display functor and $p_{\mathbf{C}}$ is a display discrete opfibration, $F((p_{\mathbf{C}})_0)$ is a display map, and, since $p_{\mathbf{D}}$ is a universe, we may then choose a strict 2-pullback square

$$\begin{array}{ccc} F(((\mho_{\bullet})_{\mathbf{C}})_0) & \xrightarrow{\phi_{\bullet}} & (\mho_{\bullet})_{\mathbf{D}} \\ {\scriptstyle F((p_{\mathbf{C}})_0)}\downarrow & \lrcorner & \downarrow{\scriptstyle p_{\mathbf{D}}} \\ F((\mho_{\mathbf{C}})_0) & \xrightarrow[\phi]{} & \mho_{\mathbf{D}} \end{array}$$

in $\mathbf{Cat}(\mathbf{D})$.

Given a type $A : \Gamma \to \mho_{\mathbf{C}} \in \mathbf{Cat}(\mathbf{C})$, we define $FA : F\Gamma \to \mho_{\mathbf{D}} \in \mathbf{Cat}(\mathbf{D})$ as $\phi \circ F(A_0)$. Similarly, given a term $a : \Gamma \to (\mho_{\bullet})_{\mathbf{C}} \in \mathbf{Cat}(\mathbf{C})$, we define $Fa : F\Gamma \to (\mho_{\bullet})_{\mathbf{D}} \in \mathbf{Cat}(\mathbf{D})$ as

$\phi_\bullet \circ F(a_0)$.

At this stage, we obtain, for any type $A : \Gamma \to \mathfrak{U}_{\mathbf{C}} \in \mathbf{Cat}(\mathbf{C})$, an isomorphism $(\tau_F)_A : F(\Gamma.A) \cong F\Gamma.FA \in \mathbf{D}$ as the coherence isomorphism of strict 2-pullbacks in the commutative diagram

$$\begin{array}{ccccc}
F(\Gamma.A) & & \xrightarrow{\quad Fv_A \quad} & & \\
 & \overset{(\tau_F)_A}{\searrow} & & \searrow & \\
 & & F\Gamma.FA & \xrightarrow{v_{FA}} & (\mathfrak{U}_\bullet)_{\mathbf{D}} \\
 & {\scriptstyle Fp_A} \searrow & \downarrow{\scriptstyle p_{FA}} & & \downarrow{\scriptstyle p_{\mathbf{D}}} \\
 & & F\Gamma & \xrightarrow[FA]{} & \mathfrak{U}_{\mathbf{D}}
\end{array}$$

in $\mathbf{Cat}(\mathbf{D})$.

Given a function $f : A \Rightarrow B : \Gamma \to \mathfrak{U}_{\mathbf{C}} \in \mathbf{Cat}(\mathbf{C})$, we define $Ff : FA \Rightarrow FB : F\Gamma \to \mathfrak{U}_{\mathbf{D}} \in \mathbf{Cat}(\mathbf{D})$ as the unique function yielding the commutative square

$$\begin{array}{ccc}
F(\Gamma.A) & \xrightarrow[(\tau_F)_A]{\simeq} & F\Gamma.FA \\
{\scriptstyle F(\Gamma.f)}\downarrow & & \downarrow{\scriptstyle F\Gamma.Ff} \\
F(\Gamma.B) & \xrightarrow[(\tau_F)_B]{\simeq} & F\Gamma.FB
\end{array}$$

in $\mathbf{D}$, which exists because $p_{\mathbf{D}}$ is a universe.

The action of $F$ on applications is left to the reader. □

We also note the following, for future reference:

**Corollary 3.2.14.** *The forgetful strict* 2*-functors* $\mathbf{CNM}_{\mathfrak{U}} \to \mathbf{CDC}_{\mathfrak{U}}$, $\mathbf{NTyp} \to \mathbf{Typ}$*, and* $\mathbf{NDTop} \to \mathbf{DTop}$ *are* 2*-equivalences.*

## 3.3 Natural Models and the Literature

Natural models in the present sense distinguish themselves in cases where 2-morphisms of models are relevant. We give a non-exhaustive list of examples that have cropped up in the literature.

We will gloss over the issue of terminal objects.

**Example 3.3.1.** *The documents Clairambault and Dybjer (2011) and Castellan et al. (2017, Appendix B) lay out a equivalence between the* $2$*-category of small locally Cartesian closed categories and a certain* $2$*-category of small CwFs. The* $2$*-morphisms of CwFs used for the result are,* a posteriori*, equivalent to natural transformations of underlying functors, as required, but are defined* ad hoc*. This result can thus be given a more canonical form by replacing CwFs with natural models in the present sense, from which our apposite notion of* $2$*-morphism (Definition 3.2.7) falls out automatically.*

**Example 3.3.2.** *The author has previously emphasized notions of 'geometric morphism' and 'Cartesian comonad' for small natural models (Zwanziger, 2019b). These may be justified as equivalent to adjunctions and comonads in* $\mathbf{NM}$*, respectively. Equivalent notions for small CwFs were earlier touched on by Nuyts (2018).*

**Example 3.3.3.** *Riley et al. (2021) introduced the notion of a* ***CwF with a bireflector*** *(****CwB****) in order to model a modal dependent type theory for parametrized spectra. A small CwB may be justified as equivalent to a bireflector in* $\mathbf{NM}$, i.e. *an endomorphism* $N : \mathbf{C} \to \mathbf{C}$ *together with* $2$*-morphisms* $\eta : \mathrm{id}_{\mathbf{C}} \to N$ *and* $\varepsilon : N \to \mathrm{id}_{\mathbf{C}}$ *such that* $\eta \circ \varepsilon = \mathrm{id}_N$*, all in* $\mathbf{NM}$*. This is a special case of a comonad in* $\mathbf{NM}$.

**Example 3.3.4.** *Birkedal et al. (2020) introduced the notion of a* ***CwF + A*** *to model a modal dependent type theory. A small CwF + A is equivalent to an endomorphism* $f_* : \mathbf{C} \to \mathbf{C}$ *in* $\mathbf{NM}$*, together with a left adjoint* $f^* : \mathbf{C} \to \mathbf{C}$ *to the underlying functor* $f_* : \mathbf{C} \to \mathbf{C}$*. It would be natural to additionally require the structure of a morphism on* $f^*$*, yielding the notion of an endoadjunction in* $\mathbf{NM}$.

**Example 3.3.5.** *Nuyts (2018) introduced an* ad hoc $2$*-category of (small) CwFs which is equivalent to our* $\mathbf{NM}$*. He works with several examples of morphisms and adjunctions in this* $2$*-category during his study of presheaf models of type theory.*

**Example 3.3.6.** *Coquand et al. (2019) introduced a notion of* ***lex operation*** *on a (small) CwF model of type theory. The relevant part of this definition of lex operation may be justified as equivalent to a pointed endomorphism in* $\mathbf{NM}$, i.e. *an endomorphism* $L : \mathbf{C} \to \mathbf{C}$*, together with a* $2$*-morphism* $\eta : \mathrm{id}_{\mathbf{C}} \to L$*, both in* $\mathbf{NM}$.

## 3.4 Categories with Attributes

Full, split comprehension categories (FSCCs; see Jacobs 1993) are significant for us, as our 2-category $\mathbf{NM}$ of small natural models is equivalent to a suitable 2-category of small FSCCs.

Since we make no use of the general notion of comprehension category, the term full, split comprehension category is cumbersome. We thus use the term **category with attributes** (**CwA**) in place of the term FSCC. CwAs in the present sense are thus regarded as refining the eponymous gadgets of Cartmell (1978), just as natural models in the present sense refine natural models in the sense of Awodey (2018).

In this section, we define our 2-category $\mathbf{CwA}$ of such (small) CwAs and sketch the equivalence with $\mathbf{NM}$. We also characterize $\mathbf{CwA}$ as an equi-comma object.

### 3.4.1 Objects

**Definition 3.4.1.** *A **small pre-category with attributes (small pre-CwA)** $\mathbf{C}$ consists of*

- *a small display category, also denoted* $\mathbf{C}$*;*
- *a small, strict* $\mathbf{C}$*-indexed category* $\mathbf{Tp}_{\mathbf{C}} \in \mathbf{Cat}^{\mathbf{C}^{\mathrm{op}}}$*;*
- *a* $\mathbf{C}$*-indexed functor* $\mathrm{comp}_{\mathbf{C}} : \mathbf{Tp}_{\mathbf{C}} \to \boldsymbol{\mathcal{D}}_{\mathbf{C}} \in \mathbf{IndCat}(\mathbf{C})$.

When $\mathbf{C}$ is a small pre-CwA, we may write $\Gamma \vdash_{\mathbf{C}}$ for $\Gamma \in \mathbf{C}$, and say that $\Gamma$ is a **context** of $\mathbf{C}$. We may refer to the morphisms of $\mathbf{C}$ as **substitutions** of $\mathbf{C}$. Given $\Gamma \vdash_{\mathbf{C}}$, we may write $\Gamma \vdash_{\mathbf{C}} A$ (resp. $\Gamma \vdash_{\mathbf{C}} A \xrightarrow{f} B$) for $A : \Gamma \to \mathbf{Tp}_{\mathbf{C}}$ (resp. $f : A \Rightarrow B : \Gamma \to \mathbf{Tp}_{\mathbf{C}}$), and say that $A$ is a **type** of $\mathbf{C}$ in context $\Gamma$ (resp. that $f$ is a **function** of $\mathbf{C}$ from $A$ to $B$ in context $\Gamma$).

When $\mathbf{C}$ is again a small pre-CwA, we denote by

$$\begin{array}{ccc} \mathbf{Tm}_{\mathbf{C}} & \xrightarrow{(\mathrm{comp}_{\bullet})_{\mathbf{C}}} & (\boldsymbol{\mathcal{D}}_{\bullet})_{\mathbf{C}} \\ {\scriptstyle p_{\mathbf{C}}}\downarrow & \lrcorner & \downarrow{\scriptstyle \varpi_{\mathbf{C}}} \\ \mathbf{Tp}_{\mathbf{C}} & \xrightarrow[\mathrm{comp}_{\mathbf{C}}]{} & \boldsymbol{\mathcal{D}}_{\mathbf{C}} \end{array}$$

the strict 2-pullback in $\mathbf{IndCat}(\mathbf{C})$.

**Definition 3.4.2.** *We say that a small pre-CwA* $\mathbf{C}$ *is a **small category with attributes (small CwA)** if the* $\mathbf{C}$*-indexed functor* $\mathrm{comp}_{\mathbf{C}} : \mathbf{Tp}_{\mathbf{C}} \to \boldsymbol{\mathcal{D}}_{\mathbf{C}}$ *is a* $\mathbf{C}$*-indexed equivalence.*

## 3.4.2 Morphisms and 2-Morphisms

**Definition 3.4.3.** *Let* $\mathbf{C}$ *and* $\mathbf{D}$ *be small pre-CwAs. A* ***morphism*** *of small pre-CwAs* $F : \mathbf{C} \to \mathbf{D}$ *consists of:*

- *a display functor, also denoted* $F : \mathbf{C} \to \mathbf{D}$*;*
- *a strict* $\mathbf{C}$*-indexed functor* $\mathbf{Tp}_F : \mathbf{Tp}_{\mathbf{C}} \to F^*\mathbf{Tp}_{\mathbf{D}} \in \mathbf{Cat}^{\mathbf{C}^{\mathrm{op}}}$*;*
- *a* $\mathbf{C}$*-indexed natural transformation*

$$\begin{array}{ccc} \mathbf{Tp}_{\mathbf{C}} & \xrightarrow{\mathbf{Tp}_F} & F^*\mathbf{Tp}_{\mathbf{D}} \\ {\scriptstyle \mathrm{comp}_{\mathbf{C}}}\downarrow & \tau_F \Uparrow & \downarrow{\scriptstyle F^* \mathrm{comp}_{\mathbf{D}}} \\ \mathcal{D}_{\mathbf{C}} & \xrightarrow[\mathcal{D}_F]{} & F^*\mathcal{D}_{\mathbf{D}} \end{array}.$$

When $F : \mathbf{C} \to \mathbf{D}$ is a morphism of small pre-CwAs, we denote by $\mathbf{Tm}_F : \mathbf{Tm}_{\mathbf{C}} \to F^*\mathbf{Tm}_{\mathbf{D}} \in \mathbf{Cat}^{\mathbf{C}^{\mathrm{op}}}$ and

$$\begin{array}{ccc} \mathbf{Tm}_{\mathbf{C}} & \xrightarrow{\mathbf{Tm}_F} & F^*\mathbf{Tm}_{\mathbf{D}} \\ {\scriptstyle (\mathrm{comp}_\bullet)_{\mathbf{C}}}\downarrow & (\tau_\bullet)_F \Uparrow & \downarrow{\scriptstyle F^* (\mathrm{comp}_\bullet)_{\mathbf{D}}} \\ (\mathcal{D}_\bullet)_{\mathbf{C}} & \xrightarrow[(\mathcal{D}_\bullet)_F]{} & F^*(\mathcal{D}_\bullet)_{\mathbf{D}} \end{array}$$

the unique data such that the pasting diagram

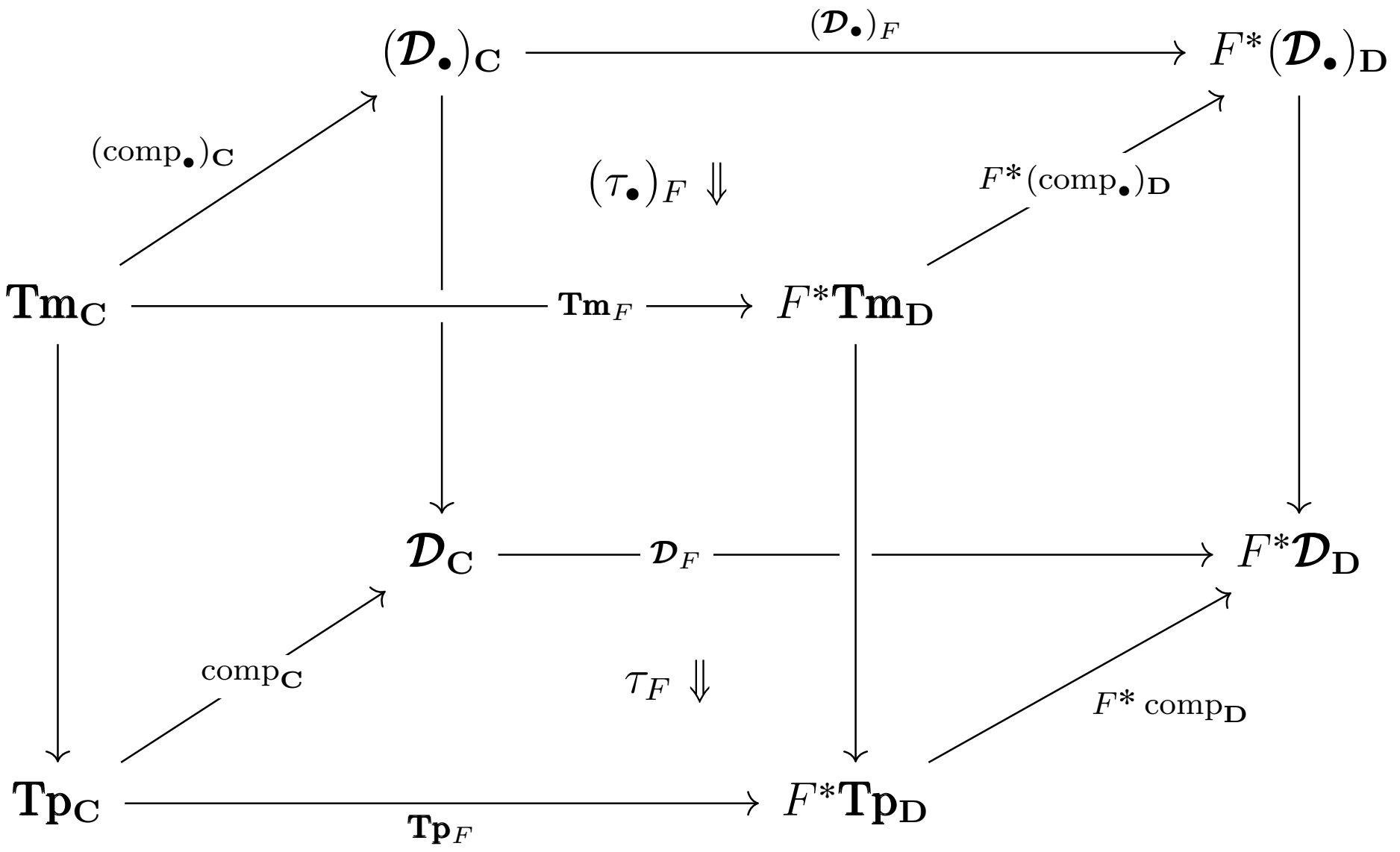

commutes in $\mathbf{IndCat}(\mathbf{C})$.

**Definition 3.4.4.** *Let* $\mathbf{C}$ *and* $\mathbf{D}$ *be small pre-CwAs, and* $F : \mathbf{C} \to \mathbf{D}$ *and* $G : \mathbf{C} \to \mathbf{D}$ *morphisms* *of pre-CwAs. A 2-**morphism** of small pre-CwAs $\alpha : F \Rightarrow G$ consists of:*

- *a natural transformation, also denoted $\alpha : F \Rightarrow G$;*
- *a $\mathbf{C}$-indexed natural transformation $\mathbf{Tp}_\alpha : \mathbf{Tp}_F \Rightarrow (\alpha^*\mathbf{Tp}_\mathbf{D} \circ \mathbf{Tp}_G)$, such that the pasting diagram*

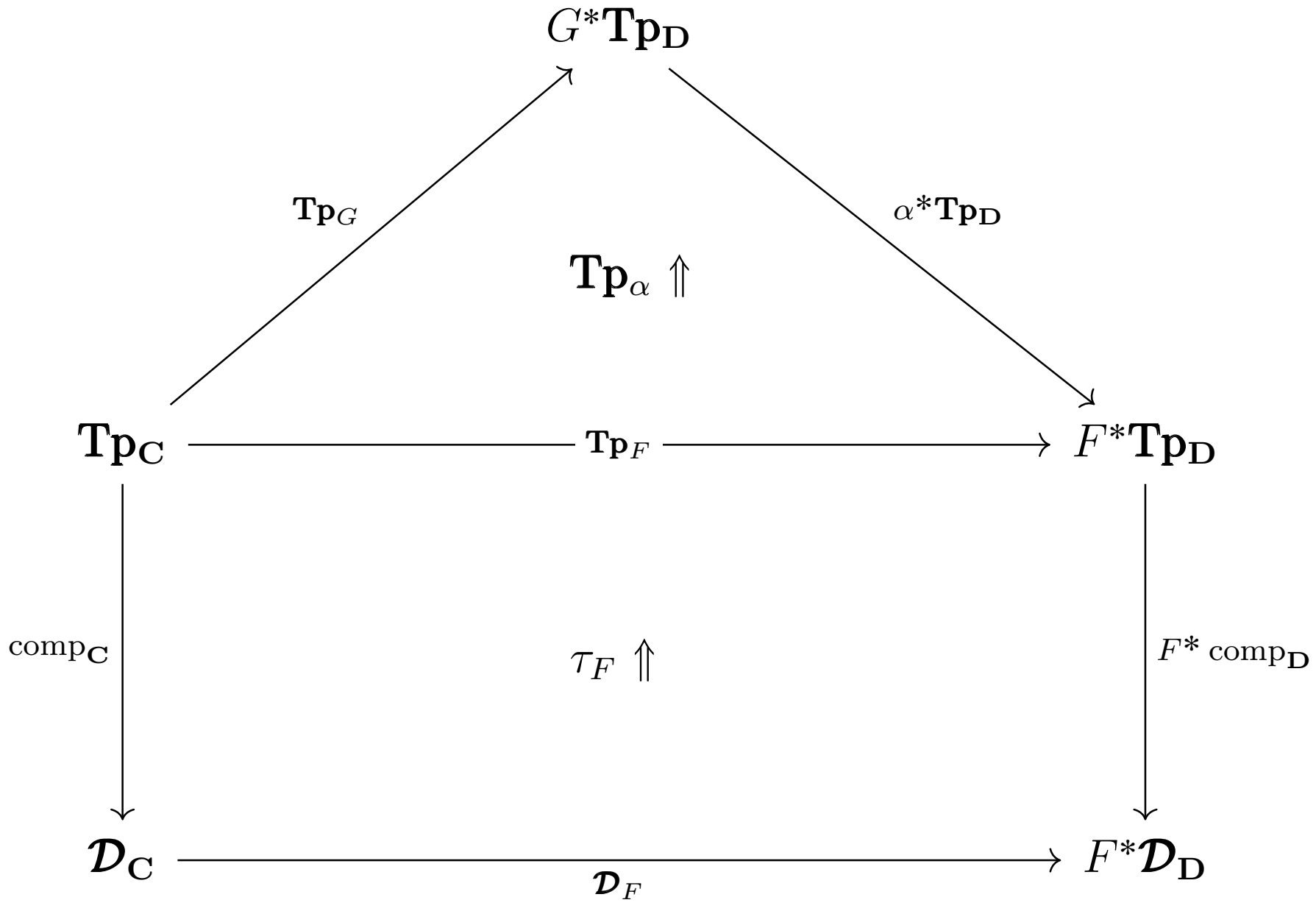

*is equal to the pasting diagram*

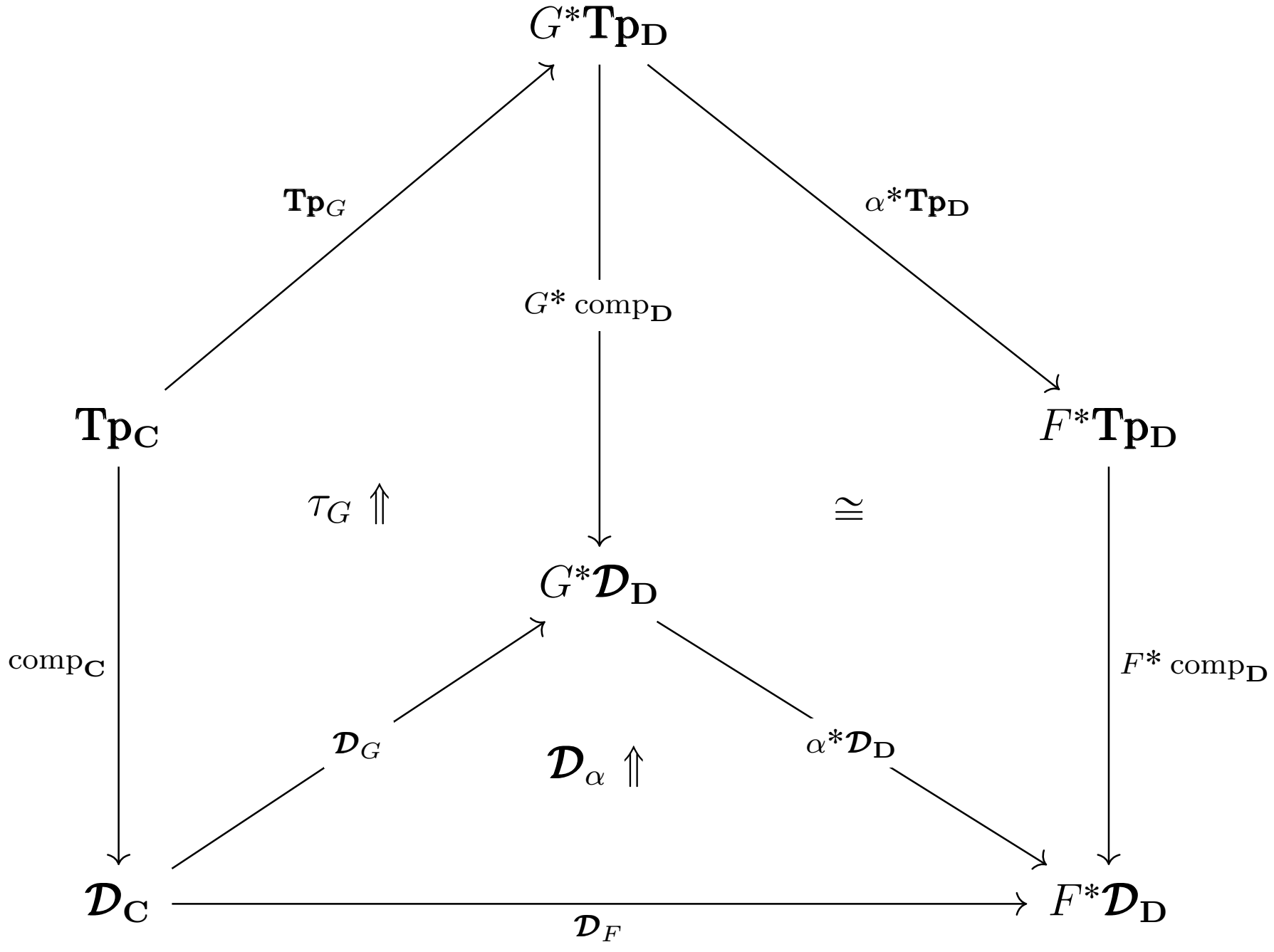

*in $\mathbf{IndCat}(\mathbf{C})$.*

When $\alpha : F \Rightarrow G : \mathbf{C} \to \mathbf{D}$ is a 2-morphism of small pre-CwAs, we denote by $\mathbf{Tm}_\alpha$ :

$\mathbf{Tm}_F \Rightarrow (\alpha^*\mathbf{Tm}_\mathbf{D} \circ \mathbf{Tm}_G)$ the unique $\mathbf{C}$-indexed natural transformation such that the pasting diagram

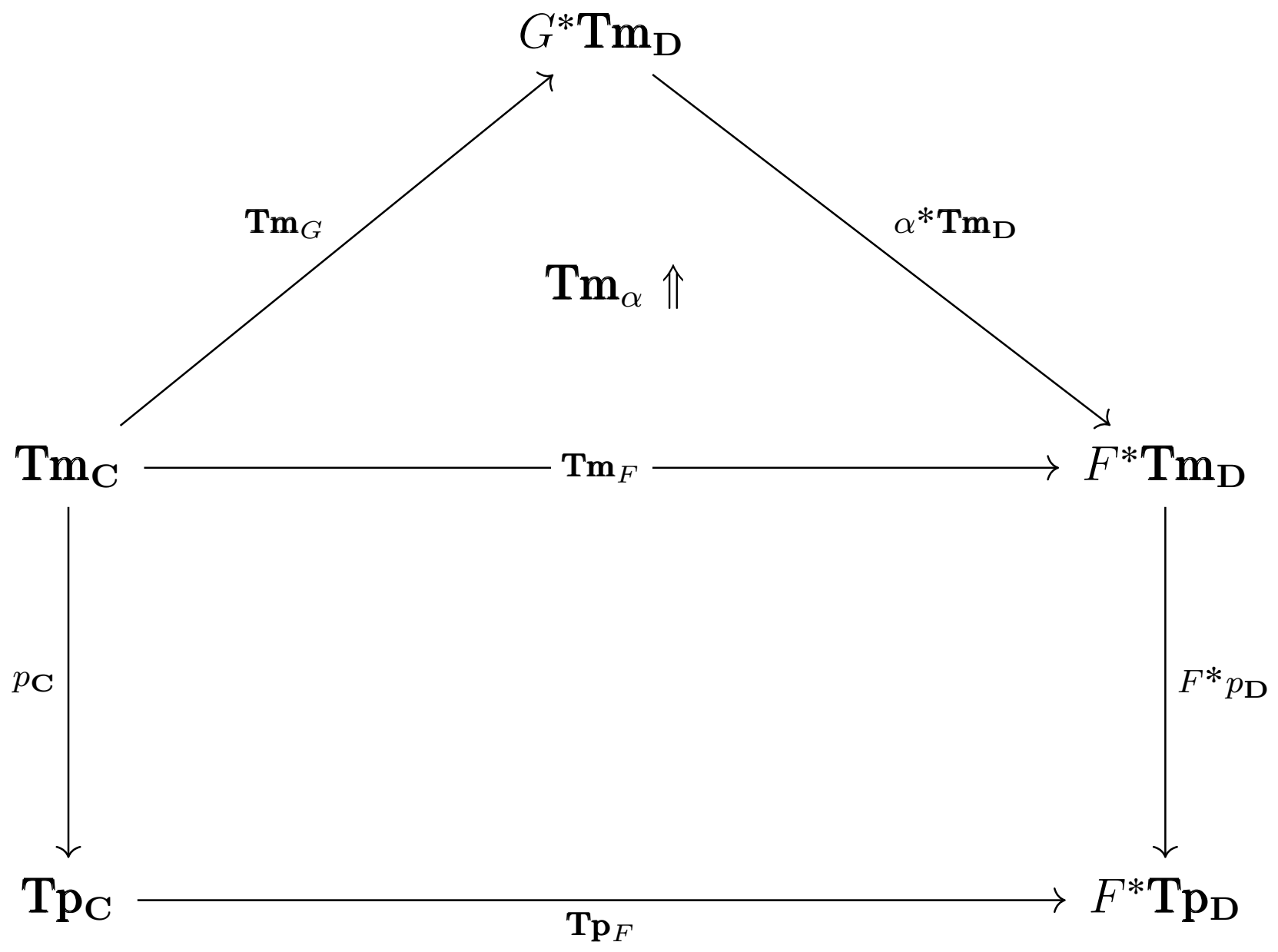

is equal to the pasting diagram

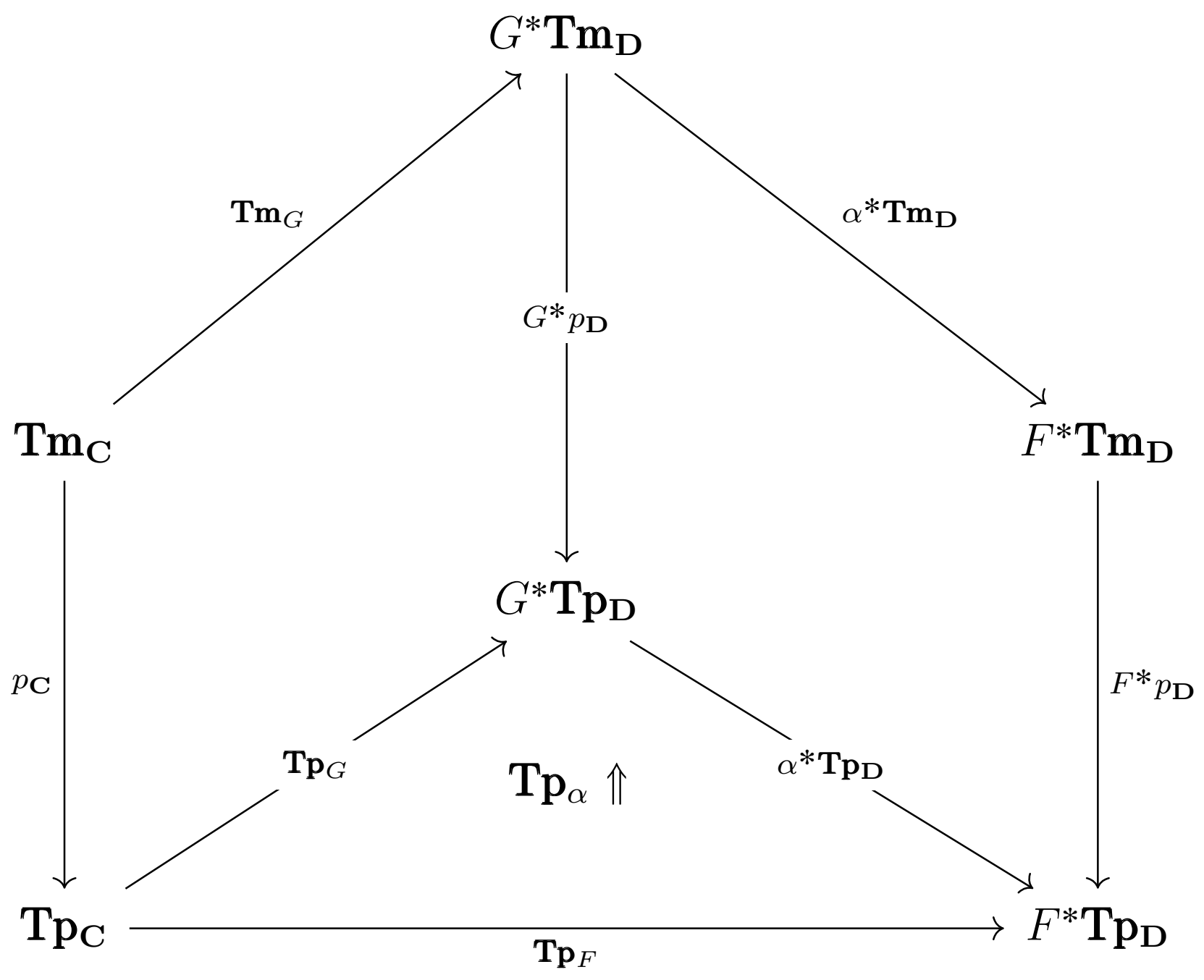

in $\mathbf{Cat}^{\mathbf{C}^{\mathrm{op}}}$, which as a matter of course is also seen to be such that the pasting diagram

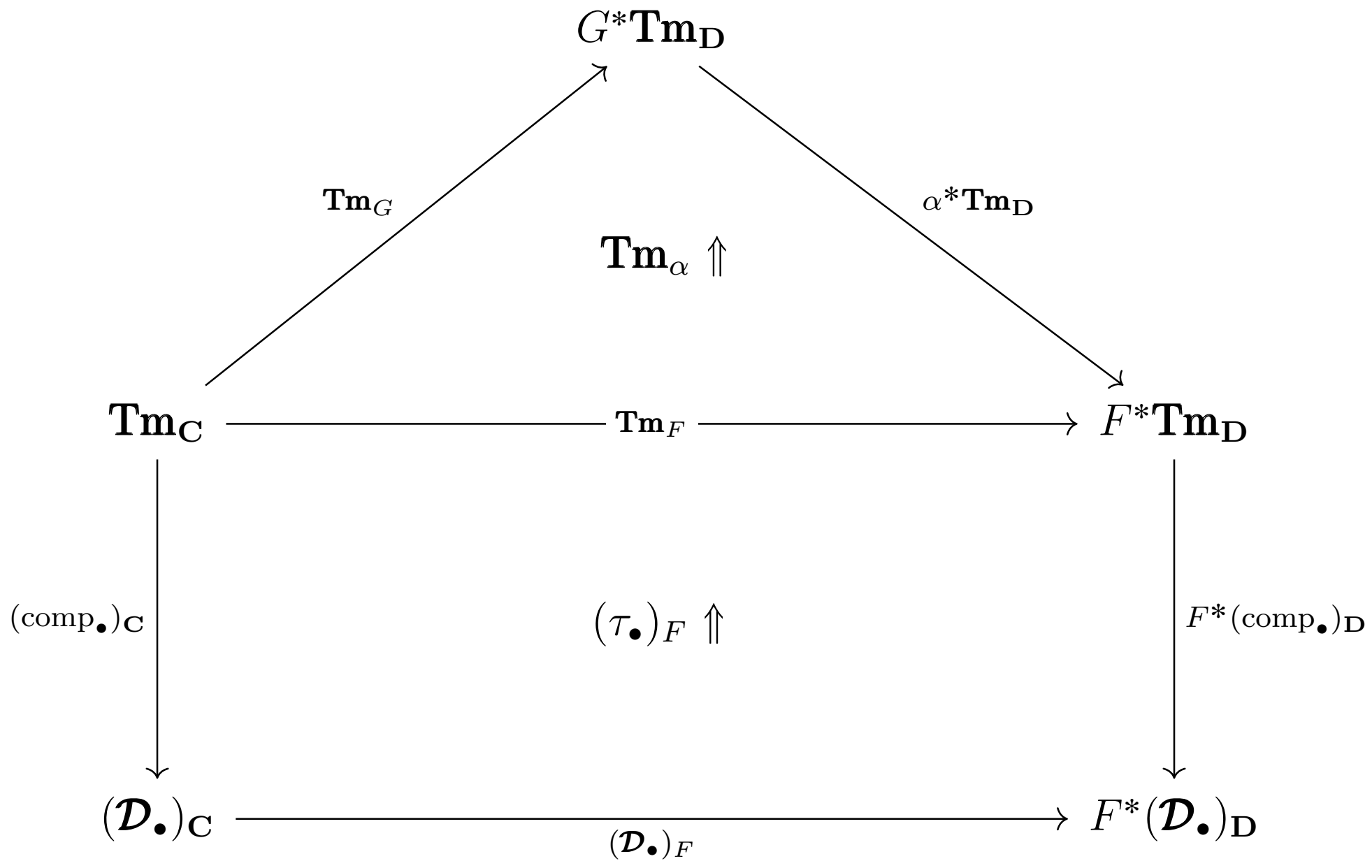

is equal to the pasting diagram

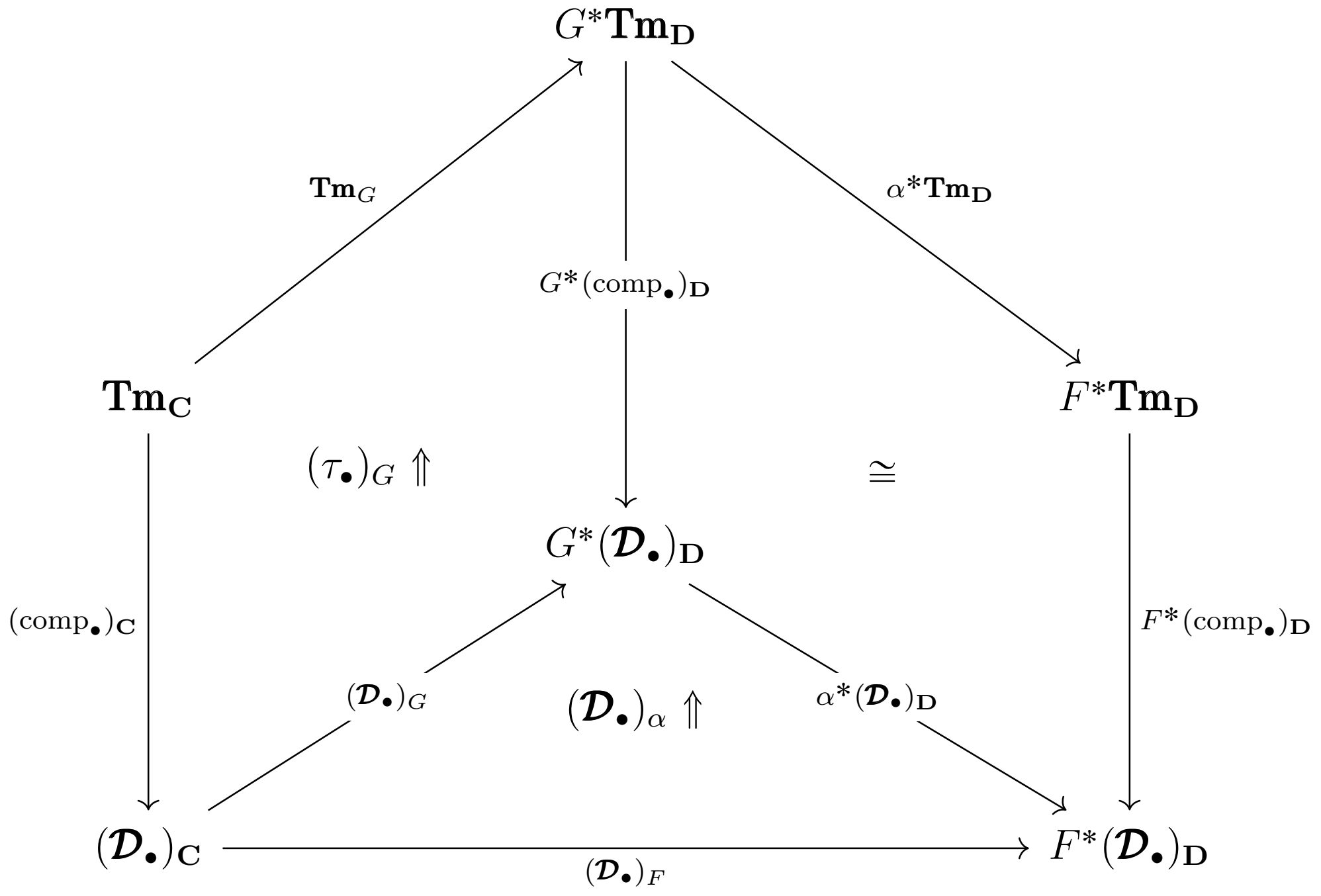

in $\mathbf{IndCat}(\mathbf{C})$.

**Definition 3.4.5.** *Let $\mathbf{C}$ and $\mathbf{D}$ be small CwAs. We say that a morphism $F : \mathbf{C} \to \mathbf{D}$ of small pre-CwAs is a **morphism of small CwAs** if the $\mathbf{C}$-indexed natural transformation*

$$\begin{array}{ccc} \mathbf{Tp}_{\mathbf{C}} & \xrightarrow{\mathbf{Tp}_F} & F^*\mathbf{Tp}_{\mathbf{D}} \\ {\scriptstyle \mathrm{comp}_{\mathbf{C}}}\big\downarrow & \tau_F \Uparrow & \big\downarrow{\scriptstyle F^*\,\mathrm{comp}_{\mathbf{D}}} \\ \boldsymbol{\mathcal{D}}_{\mathbf{C}} & \xrightarrow[\boldsymbol{\mathcal{D}}_F]{} & F^*\boldsymbol{\mathcal{D}}_{\mathbf{D}} \end{array}$$

*is a* $\mathbf{C}$*-indexed natural isomorphism.*

We denote by $\mathbf{PreCwA}$ the strict 2-category of small pre-CwAs and morphisms and 2-morphisms thereof. We denote by $\mathbf{CwA}$ the strict 2-subcategory of $\mathbf{PreCwA}$ spanned by the small CwAs and morphisms thereof.

### 3.4.3 The Equivalence of CwAs and Natural Models

We denote by $P : \mathbf{PreCwA} \to \mathbf{PreNM}$ the strict 2-functor that turns a small pre-CwA $\mathbf{C}$ into a small pre-natural model by replacing $\mathrm{comp}_{\mathbf{C}} : \mathbf{Tp}_{\mathbf{C}} \to \boldsymbol{\mathcal{D}}_{\mathbf{C}}$ by $p_{\mathbf{C}} : \mathbf{Tm}_{\mathbf{C}} \to \mathbf{Tp}_{\mathbf{C}}$, turns a morphism of small pre-CwAs $F : \mathbf{C} \to \mathbf{D}$ into a morphism of small pre-natural models by replacing

$$\begin{array}{ccc} \mathbf{Tp}_{\mathbf{C}} & \xrightarrow{\mathbf{Tp}_F} & F^*\mathbf{Tp}_{\mathbf{D}} \\ {\scriptstyle \mathrm{comp}_{\mathbf{C}}}\big\downarrow & \tau_F \Uparrow & \big\downarrow{\scriptstyle F^*\,\mathrm{comp}_{\mathbf{D}}} \\ \boldsymbol{\mathcal{D}}_{\mathbf{C}} & \xrightarrow[\boldsymbol{\mathcal{D}}_F]{} & F^*\boldsymbol{\mathcal{D}}_{\mathbf{D}} \end{array}$$

by $\mathbf{Tm}_F : \mathbf{Tm}_{\mathbf{C}} \to F^*\mathbf{Tm}_{\mathbf{D}}$, and turns a 2-morphism of small pre-CwAs $\alpha : F \Rightarrow G : \mathbf{C} \to \mathbf{D}$ into a 2-morphism of small pre-natural models by adducing $\mathbf{Tm}_\alpha : \mathbf{Tm}_F \Rightarrow (\alpha^*\mathbf{Tm}_{\mathbf{D}} \circ \mathbf{Tm}_G) : \mathbf{Tm}_{\mathbf{C}} \to F^*\mathbf{Tm}_{\mathbf{D}}$.

We denote by $C : \mathbf{PreNM} \to \mathbf{PreCwA}$ the strict 2-functor that turns a small pre-natural model $\mathbf{C}$ into a small pre-CwA by replacing $p_{\mathbf{C}} : \mathbf{Tm}_{\mathbf{C}} \to \mathbf{Tp}_{\mathbf{C}}$ by $\mathrm{comp}_{\mathbf{C}} : \mathbf{Tp}_{\mathbf{C}} \to \boldsymbol{\mathcal{D}}_{\mathbf{C}}$, turns a morphism of small pre-natural models $F : \mathbf{C} \to \mathbf{D}$ into a morphism of small pre-CwAs by replacing $\mathbf{Tm}_F : \mathbf{Tm}_{\mathbf{C}} \to F^*\mathbf{Tm}_{\mathbf{D}}$ by

$$\begin{array}{ccc} \mathbf{Tp}_{\mathbf{C}} & \xrightarrow{\mathbf{Tp}_F} & F^*\mathbf{Tp}_{\mathbf{D}} \\ {\scriptstyle \mathrm{comp}_{\mathbf{C}}}\big\downarrow & \tau_F \Uparrow & \big\downarrow{\scriptstyle F^*\,\mathrm{comp}_{\mathbf{D}}} \\ \boldsymbol{\mathcal{D}}_{\mathbf{C}} & \xrightarrow[\boldsymbol{\mathcal{D}}_F]{} & F^*\boldsymbol{\mathcal{D}}_{\mathbf{D}} \end{array},$$

and turns a 2-morphism of small pre-natural models $\alpha : F \Rightarrow G : \mathbf{C} \to \mathbf{D}$ into a 2-morphism of small pre-CwAs by forgetting $\mathbf{Tm}_\alpha : \mathbf{Tm}_F \Rightarrow (\alpha^*\mathbf{Tm}_{\mathbf{D}} \circ \mathbf{Tm}_G) : \mathbf{Tm}_{\mathbf{C}} \to F^*\mathbf{Tm}_{\mathbf{D}}$.

**Lemma 3.4.6.** *The strict* $2$*-functors* $P : \mathbf{PreCwA} \to \mathbf{PreNM}$ *and* $C : \mathbf{PreNM} \to \mathbf{PreCwA}$ *are inverse strict* $2$*-equivalences (over* $\mathbf{Cat}$*).*

*Proof.* It is evident that $PC \cong \mathrm{id}_{\mathbf{PreNM}}$. As for $\mathrm{id}_{\mathbf{PreCwA}} \cong CP$, we limit ourselves to verifying that the comprehension of a type $\Gamma \vdash_{\mathbf{C}} A$ in a small pre-CwA $\mathbf{C}$ agrees with that in $C(P(\mathbf{C}))$. Given such $\Gamma \vdash_{\mathbf{C}} A$, the comprehension in $\mathbf{C}$ is (isomorphic to) the strict 2-pullback

$$
\begin{array}{ccccc}
\cdot & & \longrightarrow & & (\boldsymbol{\mathcal{D}}_\bullet)_{\mathbf{C}} \\
\downarrow & \lrcorner & & & \downarrow{\scriptstyle \varpi_{\mathbf{C}}} \\
\Gamma & \xrightarrow[A]{} & \mathbf{Tp}_{\mathbf{C}} & \xrightarrow[\mathrm{comp}_{\mathbf{C}}]{} & \boldsymbol{\mathcal{D}}_{\mathbf{C}}
\end{array},
$$

whereas the comprehension in $C(P(\mathbf{C}))$ is (isomorphic to) the double strict 2-pullback

$$
\begin{array}{ccccc}
\cdot & \longrightarrow & \mathbf{Tm}_{\mathbf{C}} & \xrightarrow{(\mathrm{comp}_\bullet)_{\mathbf{C}}} & (\boldsymbol{\mathcal{D}}_\bullet)_{\mathbf{C}} \\
\downarrow & \lrcorner & \downarrow{\scriptstyle p_{\mathbf{C}}} & \lrcorner & \downarrow{\scriptstyle \varpi_{\mathbf{C}}} \\
\Gamma & \xrightarrow[A]{} & \mathbf{Tp}_{\mathbf{C}} & \xrightarrow[\mathrm{comp}_{\mathbf{C}}]{} & \boldsymbol{\mathcal{D}}_{\mathbf{C}}
\end{array}.
$$

But, these are isomorphic, by the two (strict 2-)pullbacks lemma. □

The following is evident:

**Theorem 3.4.7.** *The strict* $2$*-equivalence* $\mathbf{PreCwA} \simeq \mathbf{PreNM}$ *of Lemma 3.4.6 restricts to a strict* $2$*-equivalence* $\mathbf{CwA} \simeq \mathbf{NM}$.

### 3.4.4 CwA as an Equi-Comma Object

In Chapter 4, it will be useful to have a description of $\mathbf{CwA}$ as a Gray-limit.

We omit the proofs from this section, as the statements presented are, ultimately, tautologous observations.

The following Proposition 3.4.8 is presented purely for comparison with Theorem 3.4.9 (as Gray comma objects are a shade more familiar than equi-comma objects).

**Proposition 3.4.8.** $\mathbf{PreCwA}$ *arises as the Gray comma object indicated by the pasting diagram*

$$\begin{array}{ccc} \mathbf{PreCwA} & \longrightarrow & \mathbf{DC} \\ \big\downarrow & \text{comp} \Rightarrow & \big\downarrow \mathcal{D} \\ \mathbf{IndCat_s} & \hookrightarrow & \mathbf{IndCat} \end{array}$$

*in* $\mathbf{GRAY}_\ell/\mathbf{Cat}$.

**Theorem 3.4.9.** $\mathbf{CwA}$ *arises as the equi-comma object indicated by the pasting diagram*

$$\begin{array}{ccc} \mathbf{CwA} & \longrightarrow & \mathbf{DC} \\ \big\downarrow & \text{comp} \simeq & \big\downarrow \mathcal{D} \\ \mathbf{IndCat_s} & \hookrightarrow & \mathbf{IndCat} \end{array}$$

*in* $\mathbf{GRAY}/\mathbf{Cat}$.

# Chapter 4

# The CwA of Coalgebras

We construct the **CwA of coalgebras**, showing that $\mathbf{CwA}$, our strict 2-category of small CwAs (see Section 3.4), admits the construction of objects of coalgebras for comonads. Given the strict 2-equivalence of $\mathbf{CwA}$ and $\mathbf{NM}$, our strict 2-category of small natural models, (Theorem 3.4.7) we thereby obtain the **natural model of coalgebras** and that $\mathbf{NM}$ admits the construction of objects of coalgebras for comonads.

Our approach builds on results from Hermida's dissertation (Hermida, 1993).[1]

### Outline

Section 4.1 sketches the construction of the **indexed category of coalgebras**, following Hermida (1993). Section 4.2 gives the construction of the **display category of coalgebras**. Section 4.3 gives our construction of the CwA of coalgebras. Section 4.4 adapts this result to yield the natural model of coalgebras.

## 4.1 Indexed Categories of Coalgebras

As a preliminary for our construction of the CwA of coalgebras, we indicate here the construction of the **indexed category of coalgebras**, which was previously noted by Hermida (1993).

[1]Hermida worked in the setting of fibrations, and we will freely apply his results in the equivalent setting of indexed categories.

When

$$(\Box, \beta) : (\mathbf{C}, P) \to (\mathbf{C}, P)$$

is a comonad in $\mathbf{IndCat}$, we can define a comonad

$$\mathbb{B} : U^* P \to U^* P$$

in $\mathbf{IndCat}(\mathbf{C}^{\Box})$, as indicated by Definition 5.4.8 of Hermida (1993), which we term the **induced** (**indexed**) **comonad**. On $(\Gamma, \gamma) \in \mathbf{C}^{\Box}$ and $A \in P(\Gamma)$, this is given by

$$\mathbb{B}_{(\Gamma,\gamma)}(A) = P(\gamma)(\beta_{\Gamma}(A)) \quad .$$

Following Hermida (1993, Remark 5.4.12), it is then natural to define the indexed category

$$(\mathbf{C}^{\Box}, (U^* P)^{\mathbb{B}}) \quad ,$$

which we term term the (**small**) **indexed category of coalgebras**.

As also noted by Hermida, when the comonad $(\Box, \beta)$ lies in $\mathbf{IndCat_s}$, the induced comonad $\mathbb{B}$ lies in $\mathbf{Cat}^{(\mathbf{C}^{\Box})^{\mathrm{op}}}$, and the indexed category of coalgebras $(\mathbf{C}^{\Box}, (U^* P)^{\mathbb{B}})$ again lies in $\mathbf{IndCat_s}$.

The significance of the indexed category of coalgebras can be conveyed thus:

**Proposition 4.1.1.** *The strict* $2$*-category* $\mathbf{IndCat}$ *admits the construction of coalgebras over* $\mathbf{Cat}$.

*Proof.* Given a comonad $(\Box, \beta) : (\mathbf{C}, P) \to (\mathbf{C}, P)$ in $\mathbf{IndCat}$, the object of coalgebras is given by the indexed category of coalgebras $(\mathbf{C}^{\Box}, (U^* P)^{\mathbb{B}})$. □

Moreover, we have the following:

**Proposition 4.1.2.** *The strict* $2$*-category* $\mathbf{IndCat_s}$ *admits, and the forgetful strict* $2$*-functor* $\mathbf{IndCat_s} \hookrightarrow \mathbf{IndCat}$ *preserves (on the nose), the construction of coalgebras over* $\mathbf{Cat}$.

*Proof.* Given a comonad $(\Box, \beta) : (\mathbf{C}, P) \to (\mathbf{C}, P)$ in $\mathbf{IndCat_s}$, the object of coalgebras is given by the indexed category of coalgebras $(\mathbf{C}^{\Box}, (U^* P)^{\mathbb{B}})$. □

We also note at this stage a definition comparable to that of the induced indexed comonad, to be used in Chapter 7.

When $(L, \lambda) \dashv (R, \rho) : (\mathbf{C}, P) \to (\mathbf{D}, Q)$ is an adjunction in $\mathbf{IndCat}$, we can define a right adjoint

$$\mathbb{R} : L^* P \to Q$$

to

$$\lambda : Q \to L^* P$$

in $\mathbf{IndCat}(\mathbf{D})$, as indicated at Theorem 3.2.3 of Hermida (1993), which we term the **induced (indexed) right adjoint**. On $\Delta \in \mathbf{D}$ and $A \in P(L(\Delta))$, this is given by

$$\mathbb{R}_\Delta(A) = Q(\eta_\Delta)(\rho_{L\Delta}(A)) \quad .$$

We also term the adjunction $\lambda \dashv \mathbb{R}$ the **induced (indexed) adjunction**.

As also noted by Hermida, when the adjunction $(L, \lambda) \dashv (R, \rho)$ lies in $\mathbf{IndCat_s}$, the induced adjunction $\lambda \dashv \mathbb{R}$ lies in $\mathbf{Cat}^{\mathbf{D}^{\mathrm{op}}}$.

## 4.2 The Display Category of Coalgebras

We give (and flesh out slightly) the construction of the **display category of coalgebras**.

**Definition 4.2.1.** *When $\mathbf{C}$ is a display category, we will say that a comonad (of categories) $\Box : \mathbf{C} \to \mathbf{C}$ is **display** if its underlying functor is display.*

A display comonad on a small display category is then equivalent to a comonad in $\mathbf{DC}$.

We will say that a display comonad is **Cartesian** if its underlying display functor is Cartesian.

When $\Box : \mathbf{C} \to \mathbf{C}$ is a display comonad, we denote by $\mathbf{C}^\Box$ the display category with eponymous underlying category and class of displays $\mathcal{D}_{\mathbf{C}^\Box}$ given by

$$p \in \mathcal{D}_{\mathbf{C}^\Box} \Leftrightarrow Up \in \mathcal{D}_{\mathbf{C}} \quad .$$

We term $\mathbf{C}^\Box$ the **display category of coalgebras**. That this $\mathbf{C}^\Box$ is indeed a display category follows by application of Borceux (1994, Proposition 4.3.2).[2] *Cf.* also Warren (2007, Proposition 7).

The significance of the display category of coalgebras can be conveyed thus:

**Proposition 4.2.2.** *The strict* $2$*-category* $\mathbf{DC}$ *admits the construction of coalgebras over* $\mathbf{Cat}$.

*Proof.* Given a comonad $\Box : \mathbf{C} \to \mathbf{C}$ in $\mathbf{DC}$, *i.e.* a display comonad on a small display category, the object of coalgebras is given by the display category of coalgebras $\mathbf{C}^{\Box}$. $\Box$

When

$$\Box : \mathbf{C} \to \mathbf{C}$$

is a display comonad on a small display category, then

$$(\Box, \boldsymbol{\mathcal{D}}_{\Box}) : (\mathbf{C}, \boldsymbol{\mathcal{D}}_{\mathbf{C}}) \to (\mathbf{C}, \boldsymbol{\mathcal{D}}_{\mathbf{C}})$$

is a comonad in $\mathbf{IndCat}$, where

$$\boldsymbol{\mathcal{D}}_{\Box} : \boldsymbol{\mathcal{D}}_{\mathbf{C}} \to \Box^{*}\boldsymbol{\mathcal{D}}_{\mathbf{C}}$$

is the induced $\mathbf{C}$-indexed functor (Section 3.1.2). Consequently, we can form the induced $\mathbf{C}^{\Box}$-indexed comonad

$$\mathbb{B} : U^{*}\boldsymbol{\mathcal{D}}_{\mathbf{C}} \to U^{*}\boldsymbol{\mathcal{D}}_{\mathbf{C}}$$

and the indexed category of coalgebras

$$(\mathbf{C}^{\Box}, (U^{*}\boldsymbol{\mathcal{D}}_{\mathbf{C}})^{\mathbb{B}}) \quad .$$

It is then natural to ask for the relation of this indexed category of coalgebras to

$$(\mathbf{C}^{\Box}, \boldsymbol{\mathcal{D}}_{\mathbf{C}^{\Box}}) \quad ,$$

the indexed category corresponding to the display category of coalgebras $\mathbf{C}^{\Box}$. In fact, the relation is one of isomorphism:

[2]The statement in Borceux is about monads, so we dualize appropriately.

**Proposition 4.2.3.** *The strict* $2$*-functor* $\boldsymbol{\mathcal{D}} : \mathbf{DC} \to \mathbf{IndCat}$ *strictly preserves the construction of coalgebras over* $\mathbf{Cat}$.

*Proof.* Given a $\Box$-coalgebra $(B, \beta) \in \mathbf{C}^{\Box}$ and display map $p : E \to B$ of $\mathbf{C}$, we observe that $\Box$-coalgebra structures on $E$ making $p$ into a morphism of $\Box$-coalgebras (equivalently objects of $\boldsymbol{\mathcal{D}}_{\mathbf{C}^{\Box}}(B, \beta)$ laying over $p$) are in correspondence with $\mathbb{B}_{(B,\beta)}$-coalgebra structures on $p$ (*i.e.* objects of $(U^{*}\boldsymbol{\mathcal{D}}_{\mathbf{C}})^{\mathbb{B}}(B, \beta)$ laying over $p$), mediated by the diagram

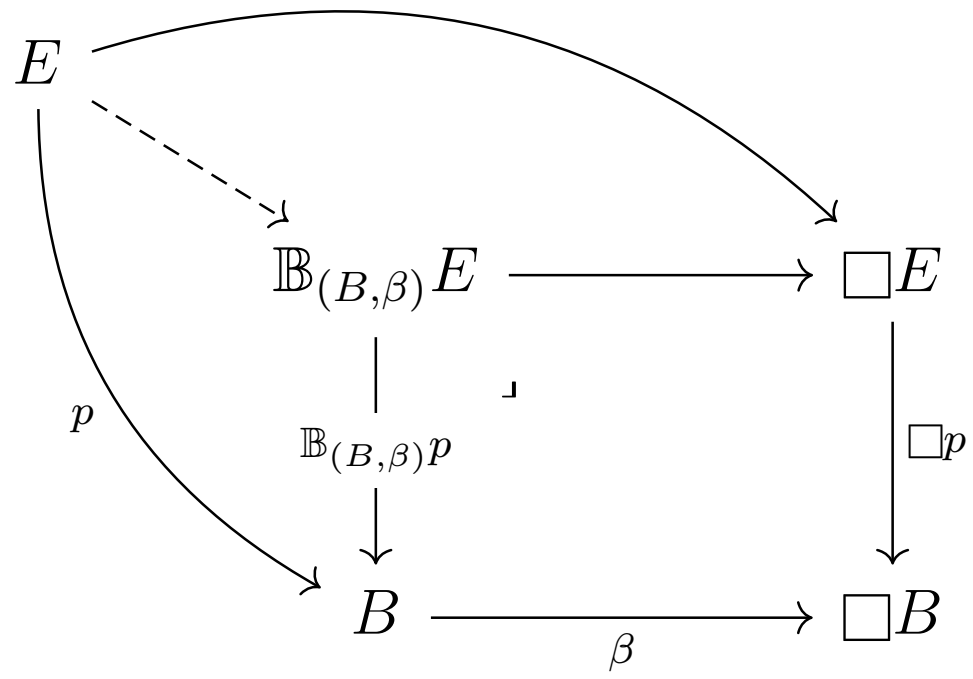

in $\mathbf{C}$. $\Box$

Finally, we note the following for future reference:

**Proposition 4.2.4.** *Let* $\mathbf{C}$ *be a small, Cartesian display category and* $\Box : \mathbf{C} \to \mathbf{C}$ *a Cartesian display comonad. Then, the display category of coalgebras* $\mathbf{C}^{\Box}$ *is again Cartesian.*

*Proof.* Since the underlying category of $\mathbf{C}$ and the underlying functor of $\Box : \mathbf{C} \to \mathbf{C}$ are Cartesian, so is the underlying category of $\mathbf{C}^{\Box}$.

As observed in the proof of Proposition 4.2.3, the $\mathbf{C}^{\Box}$-indexed functor

$$\boldsymbol{\mathcal{D}}_{U} : \boldsymbol{\mathcal{D}}_{\mathbf{C}^{\Box}} \to U^{*}\boldsymbol{\mathcal{D}}_{\mathbf{C}}$$

induced by the display functor $U : \mathbf{C}^{\Box} \to \mathbf{C}$ arises up to isomorphism as the forgetful morphism

$$(U^{*}\boldsymbol{\mathcal{D}}_{\mathbf{C}})^{\mathbb{B}} \to U^{*}\boldsymbol{\mathcal{D}}_{\mathbf{C}}$$

from the object of coalgebras $(U^{*}\boldsymbol{\mathcal{D}}_{\mathbf{C}})^{\mathbb{B}}$ for the $\mathbf{C}^{\Box}$-indexed comonad

$$\mathbb{B} : U^{*}\boldsymbol{\mathcal{D}}_{\mathbf{C}} \to U^{*}\boldsymbol{\mathcal{D}}_{\mathbf{C}}$$

induced by the comonad

$$(\Box, \boldsymbol{\mathcal{D}}_{\Box}) : (\mathrm{C}, \boldsymbol{\mathcal{D}}_{\mathbf{C}}) \to (\mathrm{C}, \boldsymbol{\mathcal{D}}_{\mathbf{C}})$$

in $\mathbf{IndCat}$.

Thus, at any coalgebra $B \in \mathbf{C}^{\Box}$, the forgetful functor

$$(\boldsymbol{\mathcal{D}}_U)_B : \boldsymbol{\mathcal{D}}_{\mathbf{C}^{\Box}}(B) \to \boldsymbol{\mathcal{D}}_{\mathbf{C}}(UB)$$

arises up to isomorphism as the forgetful functor

$$(\boldsymbol{\mathcal{D}}_{\mathbf{C}}(UB))^{\mathbb{B}_B} \to \boldsymbol{\mathcal{D}}_{\mathbf{C}}(UB)$$

from coalgebras for the comonad

$$\mathbb{B}_B : \boldsymbol{\mathcal{D}}_{\mathbf{C}}(UB) \to \boldsymbol{\mathcal{D}}_{\mathbf{C}}(UB) \quad .$$

Since $\boldsymbol{\mathcal{D}}_{\mathbf{C}}(UB)$ and $\mathbb{B}_B$ are Cartesian, so are $(\boldsymbol{\mathcal{D}}_{\mathbf{C}}(UB))^{\mathbb{B}_B}$ and the isomorphic $\boldsymbol{\mathcal{D}}_{\mathbf{C}^{\Box}}(B)$.

It is tautologous enough to extend this observation to obtain the remaining conditions. □

## 4.3 The CwA of Coalgebras

We construct the **CwA of coalgebras**.

**Theorem 4.3.1.** *The strict* 2*-category* $\mathbf{CwA}$ *admits the construction of coalgebras over* $\mathbf{Cat}$.

*Proof.* As an instance of Corollary 2.3.7, the (very large) Gray-category of strict 2-categories that admit the construction of coalgebras over $\mathbf{Cat}$ and strict 2-functors that preserve the construction of coalgebras over $\mathbf{Cat}$ is closed under the construction of equi-comma objects in $\mathbf{GRAY}/\mathbf{Cat}$. Due to the description of the forgetful strict 2-functor $\mathbf{CwA} \to \mathbf{Cat}$ as an equi-comma object in $\mathbf{GRAY}/\mathbf{Cat}$ (Theorem 3.4.9), our statement then follows from Propositions 4.1.2 and 4.2.3. □

When $\Box : \mathbf{C} \to \mathbf{C}$ is a comonad of small CwAs, we term the object of coalgebras $\mathbf{C}^{\Box} \in$

$\mathbf{CwA}$ derived from Theorem 4.3.1 the (**small**) **CwA of coalgebras**.

By unwinding the construction of the CwA of coalgebras, one sees that

$$\mathbf{Tp}_{\mathbf{C}^{\Box}} \equiv (U^{*}\mathbf{Tp}_{\mathbf{C}})^{\mathbb{B}} \quad ,$$

where

$$\mathbb{B} : U^{*}\mathbf{Tp}_{\mathbf{C}} \to U^{*}\mathbf{Tp}_{\mathbf{C}}$$

is the comonad in $\mathbf{Cat}^{(\mathbf{C}^{\Box})^{\mathrm{op}}}$ induced by the comonad

$$(\Box, \mathbf{Tp}_{\Box}) : (\mathbf{C}, \mathbf{Tp}_{\mathbf{C}}) \to (\mathbf{C}, \mathbf{Tp}_{\mathbf{C}})$$

in $\mathbf{IndCat}_{\mathbf{s}}$.

## 4.4 The Natural Model of Coalgebras

In light of Theorem 4.3.1 and the strict 2-equivalence $\mathbf{CwA} \simeq \mathrm{NM}$ over $\mathbf{Cat}$, we observe the following:

**Corollary 4.4.1.** *The strict* 2*-category* NM *admits the construction of coalgebras over* $\mathbf{Cat}$.

When $\Box : \mathbf{C} \to \mathbf{C}$ is a comonad of small natural models, we term the object $\mathbf{C}^{\Box}$ of coalgebras in NM the (**small**) **natural model of coalgebras**.

As with the CwA of coalgebras, by unwinding the construction of the natural model of coalgebras, one sees that

$$\mathbf{Tp}_{\mathbf{C}^{\Box}} = (U^{*}\mathbf{Tp}_{\mathbf{C}})^{\mathbb{B}} \quad ,$$

where

$$\mathbb{B} : U^{*}\mathbf{Tp}_{\mathbf{C}} \to U^{*}\mathbf{Tp}_{\mathbf{C}}$$

is the comonad in $\mathbf{Cat}^{(\mathbf{C}^{\Box})^{\mathrm{op}}}$ induced by the comonad

$$(\Box, \mathbf{Tp}_{\Box}) : (\mathbf{C}, \mathbf{Tp}_{\mathbf{C}}) \to (\mathbf{C}, \mathbf{Tp}_{\mathbf{C}})$$

in $\mathbf{IndCat_s}$.

For use in Chapter 7, we will denote by

$$\mathbb{B} : U^*\mathbf{Tm_C} \to U^*\mathbf{Tm_C}$$

the comonad in $\mathbf{Cat}^{(\mathbf{C}^\square)^{\mathrm{op}}}$ induced by the comonad

$$(\square, \mathbf{Tm}_\square) : (\mathbf{C}, \mathbf{Tm_C}) \to (\mathbf{C}, \mathbf{Tm_C})$$

in $\mathbf{IndCat_s}$, and by

$$\mathbb{F} : U^*\mathbf{Tp_C} \to \mathbf{Tp}_{\mathbf{C}^\square}$$

in $\mathbf{Cat}^{(\mathbf{C}^\square)^{\mathrm{op}}}$ the induced right adjoint corresponding to the adjunction

$$(U, \mathbf{Tp}_U) \dashv (F, \mathbf{Tp}_F) : (\mathbf{C}, \mathbf{Tp_C}) \to (\mathbf{C}^\square, \mathbf{Tp}_{\mathbf{C}^\square})$$

in $\mathbf{IndCat_s}$.

We must also refine our results to Cartesian natural models for future reference.

When $\mathbf{C}$ is a Cartesian small natural model, we say that a comonad of small natural models $\square : \mathbf{C} \to \mathbf{C}$ is **Cartesian** if the underlying morphism of small natural models is Cartesian.

**Proposition 4.4.2.** *Let $\mathbf{C}$ be a Cartesian small natural model and $\square : \mathbf{C} \to \mathbf{C}$ be a Cartesian comonad of small natural models. Then, the small natural model of coalgebras $\mathbf{C}^\square$ is again Cartesian.*

*Proof.* From Proposition 4.2.4. $\square$

# Chapter 5

# The Natural Display Topos of Coalgebras

One of the classical constructions of topos theory shows that, given a topos $\mathcal{E}$ and a Cartesian comonad $\Box : \mathcal{E} \to \mathcal{E}$, the category of coalgebras $\mathcal{E}^{\Box}$ is again a topos, the so-called topos of coalgebras (Kock and Wraith, 1971).

In this chapter, we show (modulo a size limitation) that, given a natural display topos $\mathcal{E}$ and a natural Cartesian display comonad $\Box : \mathcal{E} \to \mathcal{E}$, the natural model of coalgebras $\mathcal{E}^{\Box}$ (see Chapter 4) is again a natural display topos, which we term the **natural display topos of coalgebras**. We also show (modulo the same size limitation) that, given a natural typos $\mathcal{E}$ and a natural Cartesian display comonad $\Box : \mathcal{E} \to \mathcal{E}$, the natural model of coalgebras $\mathcal{E}^{\Box}$ is again a natural typos, which we term the **natural typos of coalgebras**.

We thus refine the classical result to a context involving universes (or, more precisely, type classifiers), leading to various applications. Many natural models of interest can be constructed as natural typoses or natural display toposes of coalgebras. In Chapter 6, we show that many sheaf toposes can be constructed in this way, yielding an approach to universes in sheaf toposes that extends the standard approach to presheaf toposes (Hofmann and Streicher, 1997). The natural typos of coalgebras is also implicated in the semantics of S4 DTT, as indicated in Chapter 7.

As motivation for the upcoming development, let us recall some of the classical construction (for details, see, *e.g.*, Mac Lane and Moerdijk 1992, §V.8.). Let $\Box : \mathcal{E} \to \mathcal{E}$ be a Cartesian comonad on a topos $\mathcal{E}$. As with any comonad on a category, we can construct the category of coalgebras $\mathcal{E}^{\Box}$ and the forgetful-cofree decomposition $U \dashv F : \mathcal{E} \to \mathcal{E}^{\Box}$. We may verify that $\mathcal{E}^{\Box}$ admits a subobject classifier by constructing one from a subobject classifier $\mathbf{\Omega}_{\mathcal{E}}$ of $\mathcal{E}$. Curiously, the comonad

$$\Box : \mathcal{E} \to \mathcal{E} \tag{5.1}$$

induces a comonad

$$\beta : F\mathbf{\Omega}_{\mathcal{E}} \to F\mathbf{\Omega}_{\mathcal{E}} \tag{5.2}$$

of posets *in* $\mathcal{E}^{\Box}$. A subobject classifier $\mathbf{\Omega}_{\mathbf{C}^{\Box}}$ for $\mathbf{C}^{\Box}$ is then constructed as the internal poset of coalgebras (or, equivalently, the internal poset of fixed points) $(F\mathbf{\Omega}_{\mathcal{E}})^{\beta}$.

We will undertake a roughly analogous construction, in which we work with internal categories rather than internal posets. We will let $\Box : \mathcal{E} \to \mathcal{E}$ be a natural Cartesian display comonad on a small natural display topos $\mathcal{E}$. As with any comonad of small natural models, we can construct the small natural model of coalgebras $\mathcal{E}^{\Box}$ and the forgetful-cofree decomposition $U \dashv F : \mathcal{E} \to \mathcal{E}^{\Box}$ (Corollary 4.4.1). We will verify that $\mathcal{E}^{\Box}$ admits a type classifier by constructing one from a type classifier $\mathbf{\mho}_{\mathcal{E}}$ of $\mathcal{E}$. Now, as in (5.1), the comonad

$$\Box : \mathcal{E} \to \mathcal{E}$$

will induce, as in (5.2), a comonad

$$\beta : F\mathbf{\mho}_{\mathcal{E}} \to F\mathbf{\mho}_{\mathcal{E}}$$

of categories *in* $\mathcal{E}^{\Box}$. A type classifier $\mathbf{\mho}_{\mathcal{E}^{\Box}}$ for $\mathcal{E}^{\Box}$ will then be constructed as the internal category of coalgebras $(F\mathbf{\mho}_{\mathcal{E}})^{\beta}$.

**Outline** In Section 5.1, we construct the small natural display topos and small natural typos of coalgebras. In Section 5.2, we use these to construct the small display topos and small typos of coalgebras.

## 5.1 Main Theorem

**Lemma 5.1.1.** *Let* $\mathbf{C}$ *and* $\mathbf{D}$ *be Cartesian categories and* $U : \mathbf{D} \to \mathbf{C}$ *a Cartesian, comonadic functor. Then, if* $\mathbf{C}$ *is Cartesian closed, so is* $\mathbf{D}$.

*Proof.* See Theorem A4.2.1(i) of Johnstone (2002). □

**Lemma 5.1.2.** *Let* $\mathbf{C}$ *be a small, Cartesian display category and* $\Box : \mathbf{C} \to \mathbf{C}$ *a Cartesian display comonad. Then, if* $\mathbf{C}$ *is Cartesian closed, then so is* $\mathbf{C}^{\Box}$.

*Proof.* Since the underlying category of $\mathbf{C}$ is Cartesian closed, so is that of $\mathbf{C}^{\Box}$ (Lemma 5.1.1).

As observed in the proof of Proposition 4.2.4, at any coalgebra $B \in \mathbf{C}^{\Box}$, the forgetful functor

$$\mathcal{D}_{\mathbf{C}^{\Box}}(B) \to \mathcal{D}_{\mathbf{C}}(UB)$$

is comonadic and Cartesian. Thus, since $\mathcal{D}_{\mathbf{C}}(UB)$ is Cartesian closed, so is $\mathcal{D}_{\mathbf{C}^{\Box}}(B)$ (Lemma 5.1.1 again).

Since the construction of exponentials in $\mathcal{D}_{\mathbf{C}^{\Box}}(B)$ uses the Cartesian closed structure of $\mathcal{D}_{\mathbf{C}}(UB)$, which is preserved by the inclusion

$$\mathcal{D}_{\mathbf{C}}(UB) \hookrightarrow \mathbf{C}/UB \quad ,$$

and the Cartesian structure of $\mathcal{D}_{\mathbf{C}^{\Box}}(B)$, which is preserved by the inclusion

$$\mathcal{D}_{\mathbf{C}^{\Box}}(B) \hookrightarrow \mathbf{C}^{\Box}/B \quad ,$$

it is preserved the inclusion

$$\mathcal{D}_{\mathbf{C}^{\Box}}(B) \hookrightarrow \mathbf{C}^{\Box}/B \quad .$$

The exponentials in $\mathcal{D}_{\mathbf{C}^{\Box}}(B)$ are similarly stable under reindexing because the Cartesian closed structure in $\mathcal{D}_{\mathbf{C}}(UB)$ and the Cartesian structure in $\mathcal{D}_{\mathbf{C}^{\Box}}(B)$ are. □

**Theorem 5.1.3.** *Let* $\mathbf{C}$ *be a Cartesian small natural model and* $\Box : \mathbf{C} \to \mathbf{C}$ *a Cartesian comonad of small natural models.*

1. *If* $\mathbf{C}$ *is Cartesian closed, then so is* $\mathbf{C}^{\Box}$.
2. *If* $\mathbf{C}$ *admits a type classifier, then so does* $\mathbf{C}^{\Box}$.
3. *If* $\mathbf{C}$ *is a small natural typos, then so is* $\mathbf{C}^{\Box}$.
4. *If* $\mathbf{C}$ *is a small natural display topos, then so is* $\mathbf{C}^{\Box}$.

*Proof.* 1. From Lemma 5.1.2.

2. Let $\mathfrak{V}_{\mathbf{C}} \in \mathbf{Cart}(\mathbf{C})$ be a type classifier. We have in $\mathbf{Cart}^{(\mathbf{C}^{\Box})^{\mathrm{op}}}$ the chain of isomorphisms

$$\begin{aligned} F\mathfrak{V}_{\mathbf{C}} &\cong F_!\mathfrak{V}_{\mathbf{C}} \\ &\cong U^*\mathfrak{V}_{\mathbf{C}} \\ &\cong U^*\mathbf{Tp}_{\mathbf{C}} \quad . \end{aligned}$$

The comonad $\mathbb{B} : U^*\mathbf{Tp}_{\mathbf{C}} \to U^*\mathbf{Tp}_{\mathbf{C}} \in \mathbf{Cart}^{(\mathbf{C}^{\Box})^{\mathrm{op}}}$ induced by $\Box : \mathbf{C} \to \mathbf{C}$ thus internalizes to a comonad which we denote by $\beta : F\mathfrak{V}_{\mathbf{C}} \to F\mathfrak{V}_{\mathbf{C}} \in \mathbf{Cart}(\mathbf{C}^{\Box})$. Using the Cartesian structure of $\mathbf{C}^{\Box}$, we define a type classifier $\mathfrak{V}_{\mathbf{C}^{\Box}} \in \mathbf{Cart}(\mathbf{C}^{\Box})$ by setting $\mathfrak{V}_{\mathbf{C}^{\Box}} :\equiv (F\mathfrak{V}_{\mathbf{C}})^{\beta}$. We then verify

$$\begin{aligned} \mathfrak{V}_{\mathbf{C}^{\Box}} &\equiv (F\mathfrak{V}_{\mathbf{C}})^{\beta} \\ &\cong (U^*\mathbf{Tp}_{\mathbf{C}})^{\mathbb{B}} \\ &= \mathbf{Tp}_{\mathbf{C}^{\Box}} \quad . \end{aligned}$$

3. The set of displays $\mathcal{D}_{\mathbf{C}^{\Box}}$ inherits closure under composition from $\mathcal{D}_{\mathbf{C}}$.
4. Since $U : \mathbf{C}^{\Box} \to \mathbf{C}$ is Cartesian, it takes monos to monos. Consequently, all monos in $\mathbf{C}^{\Box}$ are displays because those in $\mathbf{C}$ are. Moreover, the morphism $\Omega_{\mathbf{C}^{\Box}} \to 1_{\mathbf{C}^{\Box}}$ is display, as it arises as the composite of displays

$$\Omega_{\mathbf{C}^{\Box}} \rightarrowtail F\Omega_{\mathbf{C}} \to F1_{\mathbf{C}} \cong 1_{\mathbf{C}^{\Box}} \quad .$$

□

In the case where $\mathcal{E}$ is a small natural typos and $\Box : \mathcal{E} \to \mathcal{E}$ a Cartesian comonad of small natural models, we may also say that $\square$ is a **natural Cartesian display comonad**.

## 5.2 Corollary

**Corollary 5.2.1.** *Let* $\mathbf{C}$ *be a small, Cartesian display category and* $\square : \mathbf{C} \to \mathbf{C}$ *a Cartesian display comonad.*

1. *If* $\mathbf{C}$ *admits a universe, then so does* $\mathbf{C}^{\square}$*.*
2. *If* $\mathbf{C}$ *is a small typos, then so is* $\mathbf{C}^{\square}$*.*
3. *If* $\mathbf{C}$ *is a small display topos, then so is* $\mathbf{C}^{\square}$*.*

*Proof.*

1. Under the 2-equivalence of $\mathbf{CDC}_{\mho}$ and $\mathbf{CNM}_{\mho}$ (Corollary 3.2.14), we may transfer the comonad in the former to a comonad in the latter. This proceeds by making necessary choices for the extra structure. We then apply Theorem 5.1.3, obtaining a Cartesian small natural model that admits a type classifier. We make the return journey to $\mathbf{CDC}_{\mho}$, leaving an underlying small, Cartesian display category that admits a universe, which is precisely the display category of coalgebras.
2. The prior argument refines to the 2-equivalence of small typoses and small natural typoses (again at Corollary 3.2.14).
3. The prior argument refines to the 2-equivalence of small display toposes and small natural display toposes (again at Corollary 3.2.14).

$\square$

# Chapter 6

# Application: Type Classifiers for Sheaves

The type classifier in the natural model of presheaves on a category has been elucidated, *mutato mutandis*, by Hofmann and Streicher (1997). Type classifiers are problematic in the more general case of sheaves on a site, however. Indeed, the natural model of sheaves on a site usually lacks a type classifier. Candidates for the type classifier will, in general, form a *stack*, rather than a sheaf, in which amalgamations for matching families are unique only up to isomorphism. In the case of sheaves on a topological space, the underlying observation was already made by Grothendieck (1960, §3.3).

In this chapter, we generalize the presheaf case in an *a priori* different direction, introducing the notion of **Hofmann-Streicher natural display topos**, a particular kind of natural display topos of coalgebras (see Theorem 5.1.3). It turns out that any natural model of sheaves on a site with **enough points** can be replaced by an equivalent Hofmann-Streicher natural display topos (which, in particular, admits a type classifier). This covers many cases of interest; for example, any natural model sheaves on a topological space has enough points.

## 6.1 Preliminaries on Geometric Morphisms

In topos theory, when we say that a sheaf topos has enough points, we implicitly use the notion of geometric morphism of toposes, or, more generally, geometric morphism of Cartesian categories (*cf.* Streicher, 1997). These notions of geometric morphism can be handily refined to our set-ting, yielding **geometric morphisms of natural display toposes**, or, more generally, **geometric morphisms of Cartesian natural models**. We discuss such geometric morphisms briefly now.

**Definition 6.1.1.** *When $\mathcal{E}$ and $\mathcal{F}$ are small, Cartesian natural models, a* ***geometric morphism of small, Cartesian natural models*** $f : \mathcal{E} \to \mathcal{F}$ *will consist of an adjunction* $f^* \dashv f_* : \mathcal{E} \to \mathcal{F}$ *in* $\mathrm{CNM}$*. The morphism* $f_*$ *of* $\mathrm{CNM}$ *is called the* ***direct image****, while* $f^*$ *is called the* ***inverse image****. When $\mathcal{E}$ and $\mathcal{F}$ are small natural display toposes and* $f : \mathcal{E} \to \mathcal{F}$ *is a geometric morphism of small, Cartesian natural models, we will also say that* $f$ *is a* ***geometric morphism of small natural display toposes****.*

**Definition 6.1.2.** *We will say that a geometric morphism of small, Cartesian natural models* $p : \mathcal{E} \to \mathcal{F}$ *is a* ***surjection****, and write* $p : \mathcal{E} \twoheadrightarrow \mathcal{F}$*, if the underlying functor of* $p^*$ *is faithful.*

The following proposition gives alternate characterizations of such geometric surjections in the topos-theoretic case:

**Proposition 6.1.3.** *Let* $p : \mathcal{E} \to \mathcal{F}$ *be a geometric morphism of small, Cartesian natural models, such that the underlying display categories of $\mathcal{E}$ and $\mathcal{F}$ are display toposes. Then, the following are equivalent:*

1. $p$ *is a surjection;*
2. *the natural Cartesian display functor* $p^*$ *is comonadic;*
3. *the underlying functor* $p^*$ *is comonadic.*

In the statement of Proposition 6.1.3, we do not require that $\mathcal{E}$ or $\mathcal{F}$ are natural display toposes, to cover cases where these are natural models of sheaves that lack a type classifier.

*Proof.* The equivalence of (1) and (3) is due to the usual equivalence of surjective geometric morphisms of toposes and Cartesian comonadic functors between toposes (for which, see Mac Lane and Moerdijk, 1992).

That (2) implies (3) is easy. Finally, assume (3). Since the functor $p^* : \mathcal{F} \to \mathcal{E}$ is comonadic, it is faithful. For all $B \in \mathcal{F}$, the induced left adjoint $\boldsymbol{\mathcal{D}}_{\mathcal{F}}(B) \to \boldsymbol{\mathcal{D}}_{\mathcal{E}}(p^*B)$ is then clearly faithful, making it the inverse image of a geometric surjection of toposes $\boldsymbol{\mathcal{D}}_{\mathcal{E}}(p^*B) \to \boldsymbol{\mathcal{D}}_{\mathcal{F}}(B)$. Again by the usual equivalence, this inverse image $\boldsymbol{\mathcal{D}}_{\mathcal{F}}(B) \to \boldsymbol{\mathcal{D}}_{\mathcal{E}}(p^*B)$ is therefore comonadic. It follows that the induced left adjoint $\mathbf{Tp}_{\mathcal{F}}(B) \to \mathbf{Tp}_{\mathcal{E}}(p^*B)$ is likewise comonadic. In the end, this is sufficient to establish (2). □

## 6.2 Natural Models of Presheaves and Sheaves

We now introduce natural models of presheaves and natural models of sheaves. We will see that natural models of sheaves do not in general admit a type classifier.

In the rest of this chapter, let $\mathcal{U}$ and $\mathcal{S}$ be Grothendieck universes such that $\mathcal{U} \subsetneq \mathcal{S}$ and $\mathcal{U} \in \mathcal{S}$. We will also denote by $\mathcal{U}$ (resp. $\mathcal{S}$) the category of $\mathcal{U}$-small (resp. $\mathcal{S}$-small) sets and functions. We done by $\mathbf{Cat}_{<\mathcal{U}}$ (resp. $\mathbf{Cat}_{<\mathcal{S}}$) the 2-category of $\mathcal{U}$-small (resp. $\mathcal{S}$-small) categories, functors, and natural transformations.

### 6.2.1 Presheaves

When $\mathbf{C}$ is a $\mathcal{U}$-small category, we denote by $\mathcal{S}^{\mathbf{C}^{\mathrm{op}}}$ the small natural model with eponymous underlying category and small natural model structure given by

$$\mathbf{Tp}(P) = \mathcal{U}^{(\int P)^{\mathrm{op}}}$$

and

$$\mathbf{Tm}(P) = 1/\mathcal{U}^{(\int P)^{\mathrm{op}}} \quad ,$$

where $P \in \mathcal{S}^{\mathbf{C}^{\mathrm{op}}}$, with the evident reindexing and projection operations.

As one might hope, this $\mathcal{S}^{\mathbf{C}^{\mathrm{op}}}$ is a small natural display topos. In this connection, we note only that the type classifier $\mathfrak{U}$ for $\mathcal{S}^{\mathbf{C}^{\mathrm{op}}}$ may be described as the strict $\mathbf{C}$-indexed category given on objects $I \in \mathbf{C}$ by

$$\mathfrak{U}(I) = \mathcal{U}^{(\mathbf{C}/I)^{\mathrm{op}}} \quad ,$$

with the evident action on morphisms (see Hofmann and Streicher, 1997).

Naturally, the foregoing is greatly simplified when the indexing category is a set. When $X$ is a $\mathcal{U}$-small set, we will write $\mathcal{S}^X$ in place of $\mathcal{S}^{X^{\mathrm{op}}}$. The type classifier of $\mathcal{S}^X$ may alternatively be described as the (necessarily strict) $X$-indexed category constant at $\mathcal{U}$.

### 6.2.2 Sheaves

When $(\mathbf{C}, J)$ is an $\mathcal{S}$-small site, we denote by $\mathrm{Sh}_{<\mathcal{S}}(\mathbf{C}, J)$ (resp. $\mathrm{Sh}_{<\mathcal{U}}(\mathbf{C}, J)$) the category of $\mathcal{S}$-small (resp. $\mathcal{U}$-small) sheaves on $(\mathbf{C}, J)$.

When $S \in \mathrm{Sh}_{<\mathcal{S}}(\mathbf{C}, J)$, we denote by $J_S$ the coverage of $\int S$ given by

$$F \in J_S(e) \Leftrightarrow F \text{ is a matching family that amalgamates to } e \quad .$$

When $(\mathbf{C}, J)$ is a $\mathcal{U}$-small site, we also denote by $\mathrm{Sh}_{<\mathcal{S}}(\mathbf{C}, J)$ the small natural model with eponymous underlying category and small natural model structure given by

$$\mathbf{Tp}(S) = \mathrm{Sh}_{<\mathcal{U}}(\int S, J_S)$$

and

$$\mathbf{Tm}(S) = 1/\,\mathrm{Sh}_{<\mathcal{U}}(\int S, J_S) \quad ,$$

where $S \in \mathrm{Sh}_{<\mathcal{S}}(\mathbf{C}, J)$, with the evident reindexing and projection operations.

When $X$ is a $\mathcal{U}$-small topological space, we will write $\mathrm{Sh}_{<\mathcal{S}}(X)$ for $\mathrm{Sh}_{<\mathcal{S}}(\mathcal{O}(X), J)$, in which $J$ denotes the usual coverage of $\mathcal{O}(X)$.

Although the underlying category of the natural model $\mathrm{Sh}_{<\mathcal{S}}(\mathbf{C}, J)$ is a topos, it will not ordinarily admit a type classifier, and thus, unfortunately, will not ordinarily be a small natural display topos. This renders in our terms a well-known issue; the counterexample in the case of $\mathrm{Sh}_{<\mathcal{S}}(X)$ for a $\mathcal{U}$-small topological space $X$ was essentially already noted by Grothendieck (1960, §3.3).

To see the issue, we can try deriving what a type classifier $\mho$ for $\mathrm{Sh}_{<\mathcal{S}}(\mathbf{C}, J)$ must be if it exists. In this case, we have

$$\mho \in \mathbf{Cat}(\mathrm{Sh}_{<\mathcal{S}}(\mathbf{C}, J)) \hookrightarrow \mathbf{Cat}(\mathcal{S}^{\mathbf{C}^{\mathrm{op}}}) \cong (\mathbf{Cat}_{<\mathcal{S}})^{\mathbf{C}^{\mathrm{op}}} \quad ,$$

so we could describe $\mho$ by its corresponding strict $\mathbf{C}$-indexed category, which we abusively also denote by $\mho$.

We must have, on objects $I \in \mathbf{C}$,

$$\begin{aligned}
\mho(I) &\cong (\mathbf{Cat}_{<\mathcal{S}})^{\mathbf{C}^{\mathrm{op}}}(\mathrm{y}\,I, \mho) && ((\mathbf{Cat}_{<\mathcal{S}})\text{-enriched Yoneda Lemma}) \\
&\cong \mathbf{Cat}(\mathcal{S}^{\mathbf{C}^{\mathrm{op}}})(\mathrm{y}\,I, \mho) && (\mathbf{Cat}(\mathcal{S}^{\mathbf{C}^{\mathrm{op}}}) \cong (\mathbf{Cat}_{<\mathcal{S}})^{\mathbf{C}^{\mathrm{op}}}) \\
&= \mathbf{Cat}(\mathcal{S}^{\mathbf{C}^{\mathrm{op}}})(\mathrm{y}\,I, i\mho) && (i : \mathrm{Sh}_{<\mathcal{S}}(\mathbf{C}, J) \hookrightarrow \mathcal{S}^{\mathbf{C}^{\mathrm{op}}}) \\
&\cong \mathbf{Cat}(\mathrm{Sh}_{<\mathcal{S}}(\mathbf{C}, J))(a\,\mathrm{y}\,I, \mho) && (a \dashv i : \mathrm{Sh}_{<\mathcal{S}}(\mathbf{C}, J) \to \mathcal{S}^{\mathbf{C}^{\mathrm{op}}}) \\
&\cong \mathbf{Tp}(a\,\mathrm{y}\,I) \\
&\equiv \mathrm{Sh}_{<\mathcal{U}}(\textstyle\int a\,\mathrm{y}\,I, J_{a\,\mathrm{y}\,I}) \quad ,
\end{aligned}$$

together with the evident action on arrows.

However, such an $\mho$ fails to satisfy the appropriate sheaf condition in general, contradicting the assumption that it lies in $\mathbf{Cat}(\mathrm{Sh}_{<\mathcal{S}}(\mathbf{C}, J))$! It may only be a stack, in which amalgamations are unique up to isomorphism.

This may already be observed in the case where $\mathrm{Sh}_{<\mathcal{S}}(\mathbf{C}, J)$ is $\mathrm{Sh}_{<\mathcal{S}}(X)$, for a $\mathcal{U}$-small topological space $X$. Here, simplifying, we have, for each $V \in \mathcal{O}(X)$,

$$\begin{aligned}
\mho(V) &\cong \mathrm{Sh}_{<\mathcal{U}}(\textstyle\int a\,\mathrm{y}\,V, J_{a\,\mathrm{y}\,V}) \\
&\cong \mathrm{Sh}_{<\mathcal{U}}(\textstyle\int \mathrm{y}\,V, J_{\mathrm{y}\,V}) && \text{(subcanonicality)} \\
&\dots \\
&\cong \mathrm{Sh}_{<\mathcal{U}}(V) \quad ,
\end{aligned}$$

together with the evident restriction mappings. The reader may verify that such an $\mho$ violates the sheaf condition in general, for instance when $X$ is the discrete space with two elements.

## 6.3 Hofmann-Streicher Natural Display Toposes

In this section, we generalize from natural models of presheaves to a wider class of natural display toposes, which we term the **Hofmann-Streicher** natural display toposes. We will see that any natural model of sheaves on a site with **enough points** can be replaced by an equivalent Hofmann-Streicher natural display topos (which, in particular, admits a type classifier). This covers many cases of interest; for example, any natural model sheaves on a topological space has enough points.

**Definition 6.3.1.** *When $(\mathbf{C}, J)$ is a $\mathcal{U}$-small site, we will say that the small, Cartesian natural model $\mathrm{Sh}_{<\mathcal{S}}(\mathbf{C}, J)$ has* ***enough points*** *if there exists some $\mathcal{U}$-small set $X$ and surjective geometric morphism of small, Cartesian natural models $p : \mathcal{S}^X \twoheadrightarrow \mathrm{Sh}_{<\mathcal{S}}(\mathbf{C}, J)$.*

Alternatively, by Proposition 6.1.3, $\mathrm{Sh}_{<\mathcal{S}}(\mathbf{C}, J)$ has enough points when there exists some $\mathcal{U}$-small set $X$ and natural Cartesian display functor $P : \mathrm{Sh}_{<\mathcal{S}}(\mathbf{C}, J) \to \mathcal{S}^X$ such that $P$ is comonadic.

**Definition 6.3.2.** *We will say that a small natural display topos $\mathcal{E}$ is* ***Hofmann-Streicher*** *if it is equivalent to $\mathrm{Sh}_{<\mathcal{S}}(\mathbf{C}, J)$ for some $\mathcal{U}$-small site $(\mathbf{C}, J)$ and there exists some $\mathcal{U}$-small set $X$ and natural Cartesian display functor $P : \mathcal{E} \to \mathcal{S}^X$ such that $P$ is strictly comonadic.*

**Example 6.3.3.** *Let $\mathbf{C}$ be a $\mathcal{U}$-small category. Then, the inclusion functor*

$$|\mathbf{C}| \hookrightarrow \mathbf{C}$$

*induces a surjective geometric morphism of small natural display toposes*

$$\mathcal{S}^{|\mathbf{C}|} \twoheadrightarrow \mathcal{S}^{\mathbf{C}^{\mathrm{op}}} \quad ,$$

*such that the inverse image exhibits $\mathcal{S}^{\mathbf{C}^{\mathrm{op}}}$ as Hofmann-Streicher.*

**Proposition 6.3.4.** *When $(\mathbf{C}, J)$ is a $\mathcal{U}$-small site, the small, Cartesian natural model $\mathrm{Sh}_{<\mathcal{S}}(\mathbf{C}, J)$ has enough points if and only if it is equivalent to a Hofmann-Streicher small natural display topos.*

*Proof.* If $\mathrm{Sh}_{<\mathcal{S}}(\mathbf{C}, J)$ has enough points, choose a surjective geometric morphism of small, Cartesian natural models $p : \mathcal{S}^X \twoheadrightarrow \mathrm{Sh}_{<\mathcal{S}}(\mathbf{C}, J)$, where $X$ is $\mathcal{U}$-small. Then, $\mathrm{Sh}_{<\mathcal{S}}(\mathbf{C}, J)$ is equivalent to the small natural display topos of coalgebras $(\mathcal{S}^X)^{(p^* p_*)}$. □

**Example 6.3.5.** *Let $X$ be a $\mathcal{U}$-small topological space. Then, the continuous inclusion function*

$$|X| \hookrightarrow X$$

*induces a surjective geometric morphism of small, Cartesian natural models*

$$\mathrm{Sh}_{<\mathcal{S}}(|X|) \twoheadrightarrow \mathrm{Sh}_{<\mathcal{S}}(X) \quad .$$

*Since $\mathcal{S}^{|X|} \simeq \mathrm{Sh}_{<\mathcal{S}}(|X|)$, we obtain, by composition, a surjective geometric morphism of small, Cartesian natural models*

$$p : \mathcal{S}^{|X|} \twoheadrightarrow \mathrm{Sh}_{<\mathcal{S}}(X) \quad ,$$

*which exhibits $\mathrm{Sh}_{<\mathcal{S}}(X)$ as having enough points. Thus, by Proposition 6.3.4, we can replace $\mathrm{Sh}_{<\mathcal{S}}(X)$ by the equivalent Hofmann-Streicher small natural display topos $(\mathcal{S}^{|X|})^{(p^* p_*)}$.*

# Chapter 7

# Application: Semantics of S4 Dependent Type Theory

In this brief chapter, we sketch an S4 DTT and its canonical interpretation with respect to a comonad of natural models. This approach is apparently the first modeling S4 DTT without caveats such as idempotence or restricted dependency (as in de Paiva and Ritter, 2016).

**Outline** In Section 7.1, we delineate our S4 DTT. In Section 7.2, we delineate the interpretation of our S4 DTT in a comonad of natural models.

## 7.1 S4 Dependent Type Theory

In this section, we delineate an S4 DTT, **S4DTT**.

The system **S4DTT** extracts the S4 operator from the more complex modal DTT of Shulman (2018). Rules for a dependent S4 operator were first given in Nanevski et al. (2008).

In **S4DTT**, we distinguish between **modal variable assumptions**, which are typed with a double colon, *i.e.* written as

$$u :: A \quad ,$$

and **ordinary variable assumptions**, written as

$$x : A \quad .$$

We follow Nanevski et al. (2008) in distinguishing between **modal context judgements**, which we write as

$$\Delta \vdash \quad ,$$

in which $\Delta$ is a list of modal variables, and **ordinary context judgements**, which we write as

$$\Delta \mid \Gamma \vdash \quad ,$$

in which $\Delta$ is a list of modal variables and $\Gamma$ is a list of ordinary variables.

The hypothetical judgements for types and terms take the form

$$\Delta \mid \Gamma \vdash A$$

and

$$\Delta \mid \Gamma \vdash t : A \quad ,$$

respectively, in which

$$\Delta \mid \Gamma \vdash$$

forms an ordinary context judgement.

We will say that a type or a term is **modal** if it appears in a hypothetical judgement of the form

$$\Delta \mid \cdot \vdash A$$

and

$$\Delta \mid \cdot \vdash t : A \quad ,$$

respectively, *i.e.* with no ordinary variables appearing in its context.

It will be a consequence of the rules of **S4DTT** that modal terms may be substituted for variables of like type that appear as modal assumptions, and, naturally, that not-necessarily-modal terms may be substituted for variables of like type that appear as ordinary assumptions.

Table 7.1 describes the context and variable rules for **S4DTT**. **Modal** rules are marked with a $\Box$, while the **ordinary** rules are unmarked.

$$\frac{}{\cdot \vdash}\text{ Emp.}^{\Box} \qquad \frac{\Delta \mid \cdot \vdash B}{\Delta, u :: B \vdash}\text{ Ext.}^{\Box} \qquad \frac{\Delta, u :: A, \Delta' \mid \Gamma \vdash}{\Delta, u :: A, \Delta' \mid \Gamma \vdash u : A}\text{ Var.}^{\Box}$$

$$\frac{\Delta \vdash}{\Delta \mid \cdot \vdash}\text{ Emp.} \qquad \frac{\Delta \mid \Gamma \vdash B}{\Delta \mid \Gamma, x : B \vdash}\text{ Ext.} \qquad \frac{\Delta \mid \Gamma, x : A, \Gamma' \vdash}{\Delta \mid \Gamma, x : A, \Gamma' \vdash x : A}\text{ Var.}$$

Table 7.1: The Context Rules for **S4DTT**

Table 7.2 describes the rules for the S4 box operator of **S4DTT**.

$$\frac{\Delta \mid \cdot \vdash B}{\Delta \mid \Gamma \vdash \Box B}\ \Box\text{-Form.} \qquad \frac{\Delta \mid \cdot \vdash t : B}{\Delta \mid \Gamma \vdash t^{\Box} : \Box B}\ \Box\text{-Intro.}$$

$$\frac{\Delta \mid \Gamma, x : \Box A \vdash B \qquad \Delta \mid \Gamma \vdash s : \Box A \qquad \Delta, u :: A \mid \Gamma \vdash t : B[u^{\Box}/x]}{\Delta \mid \Gamma \vdash (\mathsf{let}\ u^{\Box} := s\ \mathsf{in}\ t) : B[s/x]}\ \Box\text{-Elim.}$$

$$\frac{\Delta \mid \Gamma, x : \Box A \vdash B \qquad \Delta \mid \cdot \vdash s : A \qquad \Delta, u :: A \mid \Gamma \vdash t : B[u^{\Box}/x]}{\Delta \mid \Gamma \vdash (\mathsf{let}\ u^{\Box} := s^{\Box}\ \mathsf{in}\ t) \equiv t[s/u] : B[s^{\Box}/x]}\ \Box\text{-}\beta\text{-Conv.}$$

$$\frac{\Delta \mid \Gamma, x : \Box A \vdash B \qquad \Delta \mid \Gamma \vdash s : \Box A \qquad \Delta \mid \Gamma, x : \Box A \vdash t : B}{\Delta \mid \Gamma \vdash (\mathsf{let}\ u^{\Box} := s\ \mathsf{in}\ t[u^{\Box}/x]) \equiv t[s/x] : B[s/x]}\ \Box\text{-}\eta\text{-Conv.}$$

Table 7.2: The Rules for $\Box$

## 7.2 Interpretation

In this section, we delineate a canonical interpretation of **S4DTT** in a comonad of natural models.

For the purpose of interpreting empty contexts, we will henceforth assume that the category of contexts of a natural model admits a terminal object, and that morphisms of natural models preserve terminal objects in the category of contexts. With the forbearance of the reader, we will tacitly adjust prior results about natural models to this setting.

We use the method of partial interpretation invented by Streicher (1991) and adapted to modal settings by Birkedal et al. (2020) and Riley et al. (2021).

We avail ourselves in Definition 7.2.1 of various notions entailed by a comonad $\square : \mathbf{C} \to \mathbf{C}$ of (small) natural models, including the forgetful-cofree decomposition $U \dashv F : \mathbf{C} \to \mathbf{C}^{\square}$ (derived from Corollary 4.4.1) and (as discussed in Section 4.4) the induced comonads

$$\mathbb{B} : U^{*}\mathbf{Tp}_{\mathbf{C}} \to U^{*}\mathbf{Tp}_{\mathbf{C}}$$

and

$$\mathbb{B} : U^{*}\mathbf{Tm}_{\mathbf{C}} \to U^{*}\mathbf{Tm}_{\mathbf{C}}$$

and induced right adjoint

$$\mathbb{F} : U^{*}\mathbf{Tp}_{\mathbf{C}} \to \mathbf{Tp}_{\mathbf{C}^{\square}}$$

in $\mathbf{Cat}^{(\mathbf{C}^{\square})^{\mathrm{op}}}$.

To simplify the presentation, we omit from Definition 7.2.1 any clauses concerning weakening.

**Definition 7.2.1.** *Let $\square : \mathbf{C} \to \mathbf{C}$ be a comonad of small natural models. A partial interpretation function $[\![-]\!]$ is given by recursion on the raw syntax of **S4DTT** as follows:*

- $[\![\cdot]\!] = 1_{\mathbf{C}^{\square}}$;
- $[\![\Delta, u :: B]\!] = [\![\Delta]\!].\mathbb{F}[\![\Delta \mid \cdot\, ; B]\!]$;
- $[\![\Delta, u :: A \mid \cdot\, ; u : A]\!] = (\varepsilon_{[\![\Delta \mid \cdot\, ; A]\!]}[Up_{\mathbb{F}[\![\Delta \mid \cdot\, ; A]\!]}])(Uv_{\mathbb{F}[\![\Delta \mid \cdot\, ; A]\!]})$;
- $[\![\Delta \mid \cdot\, ]\!] = U[\![\Delta]\!]$;
- $[\![\Delta \mid \Gamma, x : B]\!] = [\![\Delta \mid \Gamma]\!].[\![\Delta \mid \Gamma; B]\!]$;

- $[\![\Delta \mid \Gamma, x : A; x : A]\!] = v_{[\![\Delta \mid \Gamma; A]\!]}$;
- $[\![\Delta \mid \cdot\, ; \Box B]\!] = \mathbb{B}[\![\Delta \mid \cdot\, ; B]\!]$;
- $[\![\Delta \mid \cdot\, ; t^{\Box}]\!] = \mathbb{B}[\![\Delta \mid \cdot\, ; t]\!]$;
- $[\![\Delta \mid \cdot\, ; \textsf{let}\, u^{L} := s\ \textsf{in}\, t]\!] = [\![\Delta, u :: A \mid \cdot\, ; t]\!][\overline{[\![\Delta \mid \cdot\, ; s]\!]}]$.

**Conjecture 7.2.2.** *The partial interpretation function given in Definition 7.2.1 is sound in the following sense:*

- *if* $\Delta \vdash$*, then* $[\![\Delta]\!] \in \mathbf{C}^{\Box}$*;*
- *if* $\Delta \mid \Gamma \vdash$*, then* $[\![\Delta \mid \Gamma]\!] \in \mathbf{C}$*;*
- *if* $\Delta \mid \Gamma \vdash B$*, then* $[\![\Delta \mid \Gamma; B]\!] \in \mathbf{Tp}_{\mathbf{C}}([\![\Delta \mid \Gamma]\!])$*;*
- *if* $\Delta \mid \Gamma \vdash t : B$*, then* $[\![\Delta \mid \Gamma; t]\!] \in \mathbf{Tm}_{\mathbf{C}}([\![\Delta \mid \Gamma]\!])$ *and* $p_{\mathbf{C}} \circ [\![\Delta \mid \Gamma; t]\!] = [\![\Delta \mid \Gamma; B]\!] \in \mathbf{Tp}_{\mathbf{C}}([\![\Delta \mid \Gamma]\!])$*;*
- *if* $\Delta \equiv \Delta' \vdash$*, then* $[\![\Delta]\!] = [\![\Delta']\!] \in \mathbf{C}^{\Box}$*;*
- *if* $(\Delta \mid \Gamma) \equiv (\Delta' \mid \Gamma') \vdash$*, then* $[\![\Delta \mid \Gamma]\!] = [\![\Delta' \mid \Gamma']\!] \in \mathbf{C}$*;*
- *if* $\Delta \mid \Gamma \vdash B \equiv B'$*, then* $[\![\Delta \mid \Gamma; B]\!] = [\![\Delta \mid \Gamma; B']\!] \in \mathbf{Tp}_{\mathbf{C}}([\![\Delta \mid \Gamma]\!])$*;*
- *if* $\Delta \mid \Gamma \vdash t \equiv t' : B$*, then* $[\![\Delta \mid \Gamma; t]\!] = [\![\Delta \mid \Gamma; t']\!] \in \mathbf{Tm}_{\mathbf{C}}([\![\Delta \mid \Gamma]\!])$.

*Proof.* Deferred to future work. □

### Comparison to Other Work

Our approach to the interpretation of comonadic dependent type theory revises Zwanziger (2019b) and parallels Riley et al. (2021), which deals with a special case. An interpretation for an S4 DTT was previously given by de Paiva and Ritter (2016). However, in that system, types are limited to depending on modal variables, *i.e.* they do not depend on ordinary variables.